\documentclass[11pt]{article}
\usepackage[a4paper,
bindingoffset=0.2in,
left=0.5in,
right=0.5in,
top=0.4in,
bottom=0.6in,
footskip=.4in]{geometry}
\usepackage{amsthm,amsmath,amssymb,amsfonts,amscd}
\usepackage{graphicx}
\usepackage{enumerate}
\usepackage[all]{xy}
\usepackage{booktabs}
\usepackage{mathrsfs}
\usepackage{hyperref}
\usepackage[nocompress]{cite}
\usepackage{caption}
\usepackage{subcaption}
\usepackage{verbatim}
\usepackage{xcolor}
\usepackage[normalem]{ulem}
\usepackage{mathtools}
\usepackage{tikz}
\usepackage{multirow}
\usetikzlibrary{arrows}
\usetikzlibrary{spy}

\mathtoolsset{showonlyrefs,showmanualtags}

\newtheorem{theorem}{Theorem}[section]

\newtheorem{prop}[theorem]{Proposition}
\newtheorem{lemma}[theorem]{Lemma}

\theoremstyle{definition}
\newtheorem{defi}[theorem]{Definition}
\newtheorem{ex}[theorem]{Example}

\theoremstyle{remark}
\newtheorem{remark}[theorem]{Remark}

\newcommand{\R}{\mathbb{R}}
\newcommand{\Z}{\mathbb{Z}}
\newcommand{\N}{\mathbb{N}}

\numberwithin{equation}{section}

\newlength{\spyimagewidth}
\newlength{\spyimageheight}

\newcommand{\spybarbara}[2][\linewidth]{%
  \setlength{\spyimagewidth}{#1}
  \setlength{\spyimageheight}{\spyimagewidth}
  \begin{tikzpicture}[spy using outlines={rectangle, magnification=2.0, size=0.45\spyimagewidth, connect spies,color=white}]
    \node[inner sep=0pt, anchor=south west] (image) at (0,0) {\includegraphics[width=\spyimagewidth]{#2}};
    \coordinate (spy point) at (0.74\spyimagewidth , 0.78\spyimagewidth );
    \coordinate (spy node) at (0.225\spyimagewidth, 0.775\spyimageheight);
    \spy on (spy point) in node at (spy node);
  \end{tikzpicture}%
}

\newcommand{\spypeppers}[2][\linewidth]{%
  \setlength{\spyimagewidth}{#1}
  \setlength{\spyimageheight}{\spyimagewidth}
  \begin{tikzpicture}[spy using outlines={rectangle, magnification=2.0, size=0.45\spyimagewidth, connect spies,color=white}]
    \node[inner sep=0pt, anchor=south west] (image) at (0,0) {\includegraphics[width=\spyimagewidth]{#2}};
    \coordinate (spy point) at (0.75\spyimagewidth , 0.65\spyimagewidth );
    \coordinate (spy node) at (0.225\spyimagewidth, 0.775\spyimageheight);
    \spy on (spy point) in node at (spy node);
  \end{tikzpicture}%
}

\newcommand{\spypeppersblur}[2][\linewidth]{%
  \setlength{\spyimagewidth}{#1}
  \setlength{\spyimageheight}{\spyimagewidth}
  \begin{tikzpicture}[spy using outlines={rectangle, magnification=2.0, size=0.45\spyimagewidth, connect spies,color=white}]
    \node[inner sep=0pt, anchor=south west] (image) at (0,0) {\includegraphics[width=\spyimagewidth]{#2}};
    \coordinate (spy point) at (0.75\spyimagewidth , 0.65\spyimagewidth );
    \coordinate (spy node) at (0.225\spyimagewidth, 0.775\spyimageheight);
    \spy on (spy point) in node at (spy node);
  \end{tikzpicture}%
}

\newcommand{\spyboat}[2][\linewidth]{%
  \begingroup
  \setlength{\spyimagewidth}{#1}%
  \setlength{\spyimageheight}{\spyimagewidth}%
  \begin{tikzpicture}[spy using outlines={rectangle, magnification=2.0, size=0.45\spyimagewidth, connect spies, color=white}]
    \node[inner sep=0pt, anchor=south west] (image) at (0,0)
      {\includegraphics[width=\spyimagewidth]{#2}};
    \coordinate (spy point) at (0.75\spyimagewidth, 0.65\spyimageheight);
    \coordinate (spy node) at (0.225\spyimagewidth, 0.225\spyimageheight);
    \spy on (spy point) in node at (spy node);
  \end{tikzpicture}%
  \endgroup
}

\newcommand{\spycactus}[2][\linewidth]{%
  \setlength{\spyimagewidth}{#1}
  \setlength{\spyimageheight}{\spyimagewidth}
  \begin{tikzpicture}[spy using outlines={rectangle, magnification=2.0, size=0.45\spyimagewidth, connect spies,color=white}]
    \node[inner sep=0pt, anchor=south west] (image) at (0,0) {\includegraphics[width=\spyimagewidth]{#2}};
    \coordinate (spy point) at (0.15\spyimagewidth , 0.8\spyimagewidth );
    \coordinate (spy node) at (0.775\spyimagewidth, 0.805\spyimageheight);
    \spy on (spy point) in node at (spy node);
  \end{tikzpicture}%
}

\newcommand{\spycactusblur}[2][\linewidth]{%
  \setlength{\spyimagewidth}{#1}
  \setlength{\spyimageheight}{\spyimagewidth}
  \begin{tikzpicture}[spy using outlines={rectangle, magnification=2.0, size=0.45\spyimagewidth, connect spies,color=white}]
    \node[inner sep=0pt, anchor=south west] (image) at (0,0) {\includegraphics[width=\spyimagewidth]{#2}};
    \coordinate (spy point) at (0.525\spyimagewidth , 0.375\spyimagewidth );
    \coordinate (spy node) at (0.775\spyimagewidth, 0.805\spyimageheight);
    \spy on (spy point) in node at (spy node);
  \end{tikzpicture}%
}

\newcommand{\spyTGV}[2][\linewidth]{%
  \setlength{\spyimagewidth}{#1}
  \setlength{\spyimageheight}{\spyimagewidth}
  \begin{tikzpicture}[spy using outlines={rectangle, magnification=3.5, size=0.45\spyimagewidth, connect spies,color=white}]
    \node[inner sep=0pt, anchor=south west] (image) at (0,0) {\includegraphics[width=\spyimagewidth]{#2}};
    \coordinate (spy point) at (0.275\spyimagewidth , 0.725\spyimagewidth );
    \coordinate (spy node) at (0.775\spyimagewidth, 0.225\spyimageheight);
    \spy on (spy point) in node at (spy node);
  \end{tikzpicture}%
}

\title{A Weighted Second-Order Nonlocal Total Variation Regularizer for Image Restoration}

\author{Francesc Alcover\thanks{DMI \& IAC3, University of the Balearic Islands, Cra. de Valldemossa, km. 7.5, E-07122 Palma, Illes Balears, Spain. E-mail: \texttt{francesc.alcover@uib.cat}, \texttt{joan.duran@uib.es}, \texttt{catalina.sbert@uib.es}} \and Joan Duran\footnotemark[1] \and Catalina Sbert\footnotemark[1]}

\date{\today}

\begin{document}

\maketitle

\begin{abstract}
Despite the remarkable progress of deep learning in image processing, variational methods remain a fundamental and actively developing framework, offering interpretability, mathematical rigor, and well-founded guarantees. In this work, we introduce a weighted second-order nonlocal total variation regularizer for image restoration, which is based on measuring gradient differences across the image domain and allows for weights that are neither smooth nor symmetric and may vanish outside a search window. We first provide a rigorous definition of the regularizer and prove that it enjoys several desirable properties under realistic assumptions on the underlying weight function. We then consider a nonlocal variant of the variational image restoration model of Chambolle and Lions, and prove the existence and uniqueness of solutions to the associated minimization problem. On the computational side, we develop a discrete formulation of the model and derive the corresponding Chambolle-Pock primal-dual algorithm. Finally, we conduct an extensive experimental evaluation, comparing the proposed regularizer with a range of local and nonlocal, first- and second-order regularizers across several image restoration tasks. The results demonstrate its strong noise-suppression and texture-preservation capabilities, as well as its robustness to nonlocal weights computed under non-ideal conditions.
\end{abstract}



\section{Introduction}
While deep learning methods currently dominate many image processing tasks due to their powerful representation capabilities, model-based methodologies remain an active line of research \cite{zaabouli2025novel, liu2026image, yang2026image, liu2024non, pollock20253dmpr, cox2025model, chen2024fractional}. In particular, variational models provide a mathematically grounded framework for incorporating prior knowledge through regularizer functionals, offering interpretability, stability, and well-understood analytical properties. The connection between model-based and learning-based approaches has been further reinforced by the emergence of deep unfolding techniques \cite{xu2023deep}, which reinterpret iterative optimization algorithms as trainable neural network architectures. This perspective highlights the importance of well-designed image regularizers, even within data-driven settings.

The performance of a variational model depends considerably on the choice of the aforementioned regularizer functional, and the literature offers a rich catalogue of them. The earliest proposals, going back to Tikhonov \cite{tikhonov1963solution,tikhonov1977solutions}, penalize the $L^2$ norm of the gradient. They lead to well-posed problems which are simple to solve, but their solutions are necessarily smooth and therefore blur edges and other discontinuities. Replacing this penalization by the $L^1$ norm of the gradient, as proposed by Rudin, Osher and Fatemi \cite{Rudin-Osher-Fatemi} and further analyzed by Chambolle and Lions \cite{chambolle1997image}, yields the total variation, which preserves discontinuities and has become the standard first-order regularizer. Its well-known disadvantage is the staircase effect, by which regions presenting smooth variations are reconstructed as piecewise constant patterns. Several strategies have been proposed to mitigate it: replacing the $L^1$ penalization by a function with linear growth at infinity and quadratic behavior near the origin \cite{aubert1997variational,aubert2006mathematical}, introducing a spatially varying weight which relaxes the penalization where edges are likely to occur \cite{strong1997spatially,ChenRao2003,davoli2024}, or moving to higher-order variants whose minimizers are piecewise polynomial rather than piecewise constant \cite{Bergounioux,bredies2010total,Hu,kim2020nonlocal}. 

A rather different line of work exploits the self-similarity of natural images. Following the nonlocal filters of Yaroslavsky \cite{Yaroslavsky,YaroslavskyEden} and the nonlocal means algorithm of Buades, Coll and Morel \cite{Buades-Coll-Morel,Buades-Coll-Morel-bis}, variational counterparts were proposed by Kindermann et~al. \cite{Kindermann} and by Gilboa and Osher \cite{Gilboa-Osher,Gilboa-Osher-bis}, in which the penalization no longer acts on the derivatives at a point but on pairs of points, weighted by the similarity of their surrounding patches. Such regularizers recover textures and repetitive structures that local models cannot, at the price of a higher computational cost and of a weight function which, in practice, must be computed from the degraded observation itself. Most of them are nonlocal versions of first-order regularizers, although nonlocal counterparts of higher-order ones have also been proposed, such as the nonlocal total generalized variation of Ranftl et~al. \cite{ranftl2014non} and the unweighted second-order functional introduced in \cite{Kindermann}, which penalizes differences of gradients between distant points.

In our previous work \cite{alcover2026nonlocal}, we introduced a weighted second-order nonlocal regularizer based on the classical total variation, extending the work of Kindermann et al. \cite{Kindermann}, which considered the constant weight case and performed a preliminary theoretical study. In this paper, we give the following relaxed definition of the regularizer: Let $\Omega\subset \R^N$ be a nonempty, connected, open set with finite measure. We say that $\omega:\Omega\times\Omega \to \mathbb{R}$ is an \emph{admissible weight function} if $\omega$ is nonnegative and continuous a.e.~in $\Omega\times \Omega$. Given a fixed admissible weight function $\omega$, we define the \emph{weighted nonlocal total variation} of a function $u \in L^1(\Omega)$ as
\begin{equation}
    \mathrm{NLTV}_\omega(u) := \sup\left\{ \int_\Omega u \; \mathrm{div}_{\mathrm{NL}}(\phi) \mathrm{d}x : \phi \in \mathcal{C}_c^1(\Omega\times\Omega;\R^N), |\phi| \le  \omega \; \mbox{ a.e.~in } \Omega\times\Omega \right\},
\end{equation}
where $\mathrm{div}_{\mathrm{NL}}$ denotes the \emph{nonlocal divergence}, defined for a function $\phi \in {\mathcal{C}_c^1(\Omega\times\Omega;\R^N)}$ as
\begin{equation}\label{intro-div}
    \mathrm{div}_{\mathrm{NL}}(\phi)(x) := \int_\Omega \mathrm{div}_x \left[ \phi(y,x) - \phi(x,y) \right] \mathrm{d}y
\end{equation}
for $x \in \Omega$. Subsequently, we define the {\it space of functions with weighted nonlocal bounded variation} as
\begin{equation}
\mathrm{NLBV}_\omega(\Omega):= \{ u \in L^1(\Omega):\mathrm{NLTV}_\omega(u)<+\infty\}.
\end{equation}

Note that the weight enters this new definition through the constraint on the test fields rather than through the nonlocal divergence, as it does in \cite{alcover2026nonlocal} (see Lemma \ref{OldNLTVDefEquivalence}), so that $\omega$ is never differentiated. This allows us to work with weights which are merely continuous almost everywhere, and in particular with the search windows used in practice \eqref{NLMeansDefCont}, whose indicator function is discontinuous across the boundary of the window.

In \cite{alcover2026nonlocal}, under the additional assumptions that $\omega$ is $\mathcal{C}^1$ and symmetric and using a stronger definition, we conducted a detailed theoretical study of the $\mathrm{NLBV}_\omega$ space and its $\mathrm{NLTV}_\omega$ semi-norm. In particular, we proved the lower semi-continuity of $\mathrm{NLTV}_\omega$ with respect to the weak $L^1(\Omega)$ convergence, its invariance under the addition of affine functions, and the fact that $\mathrm{NLBV}_\omega(\Omega)$ is a Banach space when endowed with an appropriate norm. Moreover, we proved that for regular enough $u$ (e.g., belonging to $W^{1,1}(\Omega)$),
\begin{equation} \label{NLTVFormulaIntro}
    \mathrm{NLTV}_\omega(u) = \int_\Omega\int_{\Omega} \omega(x,y)|\nabla u(y) - \nabla u(x)|\mathrm{d}x \mathrm{d}y.
\end{equation}

In the present work, we show that all aforementioned properties hold for the weaker definition of $\mathrm{NLTV}_\omega$ and the relaxed conditions on $\omega$. In particular, we establish formula \eqref{NLTVFormulaIntro}, provided $u$ has weak first-order derivatives. In an image processing context, this expression shows that our regularizer enforces gradient similarity according to the interactions dictated by $\omega$. For sufficiently regular functions, $\mathrm{NLTV}_\omega$ can thus be interpreted as a second-order regularizer: rather than directly penalizing the magnitude of first-order derivatives, it penalizes variations of the gradient through nonlocal interactions. It is therefore a regularizer which, on the one hand, can effectively exploit the self-similarities present in natural images and, on the other hand, favors locally affine rather than locally constant reconstructions, so that it does not suffer from the staircase effect (see Figure \ref{fig:piecewiseaff}).

After establishing some preliminary and structural results involving $\mathrm{NLBV}_\omega$ and its semi-norm, we consider a nonlocal version of the classical Chambolle-Lions image restoration model \cite{chambolle1997image} incorporating the $\mathrm{NLTV}_\omega$ regularizer. This model aims to recover an image $u$ from a degraded observation of the form
\begin{equation}
    u_0 = Au+n,
\end{equation}
where $A$ is a linear operator and $n$ denotes random noise. In the original paper \cite{chambolle1997image}, under suitable conditions on $u_0$, $A$, the dimension $N$ and the exponent of the underlying $L^p$ space, the existence and uniqueness of solutions to the following problem were established:
\begin{equation} \label{OgChamLionVar}
    \min_{u \in \mathrm{BV}(\Omega)} \mathrm{TV}(u) + \frac{\lambda}{2} \int_\Omega |Au - u_0|^2 \mathrm{d}x,
\end{equation}
where $\mathrm{BV}(\Omega)$ is the space of functions with bounded variation and $\mathrm{TV}$ denotes the total variation.

Our aim is to derive a variational model analogous to \eqref{OgChamLionVar} in the $\mathrm{NLBV}_\omega(\Omega)$ setting. In particular, for a fixed admissible weight function $\omega$, $p=2$ if $N=1$ or $p=N/(N-1)$ if $N\ge 2$, $\lambda >0$, $A: L^p(\Omega) \to L^p(\Omega)$, and $u_0 \in L^2(\Omega)$, we consider the problem

\begin{equation} \label{NLChamLionVarIntro}
    \min_{u \in L^p(\Omega) \cap \mathrm{NLBV}_\omega(\Omega)} \mathrm{NLTV}_\omega(u) + \frac{\lambda}{2} \int_\Omega |Au - u_0|^2 \mathrm{d}x.
\end{equation}
We prove that if $\Omega$ is nonempty, bounded, Lipschitz and connected, $A$ is a continuous linear operator such that
\begin{equation} \label{H2'Intro}
    \forall a\in \R^N, b \in \R: A(a\cdot x + b) = 0 \; \text{ a.e.~in } \Omega \iff a=0, b=0,
\end{equation}
and $\omega$ satisfies the following boundedness condition:
\begin{equation} \label{H0Intro}
    \begin{array}{c}
        \mbox{There exists a finite open and connected cover of } \Omega, ( \Omega_n )_{n =1}^M,\\\mbox{ and a constant } c>0 \mbox{ such that } 
     \omega \ge c \; \mbox{ a.e.~in } \Omega_n\times\Omega_n \mbox{ for all } n=1,..., M,
    \end{array}
\end{equation}
then \eqref{NLChamLionVarIntro} admits a solution. Moreover, we show that if $A$ is injective, then the solution is unique. The precise result is Theorem \ref{ExistUniqNLChamLionVar}.

Condition \eqref{H2'Intro} is a modification of the original hypothesis appearing in \cite{chambolle1997image}, and is needed to prove the coercivity of the functional in \eqref{NLChamLionVarIntro}. The modification naturally stems from the fact that $\mathrm{NLTV}_\omega$ annihilates affine functions, instead of constant ones, as classical $\mathrm{TV}$ does.

Moreover, condition \eqref{H0Intro} on $\omega$ implies that $\mathrm{NLBV}_\omega(\Omega) \hookrightarrow \mathrm{BV}(\Omega)$, which allows for the Sobolev embedding and compactness to hold in the nonlocal setting (see Lemma \ref{LemmaSobEmbwboundedbelow}). Originally, in \cite{alcover2026nonlocal} we showed that if $\omega$ was bounded from below away from zero then the embedding into $\mathrm{BV}$ held. However, such a condition is unrealistic in most applications, hence why we introduce \eqref{H0Intro}, which is satisfied by weights which make use of search windows and thus vanish on a significant part of the domain (see Remark \ref{SearchWindowSatisfiesH0}). 

Theorem \ref{ExistUniqNLChamLionVar} can be viewed as a generalization of \cite[Theorem 5.5]{Kindermann} by Kindermann et~al., which establishes the existence and uniqueness of solutions to a similar problem in the case $\omega \equiv 1$. In fact, condition \eqref{H2'Intro} is also required for such a result. The argument in the aforementioned paper relies on an approximation result, namely \cite[Lemma 5.1]{Kindermann}, whose proof is stated to be analogous to that of the classical $\mathrm{BV}$ setting. The density of smooth functions in $\mathrm{NLBV}_\omega$ is, however, delicate, as the usual argument rests on a product rule which is not available for the nonlocal gradient and the nonlocal divergence. Our proof circumvents this issue: it only uses lower semi-continuity and compactness, and therefore requires no approximation-type result.

Having established the existence and uniqueness of solutions to \eqref{NLChamLionVarIntro}, we consider several image restoration tasks, all formulated within the framework of \eqref{NLChamLionVarIntro} through appropriate choices of the operator $A$. We demonstrate experimentally that the affinity of $\mathrm{NLTV}_\omega$ for affine reconstructions, together with its ability to adapt regularization to the nonlocal weights, enables remarkable structural recovery across all the image restoration tasks considered. Furthermore, we show that enforcing similarity at the gradient level, rather than at the pixel level, as in \cite{Gilboa-Osher,Gilboa-Osher-bis,ranftl2014non}, mitigates the adverse effects of unreliable weights computed from heavily corrupted initial images. This approach therefore yields better overall performance without requiring any preprocessing of the initial data.

In summary, the main contributions of this work are the following:
\begin{itemize}
\item We formally introduce a weighted second-order nonlocal total variation ($\mathrm{NLTV}_\omega$) regularizer and its corresponding function space, expanding the work of Kindermann et~al. \cite{Kindermann} and our previous theoretical study \cite{alcover2026nonlocal}. We also introduce condition \eqref{H0Intro} on the weights, which is compatible with the use of search windows and ensures several desirable properties of the regularizer and its associated space. In particular, we prove that the resulting nonlocal space is continuously embedded into the classical $\mathrm{BV}$ space, and that the regularizer vanishes only on affine functions. Condition \eqref{H0Intro} replaces the more restrictive assumption that the weights be uniformly bounded away from zero, and we conjecture that this condition may be of use in other theoretical image processing contexts.

\item Building on the preceding theoretical analysis, under certain conditions, we prove existence and uniqueness of solutions to a general image restoration variational model incorporating the nonlocal regularizer, based on the model presented by Chambolle and Lions \cite{chambolle1997image}. These conditions are compatible with a wide range of image restoration tasks. 

\item We formulate a discrete version of the $\mathrm{NLTV}_\omega$-based Chambolle-Lions model and derive the corresponding primal-dual algorithm \cite{chambolle2011first} to compute numerical solutions. Through experiments on several image restoration tasks, we demonstrate that our model compares favorably with existing image regularizers in terms of both quantitative metrics and visual quality. In particular, $\mathrm{NLTV}_\omega$ achieves remarkable structural recovery and performs well even when the nonlocal weights are computed directly from corrupted images, where heavy degradation may render these weights unreliable.
\end{itemize}

The rest of the paper is organized as follows. Section \ref{SOTASec} provides a brief historical overview of image regularizers, focusing on those most relevant to the present work. Section \ref{NLTVPriorSec} formally defines the $\mathrm{NLTV}_\omega$ regularizer and its corresponding function space $\mathrm{NLBV}_\omega$, and presents the novel weight condition \eqref{H0Intro}. It also establishes the embedding of $\mathrm{NLBV}_\omega$ into $\mathrm{BV}(\Omega)$  and characterizes affine functions as precisely those annihilated by $\mathrm{NLTV}_\omega$. Section \ref{NLChamLionsSec} introduces the nonlocal Chambolle–Lions model \eqref{NLChamLionVarIntro}, establishes the existence and uniqueness of solutions to the corresponding variational problem, and presents its discrete formulation together with the associated primal-dual algorithm. Section \ref{ExpResSec} presents an extensive experimental comparison of the proposed model against other image regularizers across a range of image restoration tasks consistent with the theoretical framework. The comparison considers computational cost, quantitative performance metrics, and visual reconstruction quality. Finally, Section \ref{ConclusionsSec} summarizes the main findings and outlines directions for future research. Some technical remarks are presented in the Appendix, namely the proofs of the equivalence of the $\mathrm{NLTV}_\omega$ definitions presented in this article and \cite{alcover2026nonlocal} and that of formula \eqref{NLTVFormulaIntro}, both given in Appendix \ref{app:TechProofs}, and the consistency of the nonlocal weights used in the experiments with the theoretical setting, proven in Appendix \ref{app:NLMeansWeights}.  

\textbf{Acknowledgements.} F.A., J.D. and C.S. are supported by the MoMaLIP project PID2021-125711OB-I00 funded by MICIU/AEI/10.13039/501100011033 and the European Union NextGeneration EU/PRTR.

\section{Image Processing Regularizers} \label{SOTASec}

Most tasks in image processing can be formulated as ill-posed inverse problems, where uniqueness, stability and even existence of solutions are not guaranteed. A classical way to address this issue is through regularization, which reformulates the inverse problem as a variational one in which an energy functional is minimized. This energy typically combines a data fidelity term, derived from the image formation model, with a regularization term that encodes prior assumptions about the desired structure of the solution and promotes the notion of a “natural” image. While the fidelity term is usually well-defined by the acquisition process, the main difficulty lies in designing a robust \emph{regularizer} that enforces appropriate smoothness and structural constraints. Moreover, a useful image regularizer should not be tailored to a specific application, but instead be sufficiently flexible to perform well across a range of image processing tasks.

The earliest image regularizers were based on the assumption that natural images are approximately piecewise constant functions and therefore exhibit small gradient magnitudes over most of their domain. Under this assumption, Tikhonov \cite{tikhonov1963solution} proposed a regularization scheme based on minimizing the $L^2$ norm of the image gradient. However, the resulting solutions are necessarily smooth, leading to the undesirable blurring of edges and other discontinuities commonly present in images. To overcome this limitation, Rudin, Osher, and Fatemi \cite{Rudin-Osher-Fatemi} (see also Chambolle and P.-L. Lions \cite{chambolle1997image}) replaced the $L^2$ norm with the $L^1$ norm of the gradient, giving rise to the notion of \emph{total variation}. Formally, for an open and bounded set $\Omega \subset \R^N$, the total variation of a function $u \in L^1(\Omega)$ is defined as

\begin{equation} \label{TVSOTAdef}
    \mathrm{TV}(u) := \sup\left\{ \int_\Omega u \; \mathrm{div}(\varphi) \mathrm{d}x : \varphi \in \mathcal{C}_c^1(\Omega;\R^N), \|\varphi\|_\infty \le 1 \right\}.
\end{equation}
The \emph{space of functions with bounded variation} is consequently defined as
\begin{equation}
    \mathrm{BV}(\Omega) := \left\{ u \in L^1(\Omega) : \mathrm{TV}(u) < +\infty \right\}.
\end{equation}
For a comprehensive treatment of the theory of $\mathrm{BV}$ spaces and total variation, we refer the reader to, for instance, \cite{AmbrosioFuscoPallara}. An integration-by-parts argument shows that, for sufficiently smooth functions, the total variation coincides with the $L^1$ norm of its gradient, that is,

\begin{equation} \label{TVSOTAaux1}
    \mathrm{TV}(u) = \int_\Omega |\nabla u (x)| \mathrm{d}x,
\end{equation}
where $|\cdot|$ denotes the Euclidean norm in $\R^N$. More importantly, total variation remains well defined for discontinuous functions and naturally promotes piecewise constant solutions. As a consequence, in an image processing context, total variation regularization preserves sharp edges while effectively removing noise and mitigating other image degradations. This represented a major advance in the state of the art at the time and sparked a significant rise of other $\mathrm{TV}$-based models \cite{aubert1997variational,strong1997spatially, ChenRao2003, bredies2010total}.

One of the main drawbacks of $\mathrm{TV}$ regularization is the \emph{staircase effect}, in which regions that present smooth pixel variations are erroneously reconstructed as piecewise constant patterns. Many alternatives to total variation have been proposed to combat this flaw. In \cite{aubert1997variational, aubert2006mathematical}, Aubert et~al. proposed to replace \eqref{TVSOTAaux1} by

\begin{equation} \label{TVSOTAAubert}
    \int_\Omega \Phi( |\nabla u|) \mathrm{d}x,
\end{equation}
where $\Phi : \R^+ \to \R^+$ is assumed to be strictly convex, non-decreasing and coercive for the purpose of well-posedness. Additional conditions on $\Phi$ and its derivatives were imposed derived from the corresponding Euler-Lagrange equation so that the regularization favors images composed of homogeneous regions separated by sharp edges. The functions satisfying these properties behave quadratically close to zero and linearly away from it. Consequently, the resulting regularizer behaves like standard $\mathrm{TV}$ regularization in regions with high gradient values, while reducing to Tikhonov regularization in low-gradient areas, which preserves smooth transitions and reduces the staircase effect. While the authors chose $\Phi(t)= \sqrt{1+t^2}$ for their experiments, many other functions can be considered that satisfy the properties, such as Huber functions.

Another straightforward extension of $\mathrm{TV}$, first proposed by Strong et~al. in \cite{strong1997spatially}, consists of incorporating a spatially varying weight in \eqref{TVSOTAaux1}, so that the regularizer becomes
\begin{equation} \label{TVSOTAAux2}
    \int_\Omega \omega(x) |\nabla u(x)| \mathrm{d}x.
\end{equation}
For an appropriate choice of the weight $\omega$, this addition produces a regularization that better adapts to the geometry of the underlying image. For example, in the aforementioned paper by Strong et~al., they proposed a weighting scheme in which $\omega$ measures the likelihood of an edge being present at any given location. More precisely, they chose $\omega(x) = (1+|\nabla f(x)|^2)^{-3/2}$, where $f$ denotes the observed image. This choice assigns smaller weights to regions where edges are likely to occur, thereby reducing the amount of smoothing across discontinuities. Of course, many other choices of $\omega$ are possible, depending on the image processing task and the type of prior information one wishes to incorporate into the regularization, see, e.g., \cite{ChenRao2003, davoli2024}.

The regularizers discussed so far are local, as they only involve infinitesimal neighborhoods around each point. Nonlocal methods, instead, are based on the fact that natural images are self-similar, in the sense that for each small patch there exist many other similar patches, possibly located far away in the image. Therefore, such regularizers tend to better preserve repeating patterns and can better distinguish degradation effects such as noise from genuine image textures.

In the highly influential work \cite{Buades-Coll-Morel,Buades-Coll-Morel-bis}, Buades, Coll and Morel introduced a nonlocal filter for image denoising, building upon previous approaches by Yaroslavsky and Eden \cite{Yaroslavsky,YaroslavskyEden}. Variational interpretations of this filter were later proposed by Kindermann et~al. \cite{Kindermann} as well as by Gilboa and Osher \cite{Gilboa-Osher,Gilboa-Osher-bis}. In particular, in the latter reference, the regularizer \eqref{TVSOTAaux1} is replaced by the nonlocal isotropic functional
\begin{equation} \label{GilboaSOTAIso}
\int_\Omega \sqrt{\int_\Omega \omega(x,y) (u(y)-u(x))^2 \mathrm{d} y} \; \mathrm{d}x,
\end{equation}
and the nonlocal anisotropic functional
\begin{equation}\label{GilboaSOTA}
\int_{\Omega \times \Omega} \sqrt{\omega(x,y)}\lvert u(y)-u(x) \rvert \mathrm{d}y\mathrm{d}x,
\end{equation}
where $\omega$ is a positive weight function that dictates the interactions between pairs of points. As standard $\mathrm{TV}$ regularization does, the functionals above annihilate constants and therefore promote piecewise constant reconstructions. The pieces, however, are no longer spatial regions but groups of pixels linked by the weight function, as the penalization drives pixels with similar surrounding patches towards a common intensity.

As already argued, the staircase effect is a consequence of the fact that the $\mathrm{TV}$ annihilates constants and its induced regularization favors piecewise constant reconstructions. To combat this effect, a common strategy has been to consider higher-order $\mathrm{TV}$ variants which encourage piecewise polynomial functions of degree greater than one as solutions \cite{Kindermann, Bergounioux, bredies2010total,Hu, kim2020nonlocal,ranftl2014non}. Among regularizers derived from this strategy, we highlight the Total Generalized Variation ($\mathrm{TGV}$) regularizer, first introduced in \cite{bredies2010total} by Bredies et~al. Its definition may be interpreted as an extension of total variation induced by symmetrized derivatives of order possibly greater than one. Of particular interest is the second-order case, which, for a sufficiently smooth function $u$, admits the closed form formula
\begin{equation} \label{TGVDef}
    \mathrm{TGV}^2(u) = \min_w \alpha_1 \int_\Omega |\nabla u - w| \mathrm{d}x + \alpha_0 \int_\Omega |E w| \mathrm{d}x,
\end{equation}
where $w : \Omega \to \R^2$ is an auxiliary vector field, $E$ denotes the symmetrized gradient operator and $\alpha_0, \alpha_1 > 0$ are weighting parameters. Importantly, $\mathrm{TGV}^2(u)=0$ if and only if $u$ is affine, from which the regularization induced by this functional promotes piecewise affine functions. Consequently, $\mathrm{TGV}^2$ does not lead to the staircase effect.

A nonlocal counterpart of $\mathrm{TGV}^2$ was proposed by Ranftl et~al.~\cite{ranftl2014non} in the context of optical flow estimation. For sufficiently smooth functions, it reads
\begin{equation} \label{NLTGVDef}
    \mathrm{NLTGV}^2(u) = \min_w \int_\Omega \int_\Omega \omega_1(x,y) |u(x) - u(y) - \langle w(x), x - y \rangle | \mathrm{d}x \mathrm{d}y + \sum_{i=1}^2 \int_\Omega \int_\Omega \omega_0(x,y) |w^i(x) - w^i(y)| \mathrm{d}y \mathrm{d}x,
\end{equation}
where $w: \Omega \to \R^2$ is an auxiliary vector field, $w(x) = (w^1(x),w^2(x))$, and $\omega_0$ and $\omega_1$ are nonnegative weight functions. We can interpret the expression above as the $\mathrm{TGV}^2$ functional but with the total variation replaced by Gilboa and Osher's anisotropic regularizer \eqref{GilboaSOTA}. Moreover, $\omega_0$ and $\omega_1$ are now weight functions dictating nonlocal interactions instead of simply being constant parameters, which enables for a more expressive regularization. $\mathrm{NLTGV}^2$ regularization favors piecewise affine reconstructions while taking nonlocal interactions into account.

Finally, let us discuss the regularizer introduced by Kindermann et~al.~\cite{Kindermann}, which has subsequently been applied to illumination and reflectance decomposition \cite{Wang-Ng} and image segmentation \cite{zhai2011novel}, and which is the main inspiration behind the present work. In \cite{Kindermann}, the authors define the \emph{Nonlocal Total Variation} of a function $u \in L^1(\Omega)$ as
\begin{equation}
    \mathrm{NLTV}(u) = \sup \left\{ \int_\Omega u \; \mathrm{div}_{\mathrm{NL}}(\phi) \mathrm{d}x : \phi \in \mathcal{C}_c^1(\Omega\times\Omega;\R^N), \|\phi\|_\infty \le 1\right\},
\end{equation}
where the nonlocal divergence is defined as in \eqref{intro-div}. For sufficiently smooth functions, the functional admits the formula
\begin{equation}
    \mathrm{NLTV}(u)=\int_{\Omega\times\Omega} \left| \nabla u(y) - \nabla u(x) \right| \mathrm{d}x \mathrm{d}y.
\end{equation}
Unlike \eqref{GilboaSOTAIso}, \eqref{GilboaSOTA}, this regularizer enforces patch similarity at the gradient level rather than at the pixel intensity level. It therefore annihilates affine functions instead of constants, which, as argued above, avoids the staircase effect. For numerical experiments, they propose to restrict nonlocal interactions only to pixels within a given distance.

\section{The Weighted Nonlocal Total Variation Regularizer} \label{NLTVPriorSec}

In this section, we formally define the proposed weighted nonlocal total variation regularizer and study its behavior through theoretical and experimental results. 

\subsection{Formal Definition and Preliminary Results} \label{DefPreliminariesSubsec}

Let $\Omega$ be a nonempty, connected, open set of finite measure. Throughout this manuscript we refer to $\Omega$ as \emph{domain}. Some of the results presented in this section hold under weaker assumptions on $\Omega$; whenever this is the case, we will explicitly state the corresponding conditions. We begin by introducing the basic assumptions on the weight functions.

\begin{defi}[Admissible weight function]
    We say that $\omega:\Omega\times\Omega\to \R$ is an \emph{admissible weight function} if it is nonnegative and continuous a.e.~in $\Omega\times\Omega$.
\end{defi}

Below, we define the weighted nonlocal total variation regularizer formally in the continuous setting.

\begin{defi}[Weighted nonlocal total variation] \label{NLBVDef}
    Let $\omega$ be an admissible weight function. Given $u \in L^1(\Omega)$, we define its \emph{weighted nonlocal total variation (induced by $\omega$)} as

    \begin{equation} \label{NLTV_wDef}
        \mathrm{NLTV}_\omega(u) := \sup \left\{ \int_\Omega u(x) \; \mathrm{div}_\mathrm{NL} (\phi)(x) \mathrm{d}x : \phi \in \mathcal{C}^1_c (\Omega \times \Omega; \mathbb{R}^N), | \phi | \le \omega \; \mbox{ a.e.~in } \Omega\times\Omega \right\},
    \end{equation}
    where $\mathrm{div}_\mathrm{NL}$ denotes the \emph{nonlocal divergence}, defined for a function $\phi \in {\mathcal{C}_c^1(\Omega\times\Omega;\R^N)}$ as
    \begin{equation}
        \mathrm{div}_{\mathrm{NL}}(\phi)(x) := \int_\Omega \mathrm{div}_x \left[ \phi(y,x) - \phi (x,y) \right] \mathrm{d}y
    \end{equation}
    for $x \in \Omega$. Here, $\mathrm{div}_x$ denotes the classical divergence with respect to the $x$ variable.
    
    The \emph{space of functions with weighted nonlocal bounded variation (induced by $\omega$)} is defined as
        
    $$\mathrm{NLBV}_\omega(\Omega) := \{ u \in L^1(\Omega) : \mathrm{NLTV}_\omega(u) < +\infty \}.$$
\end{defi}

As mentioned before, the constant-weight case was previously studied by Kindermann et~al. \cite{Kindermann}, who also proved the existence of solutions to a similar problem to that of Chambolle and Lions \cite{chambolle1997image}. Their argument, however, relies on an approximation result, namely \cite[Lemma 5.1]{Kindermann}, whose proof is stated to be analogous to that in the local setting. We have encountered difficulties in reproducing it, as the nonlocal gradient and divergence operators do not satisfy a product rule, a key ingredient in deriving an estimate used in the classical proof. In this article, we establish existence for general $\omega$ without relying on any approximation-type result. Following the notation introduced by Kindermann et al., we denote by $\mathrm{NLTV}$ and $\mathrm{NLBV}$ the weighted nonlocal total variation and the space of functions with weighted nonlocal bounded variation induced by the constant weight function $\omega \equiv 1$, respectively.

In our previous work \cite{alcover2026nonlocal}, we conducted a detailed theoretical study of the $\mathrm{NLBV}_\omega$ space and the $\mathrm{NLTV}_\omega$ functional under the additional assumptions that $\omega$ is $\mathcal{C}^1$ and symmetric. In particular, we considered a different definition than \eqref{NLTV_wDef} in which we explicitly differentiated $\omega$, allowing for the corresponding function space to be closely related to a nonlocal Sobolev space, also introduced in the article. For $\mathcal{C}^1$ and symmetric weights both definitions of $\mathrm{NLTV}_\omega$ and $\mathrm{NLBV}_\omega(\Omega)$ coincide, as the next lemma states. We detail its proof in Appendix \ref{app:TechProofs} to preserve the flow of the text, as it is rather technical.

\begin{lemma} \label{OldNLTVDefEquivalence}
    Let $\omega$ be an admissible weight function which, in addition, is symmetric and belongs to $\mathcal{C}^1(\Omega\times\Omega)$. Then

    $$\mathrm{NLTV}_\omega(u)= \sup \left\{  \int_\Omega u(x) \; \mathrm{div}_\omega (\phi)(x) \mathrm{d}x : \phi \in \mathcal{C}^1_c (\Omega \times \Omega; \mathbb{R}^N), \| \phi \|_{L^\infty(\Omega\times\Omega;\R^N)} \le 1 \right\}$$
    for any $u \in L^1(\Omega)$, where $\| \cdot \|_{L^\infty(\Omega\times\Omega;\R^N)}$ denotes the uniform norm in $\mathcal{C}_c(\Omega\times\Omega;\R^N)$ and

    \begin{equation}
        \mathrm{div}_{\omega}(\phi)(x) := \int_\Omega \mathrm{div}_x \left[ \omega(x,y) ( \phi(y,x) - \phi (x,y) ) \right] \mathrm{d}y
    \end{equation}
    for $x \in \Omega$.
\end{lemma}

Note that in particular the equivalence above is satisfied by the weight $\omega \equiv 1$ in $\Omega\times\Omega$, hence we can directly apply the results presented in \cite{alcover2026nonlocal} that only involve the $\mathrm{NLBV}$ space and its $\mathrm{NLTV}$ semi-norm.

In this work, we consider continuous weights only in the ``almost everywhere" sense, which are applicable to a broader class of settings. In the following, we establish properties and results relevant to the present article. Many of them are generalizations of the results presented in \cite{alcover2026nonlocal}.

\begin{lemma}
\label{lemma:lsc_nlbv}
    Let $\omega$ be an admissible weight function. If $u\in L^1(\Omega)$ and $(u_n)_{n\in\mathbb{N}}$ is a sequence in $\mathrm{NLBV}_\omega(\Omega)$ such that $u_n \rightharpoonup u$ weakly in $L^1(\Omega)$, then

    \begin{equation}
        \mathrm{NLTV}_\omega(u) \le \liminf_{n \to +\infty} \mathrm{NLTV}_\omega(u_n).
    \end{equation}
\end{lemma}

\begin{proof}
    For any $\phi \in \mathcal{C}_c^1(\Omega\times\Omega;\mathbb{R}^N)$ such that $|\phi|\le \omega$ a.e.~in $\Omega\times\Omega$ one readily checks that $\mathrm{div}_\mathrm{NL}(\phi)$ belongs to $L^\infty(\Omega)$, which by assumption implies that 
    \begin{equation*}
    \lim_{n \to \infty}\int_\Omega u_n \mathrm{div}_\mathrm{NL}(\phi) \mathrm{d}x= \int_\Omega u \; \mathrm{div}_\mathrm{NL}(\phi) \mathrm{d}x.
    \end{equation*}
    
    As a consequence,
    \begin{equation}
        \int_\Omega u \; \mathrm{div}_\mathrm{NL}(\phi) \mathrm{d}x \le \liminf_{n \to +\infty} \mathrm{NLTV}_\omega(u_n),
    \end{equation}
    and by taking supremum over all $\phi$ the result is proven.
\end{proof}

From lower semi-continuity it follows that the space of functions with weighted nonlocal bounded variation is a Banach space for any admissible weight function.

\begin{prop} \label{Proposition_Banach_NLBV}
    Let $\omega$ be an admissible weight function. Then $\mathrm{NLBV}_\omega(\Omega)$ is a Banach space with the norm

    $$\|u\|_{\mathrm{NLBV}_\omega(\Omega)} := \|u\|_{L^1(\Omega)} + \mathrm{NLTV}_\omega(u).$$
\end{prop}

\begin{proof}
    Analogous to the local case. In particular, one can directly adapt the arguments presented in \cite[Theorem 10.1.1]{attouch2014variational}.
\end{proof}

Next, we formally state that $\mathrm{NLTV}_\omega$ satisfies the formula \eqref{NLTVFormulaIntro} for sufficiently smooth functions. The proof presented in \cite[Proposition 3.18]{alcover2026nonlocal} of the same statement does not translate well into our new setting. As it turns out, the new conditions on $\omega$ and the weaker $\mathrm{NLTV}_\omega$ definition introduce several technical nuances that prevent a straightforward adaptation of the original argument. As a consequence, the proof which we have derived is very technical, and we detail it in Appendix \ref{app:TechProofs} to preserve the flow of the main text.

\begin{prop}
\label{C1NLTV}
    Let $\omega$ be an admissible weight function. For any $u \in L^1(\Omega)$ with weak first-order derivatives we have
    \begin{equation} \label{NLTVFormula}
        \mathrm{NLTV}_\omega(u) = \int_{\Omega \times \Omega} \omega(x,y) |\nabla u (y) - \nabla u(x)| \mathrm{d}x \mathrm{d}y,
    \end{equation}
    where $\nabla u$ is the weak gradient of $u$.
\end{prop}

While the proposition above is only needed for consistency with the numerical experiments presented in Section \ref{ExpResSec}, next we prove some results which will be crucial in the theoretical arguments of the following sections. We begin with some embeddings involving the $\mathrm{NLBV}(\Omega)$ space.

\begin{lemma} \label{NLBV_wIntoNLBV}
    Let $\omega$ be an admissible weight function. If there exists $c>0$ such that
    \begin{equation}
        c \le \omega \; \mbox{ a.e.~in } \Omega\times\Omega,
    \end{equation}
    then $\mathrm{NLBV}_\omega(\Omega) \hookrightarrow \mathrm{NLBV}(\Omega)$.
\end{lemma}

\begin{proof}
    Given $\phi \in \mathcal{C}_c^1(\Omega\times\Omega;\mathbb{R}^N)$ such that $|\phi| \le 1$ in $\Omega\times\Omega$, notice that for almost every $(x,y)\in\Omega\times\Omega$

    \begin{equation}
        |c \phi(x,y)|\le c \le \omega(x,y). 
    \end{equation}

    With this in mind,

    \begin{align}
        \int_\Omega u(x) \mathrm{div}_\mathrm{NL}(\phi)(x) \mathrm{d}x &= \frac{1}{c} \int_\Omega u(x) \mathrm{div}_\mathrm{NL}(\phi)(c \phi)(x) \mathrm{d}x \le \frac{1}{c} \mathrm{NLTV}_\omega(u),
    \end{align}
    and by taking supremum over all $\phi$, we obtain that
    \begin{equation} \label{NLTV<NLTV_w}
        \mathrm{NLTV}(u) \le \frac{1}{c} \mathrm{NLTV}_\omega(u).
    \end{equation}
    
    From this, the desired embedding follows.
\end{proof}

\begin{prop} \label{NLBV_into_BV}
    Let $\omega$ be an admissible weight function. If there exists $c>0$ such that
    \begin{equation} \label{w_bounded_below}
        c \le \omega \mbox{ a.e.~in } \Omega\times\Omega,
    \end{equation}
    then $\mathrm{NLBV}_\omega(\Omega) \hookrightarrow \mathrm{BV}(\Omega)$.
\end{prop}

\begin{proof}
    Follows from the fact that $\mathrm{NLBV}(\Omega) = \mathrm{BV}(\Omega)$ with equivalence of norms by \cite[Proposition 4.3]{alcover2026nonlocal} and Lemma \ref{NLBV_wIntoNLBV}.
\end{proof}

The following result establishes one of the statements made in \cite[Remark 4.4]{alcover2026nonlocal}. Although it only involves the constant weight, we still provide a proof, as it is closely related to the nullspace of $\mathrm{NLTV}_\omega$.

\begin{prop} \label{TVaffNLTVProp}
    For any $u \in \mathrm{NLBV}(\Omega)$,
    \begin{equation}
        \inf_{g \,\mathrm{ affine}}\mathrm{TV}(u-g) \le \frac{1}{|\Omega|} \mathrm{NLTV}(u).
    \end{equation}
\end{prop}

\begin{proof}
    Fix $\theta \in \mathcal{C}_c^1(\Omega)$ satisfying $\|\theta\|_{L^\infty(\Omega)} \le 1$ and $\int_{\Omega} \theta \,\mathrm{d}z \neq 0$. Let $\varphi \in \mathcal{C}_c^1(\Omega;\R^N)$ be such that $\|\varphi\|_{L^\infty(\Omega;\R^N)}\le 1$. Then, the function
    \begin{equation}
        \phi_{\theta,\varphi}(x,y) := \theta(x) \varphi(y),
    \end{equation}
    defined for $(x,y) \in \Omega\times\Omega$, clearly belongs to $\mathcal{C}_c^1(\Omega\times\Omega;\R^N)$ and satisfies $\|\phi_{\theta,\varphi}\|_{L^\infty(\Omega\times\Omega;\R^N)} \le 1$. A straightforward computation shows that

    \begin{equation}
        \int_\Omega u(x) \; \mathrm{div}_{\mathrm{NL}}(\phi_{\theta,\varphi})(x) \mathrm{d}x = \left( \int_\Omega \theta(y) \mathrm{d}y \right) \int_\Omega u(x) \; \mathrm{div} (\varphi)(x) \mathrm{d}x - \int_\Omega u(x) \left[ \nabla \theta(x) \cdot \int_\Omega \varphi(y) \mathrm{d}y \right] \mathrm{d}x,
    \end{equation}
    from which

    \begin{align}
        \int_\Omega u(x) \; \mathrm{div} (\varphi)(x) \mathrm{d}x
        &= \frac{1}{\int_\Omega \theta(z) \mathrm{d}z} \left[ \int_\Omega u(x) \; \mathrm{div}_{\mathrm{NL}}(\phi_{\theta,\varphi})(x) \mathrm{d}x  + \left( \int_\Omega u(x) \nabla \theta(x) \mathrm{d}x \right) \cdot \int_\Omega \varphi(y) \mathrm{d}y \right].
    \end{align}

    Consider now the affine function given by

    \begin{equation}
        g_\theta(x) := \left( -\frac{1}{\int_\Omega \theta(z)\mathrm{d}z} \int_\Omega u(y) \nabla \theta(y) \mathrm{d}y\right)\cdot x
    \end{equation}
    for $x\in \Omega$. By the previous computation we have that

    \begin{align}
        \int_\Omega (u(x)-g_\theta(x)) \mathrm{div} (\varphi)(x) \mathrm{d}x &= \int_\Omega u(x) \mathrm{div} (\varphi)(x) \mathrm{d}x - \int_\Omega g_\theta(x) \mathrm{div} (\varphi)(x) \mathrm{d}x\\
        &= \frac{1}{\int_\Omega \theta(z) \mathrm{d}z} \int_\Omega u(x) \; \mathrm{div}_{\mathrm{NL}}(\phi_{\theta,\varphi})(x) \mathrm{d}x\\
        &\le \frac{1}{\left|\int_\Omega \theta(z) \mathrm{d}z\right|} \mathrm{NLTV}(u).
    \end{align}

    Therefore, by taking supremum over all $\varphi$, we have

    \begin{equation}
        \inf_{g \text{ affine}} \mathrm{TV}(u-g) \le \mathrm{TV}(u-g_\theta) \le \frac{1}{\left|\int_\Omega \theta(z) \mathrm{d}z\right|} \mathrm{NLTV}(u).
    \end{equation}

    Finally, the bound above works for all $\theta \in \mathcal{C}_c^1(\Omega)$ satisfying $\|\theta\|_{L^\infty(\Omega)} \le 1$ and $\int_{\Omega} \theta\, \mathrm{d}z \neq 0$, 
    and the optimal bound is obtained as $\theta \to 1$ in $L^1(\Omega)$, in which case

    \begin{equation}
        \frac{1}{\left|\int_\Omega \theta(z) \mathrm{d}z\right|} \to \frac{1}{|\Omega|}.
    \end{equation}

    Thus, we have the stronger bound

    \begin{equation}
        \inf_{g \text{ affine}} \mathrm{TV}(u-g) \le \frac{1}{|\Omega|} \mathrm{NLTV}(u).
    \end{equation}
\end{proof}

\begin{remark}
    The connectedness of $\Omega$ is not required for Lemma \ref{NLBV_wIntoNLBV}, Proposition \ref{NLBV_into_BV} or Proposition \ref{TVaffNLTVProp} to hold.
\end{remark}

\subsection{Embedding into BV} \label{EmbeddingIntoBVSubsec}

To prove the existence of solutions to the minimization problem presented in Section \ref{NLChamLionsSec}, we will pass through the classical $\mathrm{BV}$ space; hence the $\mathrm{NLBV}_\omega(\Omega) \hookrightarrow \mathrm{BV}(\Omega)$ embedding will be crucial. However, condition \eqref{w_bounded_below} is not realistic in most applications, as nonlocal weights often make use of search windows, resulting in the weight vanishing on a significant part of its domain. Therefore, we introduce a weaker and more realistic condition that still allows for the continuous embedding to hold:
\begin{equation}\tag{H0} \label{H0}
    \begin{array}{c}
        \mbox{There exists a finite open and connected cover of } \Omega, ( \Omega_n )_{n =1}^M,\\\mbox{ and a constant } c>0 \mbox{ such that } 
     \omega \ge c \; \mbox{ a.e.~in } \Omega_n\times\Omega_n \mbox{ for all } n=1,..., M.
    \end{array}
\end{equation}

The connectedness condition of the covering appearing in \eqref{H0} is not needed for the $\mathrm{NLBV}_\omega(\Omega) \hookrightarrow \mathrm{BV}(\Omega)$ embedding. Instead, it will be crucial in the proof of the forthcoming Lemma \ref{ZeroNLTVCharLemma}. Condition \eqref{H0} is specifically designed with weight functions with window searches in mind, as shown in the remark below.

\begin{remark} \label{SearchWindowSatisfiesH0}
    Suppose $\Omega$ is bounded and convex (for example, a square in $\R^2$). We model a general nonlocal weight as
    \begin{equation}
        \omega(x,y) = \mathcal{X}_{(\Omega\times\Omega)_{\rho}} (x,y) f(x,y),
    \end{equation}
    where, for a fixed $\rho > 0$,
    \begin{equation}
        (\Omega\times\Omega)_\rho := \{ (x,y) \in \Omega\times\Omega : |x - y| < \rho\},
    \end{equation}
    which denotes the region (or window) of interaction, and $f:\Omega\times\Omega \to \R$ is a strictly positive function that quantifies the similarity between points. In practice, allowing each point $x\in \Omega$ to interact with every other $y \in \Omega$ is computationally unfeasible, hence why we restrict interactions to a set region.
    
    To ensure that the weight $\omega$ is compatible with our setting, we assume that $f$ is continuous in $\Omega\times\Omega$. Under this assumption, $\omega$ is an admissible weight function that satisfies

    \begin{equation}
        \omega > 0 \text{ in } (\Omega\times\Omega)_\rho \quad \text{ and } \quad  \omega = 0 \text{ outside } (\Omega\times\Omega)_\rho.
    \end{equation}

    We additionally assume that $f$ is continuously extendable to
    \begin{equation}
        (\overline\Omega \times \overline\Omega)_{\overline{\rho'}} := \{ (x,y) \in \overline\Omega \times \overline\Omega : |x - y| \le \rho'\},
    \end{equation}
    for some $\rho' < \rho$, and that the extension is strictly positive in this set. Consequently, $\omega$ also admits a continuous, strictly positive extension to $(\overline\Omega \times \overline\Omega)_{\overline{\rho'}}$.

    With this setting, let us see that $\omega$ satisfies \eqref{H0}. As $\omega$ is continuous and strictly positive on $(\overline\Omega\times\overline\Omega)_{\overline{\rho'}}$ and this set is compact,
    \begin{equation}
    c := \min_{(\overline\Omega\times\overline\Omega)_{\overline{\rho'}}} \omega > 0.
    \end{equation}
    
    Moreover, by compactness of $\overline{\Omega}$, there exist $z_1,\dots,z_{M
    } \in \overline{\Omega}$
    such that $\overline{\Omega} \subset \bigcup_{n=1}^M B(z_n, \rho'/2)$. We then define
    $\Omega_n := B(z_n, \rho'/2) \cap \Omega$, which are open and cover $\Omega$. Moreover, each $\Omega_n$ is convex by the convexity of $\Omega$, and hence connected\footnote{If the convexity assumption on $\Omega$ were replaced by it having Lipschitz boundary, such $\Omega_n$ would not necessarily be connected. To circumvent this in such a case, one can instead consider $\Omega$ intersected with cones adapted to the Lipschitz constant of the graph function of $\partial\Omega$; the Lipschitz condition ensures that such intersections lie on one side of the boundary graph, making them path-connected and thus connected.}.
    
    Finally, one readily checks that $\Omega_n \times \Omega_n \subset (\overline\Omega\times\overline\Omega)_{\overline{\rho'}}$, from which $\omega \ge c$ in each $\Omega_n \times \Omega_n$. We conclude that $\omega$ satisfies \eqref{H0}.
\end{remark}

Let us prove that if \eqref{H0} is satisfied, then the desired continuous embedding into $\mathrm{BV}$ holds. To do so, we first need some technical lemmas. We use the following notation throughout: If $D \subset \Omega$ is open and $u$ is a function defined on $\Omega$, then $u \in \mathrm{BV}(D)$ means that the restriction of $u$ to $D$ belongs to $\mathrm{BV}(D)$. Similarly, $u \in \mathrm{NLBV}_\omega(D)$ means that the restriction of $u$ to $D$ belongs to $\mathrm{NLBV}_{\omega|_{D\times D}}(D)$. Furthermore, $\mathrm{TV}^D(u)$ denotes the total variation of $u$ in $D$, $\mathrm{NLTV}^D_\omega(u)$ denotes the weighted nonlocal total variation of $u$ in $D$ induced by the weight function $\omega|_{D\times D}$, and $\mathrm{NLTV}^D(u)$ denotes the weighted nonlocal total variation of $u$ in $D$ induced by the constant weight function $1$ on $D\times D$.

\begin{lemma} \label{TechUnionLemmaBV}
    Let $u\in L^1_{\mathrm{loc}}(\Omega)$ and $A,B\subset \Omega$ open and non-empty. If $u \in \mathrm{BV}(A)$ and $u\in \mathrm{BV}(B)$, then $u \in \mathrm{BV}(A\cup B)$ and
    \begin{equation}
        \mathrm{TV}^{A\cup B}(u) \le \mathrm{TV}^A(u) + \mathrm{TV}^B(u),\quad \|u\|_{\mathrm{BV}({A\cup B})} \le\|u\|_{\mathrm{BV}({A})} + \|u\|_{\mathrm{BV}({B})}.
    \end{equation}
    Moreover, if $u \in \mathrm{BV}(A\cup B)$, then $u \in \mathrm{BV}(A)$ and $u\in \mathrm{BV}(B)$, and
    \begin{equation} \label{TechUnionBVEq}
        \mathrm{TV}^A(u) \le \mathrm{TV}^{A\cup B} (u),\:\:\: \mathrm{TV}^B(u) \le \mathrm{TV}^{A\cup B} (u),\:\:\: \|u\|_{\mathrm{BV}({A})} \le \|u\|_{\mathrm{BV}({A\cup B})}, \:\:\: \|u\|_{\mathrm{BV}({B})} \le \|u\|_{\mathrm{BV}({A\cup B})}.
    \end{equation}
\end{lemma}

\begin{proof}
    The result directly follows from the fact that, for any $D\subset \Omega$ open
    \begin{equation}
        \mathrm{TV}^D(u) = \mu_u(D)
    \end{equation}
    for a certain $\R^N$-valued Radon measure $\mu_u$ in $\Omega$ (see, e.g., \cite[Proposition 3.6]{AmbrosioFuscoPallara}), from measure properties, and from the inequalities

    \begin{equation} \label{L^1setinequalities}
        \|u\|_{L^1(A\cup B)} \le \|u\|_{L^1(A)} + \|u\|_{L^1(B)}, \quad \|u\|_{L^1(A)} \le \|u\|_{L^1(A\cup B)}, \quad\|u\|_{L^1(B)} \le \|u\|_{L^1(A\cup B)}.
    \end{equation}
\end{proof}

In the nonlocal case, only the second statement of Lemma \ref{TechUnionLemmaBV} (i.e. \eqref{TechUnionBVEq}) holds, as we can only assure that $(A\times A) \cup (B\times B) \subset (A\cup B) \times (A\cup B)$, with the equality not holding in general. Moreover, we use a distributional argument rather than one based on measures. We note that, in contrast to the local setting, the definition of nonlocal divergence depends on the domain, and for $D \subset \Omega$ open, we will write $\mathrm{div}_{\mathrm{NL}}^D$ to denote the nonlocal divergence  in $D$. 

\begin{lemma} \label{TechUnionLemmaNLBV}
    Let $\omega$ be an admissible weight function and $A,B\subset \Omega$ open. Then $\mathrm{NLBV}_\omega(A\cup B) \hookrightarrow \mathrm{NLBV}_\omega(A)$ and $\mathrm{NLBV}_\omega(A\cup B) \hookrightarrow \mathrm{NLBV}_\omega(B)$. More precisely,
    \begin{gather}
        \mathrm{NLTV}_\omega^A(u) \le \mathrm{NLTV}_\omega^{A \cup B}(u), \quad\|u\|_{\mathrm{NLBV}_\omega(A)} \le \|u\|_{\mathrm{NLBV}_\omega(A\cup B)},\\
        \mathrm{NLTV}_\omega^B(u) \le \mathrm{NLTV}_\omega^{A \cup B}(u), \quad\|u\|_{\mathrm{NLBV}_\omega(B)} \le \|u\|_{\mathrm{NLBV}_\omega(A\cup B)},
    \end{gather}
    for all $u \in \mathrm{NLBV}_\omega(A\cup B)$.
\end{lemma}

\begin{proof}
    Let $\phi \in \mathcal{C}_c^1(A\times A;\R^N)$ be such that $|\phi| \le \omega$ a.e.~in $A\times A$ and extend it by zero to a function $\tilde \phi \in \mathcal{C}_c^1((A\cup B)\times (A\cup B);\R^N)$ such that $| \tilde \phi| \le \omega$ a.e.~in $(A\cup B)\times (A\cup B)$. Note that then, for every $x\in A$

    \begin{align}
        \mathrm{div}_{\mathrm{NL}}^A(\phi)(x) &= \int_A\mathrm{div}_x \left[ \phi(y,x) - \phi(x,y) \right] \mathrm{d}y=\int_{A\cup B} \mathrm{div}_x [ \tilde \phi(y,x) - \tilde \phi(x,y) ] \mathrm{d}y = \mathrm{div}_{\mathrm{NL}}^{A\cup B}(\tilde \phi)(x).
    \end{align}
    Moreover, as $\tilde \phi$ has support in $A\times A$, we can extend $\mathrm{div}_{\mathrm{NL}}^{A\cup B}(\tilde \phi)$ by zero to the whole $A\cup B$. Using this, we have
    \begin{align}
        \int_A u(x) \; \mathrm{div}_{\mathrm{NL}}^A(\phi)(x)\mathrm{d}x &= \int_{A\cup B} u(x) \; \mathrm{div}_{\mathrm{NL}}^{A\cup B}(\tilde \phi)(x)\mathrm{d}x \le \mathrm{NLTV}_\omega^{A\cup B}(u),
    \end{align}
    and by taking supremum over all $\phi$ we obtain that
    \begin{equation}
        \mathrm{NLTV}_\omega^A(u) \le \mathrm{NLTV}_\omega^{A \cup B}(u).
    \end{equation}

    Analogously, one sees that
    \begin{equation}
        \mathrm{NLTV}_\omega^B(u) \le \mathrm{NLTV}_\omega^{A \cup B}(u),
    \end{equation}
    which, together with the previous inequalities \eqref{L^1setinequalities}, completes the proof.
\end{proof}

With these lemmas proven, we move onto the desired embedding.

\begin{prop} \label{NLBVintoBVAltCondBV}
    Let $\omega$ be an admissible weight function satisfying \eqref{H0}. Then $\mathrm{NLBV}_\omega(\Omega) \hookrightarrow \mathrm{BV}(\Omega)$.
\end{prop}

\begin{proof}
    Let $(\Omega_n )_{n=1}^M$ be the finite open and connected cover of $\Omega$ and $c>0$ the constant given by \eqref{H0}. We start by using Proposition \ref{NLBV_into_BV} to obtain that, since $\omega \ge c$ a.e.~in each $\Omega_n\times\Omega_n$, then $\mathrm{NLBV}_\omega(\Omega_n) \hookrightarrow \mathrm{BV}(\Omega_n)$ for each $n=1,...,M$. Hence, for each $n=1,...,M$, there exists $C_n>0$ such that
    \begin{equation}
        \|u\|_{\mathrm{BV}(\Omega_n)} \le C_n \|u\|_{\mathrm{NLBV}_\omega(\Omega_n)}
    \end{equation}
    for all $u \in \mathrm{NLBV}_\omega(\Omega_n)$. Thus, by using Lemma \ref{TechUnionLemmaBV} iteratively,

    \begin{align}
        \|u\|_{\mathrm{BV}(\Omega)} &\le \sum_{n=1}^M \|u\|_{\mathrm{BV}(\Omega_n)} \le \sum_{n=1}^M C_n \|u\|_{\mathrm{NLBV}_\omega(\Omega_n)} \le  \left( \sum_{n=1}^M C_n \right)  \|u\|_{\mathrm{NLBV}_\omega(\Omega)},
    \end{align}
    for $u \in \mathrm{NLBV}_\omega(\Omega)$, where in the last inequality we have used Lemma \ref{TechUnionLemmaNLBV} (iteratively).
\end{proof}

\begin{remark}
    The connectedness assumption on $\Omega$ is not required for Proposition \ref{NLBVintoBVAltCondBV} to hold.
\end{remark}

\subsection{Nonlocal Total Variation and Affine Functions}

It can be readily seen from its definition that $\mathrm{NLTV}_\omega$ annihilates affine functions. What is more, if $\omega$ satisfies \eqref{H0}, the implication becomes an equivalence. We establish this characterization in the following lemma, whose proof relies on the notation introduced immediately before Lemma \ref{TechUnionLemmaBV}.

\begin{lemma}
\label{ZeroNLTVCharLemma}
    Let $\omega$ be an admissible weight function which satisfies \eqref{H0}. Then, for all $u \in \mathrm{NLBV}_\omega(\Omega)$

    \begin{equation}
        \mathrm{NLTV}_\omega(u) = 0 \iff u \text{ is affine a.e.~in } \Omega.
    \end{equation}
\end{lemma}

\begin{proof}
    The right to left implication follows immediately from Proposition \ref{C1NLTV}. Let us prove the other implication. Suppose that $u \in \mathrm{NLBV}_\omega(\Omega)$ is such that $\mathrm{NLTV}_\omega(u) = 0$ and let $(\Omega_n )_{n=1}^M$ be the finite open and connected cover of $\Omega$ and $c>0$ the constant given by \eqref{H0}. By Lemma \ref{NLBV_wIntoNLBV}, for each $n=1,...,M$ we have
    
    \begin{equation}
        \mathrm{NLTV}^{\Omega_n}(u) \le \frac{1}{c} \mathrm{NLTV}_\omega^{\Omega_n}(u).
    \end{equation}
    Moreover, using Lemma \ref{TechUnionLemmaNLBV} we know that

    \begin{equation}
        \mathrm{NLTV}_\omega^{\Omega_n}(u) \le \mathrm{NLTV}^\Omega_\omega(u).
    \end{equation}
    
    By Proposition \ref{TVaffNLTVProp} we have that for each $n=1,...,M$ there exists an affine function $g_n$ defined in $\Omega_n$ such that

    \begin{equation}
        \mathrm{TV}^{\Omega_n}(u-g_n) \le \frac{1}{|\Omega_n|} \mathrm{NLTV}^{\Omega_n}(u) \le \frac{1}{c |\Omega_n|} \mathrm{NLTV}_\omega^{\Omega_n}(u) \le \frac{1}{c |\Omega_n|} \mathrm{NLTV}^\Omega_\omega(u) = 0.
    \end{equation}
By \cite[Proposition 3.2]{AmbrosioFuscoPallara}, this means that for each $n=1,...,M$, $u-g_n$ must be equal a.e.~to a constant in $\Omega_n$, from which $u$ is equal a.e.~to an affine function, say $\tilde g_n$, in $\Omega_n$.

We claim now that we can reorder the open covering $(\Omega_n)_{n=1}^M$ in a way that for every $i = 2,...,M$, $\Omega_i \cap \Omega_j \neq \emptyset$ for some $j<i$. To prove this, consider the graph whose nodes correspond to each $\Omega_n$, and such that an edge connects two nodes $\Omega_i$ and $\Omega_j$, $i \neq j$, if and only if $\Omega_i \cap \Omega_j \neq \emptyset$. The resulting graph is connected. Indeed, suppose that there exists a non-trivial partition of $\{1,...,M\}$ into two non-empty sets $A$ and $B$ such that no edge of the graph crosses between $A$ and $B$, that is
\begin{equation}
    \Omega_i \cap \Omega_j = \emptyset \qquad \text{ for } \qquad i \in A, j \in B.
\end{equation}
Then $\Omega$ could be split into two disjoint open and non-empty subsets $\cup_{i\in A} \Omega_i$ and $\cup_{j \in B} \Omega_j$, contradicting the connectedness of $\Omega$.

Now, since the graph is connected, it must contain a walk that visits each node at least once. Reordering $\Omega_n$ according to this walk (assigning the first appearing index in the walk to nodes appearing more than once in the walk) proves the claim.

With the claim proven, as $u$ is affine in $\Omega_1$ and $\Omega_2$, $\Omega_1 \cap \Omega_2 \neq \emptyset$, and both subsets are connected, this affine function must be the same, that is $\tilde g_1 \equiv \tilde g_2$ in $\Omega_1 \cup \Omega_2$. Moreover, by the claim, either $\Omega_3 \cap \Omega_1 \neq \emptyset$ or $\Omega_3 \cap \Omega_2 \neq \emptyset$. In either case we deduce that $\tilde g_1 \equiv \tilde g_2 \equiv \tilde g_3$ in $\Omega_1 \cup \Omega_2 \cup \Omega_3$. Repeating this argument iteratively for each $n=4,...,M$, we deduce that $u$ is equal a.e.~to the same affine function in $\cup_{n=1}^M \Omega_n = \Omega$, and the proof is complete.
\end{proof}

\begin{figure}
    \centering
    \begin{subfigure}[t]{0.18\linewidth}
        \centering
        \spyTGV{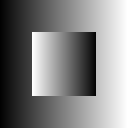}
        \caption*{\centering Synthetic image\\
        PSNR} 
    \end{subfigure}\hfill
    \begin{subfigure}[t]{0.18\linewidth}
        \centering
        \spyTGV{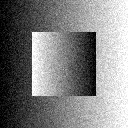}
        \caption*{\centering Noisy image\\
        24.96}
    \end{subfigure}\hfill
    \begin{subfigure}[t]{0.18\linewidth}
        \centering
        \spyTGV{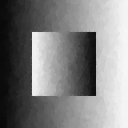}
        \caption*{\centering $\mathrm{TV}$\\
        35.39} 
    \end{subfigure}\hfill
    \begin{subfigure}[t]{0.18\linewidth}
        \centering
        \spyTGV{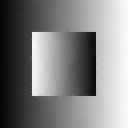}
        \caption*{\centering $\mathrm{TGV}^2$\\
        39.19}
    \end{subfigure}\hfill
    \begin{subfigure}[t]{0.18\linewidth}
        \centering
        \spyTGV{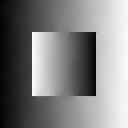}
        \caption*{\centering $\mathrm{NLTV}_\omega$\\
        \textbf{41.05}}
    \end{subfigure}\hfill
    \caption{Restoration of a synthetic noisy ($\sigma_{\mathrm{noise}} = 15$) piecewise affine image using $\mathrm{TV}$, $\mathrm{TGV}^2$ and $\mathrm{NLTV}_\omega$ regularization. $\mathrm{NLTV}_\omega$ achieves the highest PSNR and recovers the original image with remarkable accuracy, whereas $\mathrm{TGV}^2$ exhibits mild residual noise around the edges of the inner square and $\mathrm{TV}$ displays pronounced staircasing artifacts.
    }
    \label{fig:piecewiseaff}
\end{figure}

Consequently, $\mathrm{NLTV}_\omega$ regularization promotes piecewise affine functions and mitigates the staircase effect. This behavior is shown numerically in Figure \ref{fig:piecewiseaff}. In this experiment, a synthetic piecewise affine image is corrupted with gaussian noise with standard deviation $\sigma_{\mathrm{noise}} = 15$ and restored using $\mathrm{TV}$, $\mathrm{TGV}^2$ \cite{bredies2010total} and $\mathrm{NLTV}_\omega$. The peak signal-to-noise ratio (PSNR) is reported for the noisy image and each restoration result. We use the primal-dual algorithms associated with the Rudin–Osher–Fatemi (ROF) \cite{Rudin-Osher-Fatemi} models for each regularizer. For $\mathrm{NLTV}_\omega$, we employ the NLMeans weights defined in \eqref{NLMeansWeightDefDiscrete}. As expected, $\mathrm{TV}$ regularization introduces pronounced staircasing artifacts, whereas both $\mathrm{NLTV}_\omega$ and $\mathrm{TGV}^2$ regularization recover the original synthetic image with high accuracy. However, the $\mathrm{TGV}^2$ reconstruction exhibits more residual noise around the edges of the inner square than the $\mathrm{NLTV}_\omega$ reconstruction, which is specially visible in the zoomed-in crop. These qualitative observations are consistent with the PSNR values: $\mathrm{NLTV}_\omega$ achieves the highest PSNR, outperforming $\mathrm{TGV}^2$, the second-best method, by around $2 \; \mathrm{dB}$. In contrast, the $\mathrm{TV}$ reconstruction yields a substantially lower PSNR.

\section{Nonlocal Chambolle-Lions model} \label{NLChamLionsSec}

In this section, we consider the following minimization problem:
\begin{equation} \label{NLChamLionVar}
    \min_{u \in L^p(\Omega) \cap \mathrm{NLBV}_\omega(\Omega)} \mathrm{NLTV}_\omega(u) + \frac{\lambda}{2} \int_\Omega |Au - u_0|^2 \mathrm{d}x,
\end{equation}
where $\Omega \subset \R^N$ is a nonempty, bounded, connected and Lipschitz set, $\omega$ is an admissible weight function, $A: L^p(\Omega) \to L^p(\Omega)$, and $p$ is a certain exponent in $[1,+\infty)$ depending on $N$. Recall that this problem is essentially a nonlocal version of the classical model proposed by Chambolle and Lions \cite{chambolle1997image}.  The original hypothesis appearing in \cite{chambolle1997image} are the following: $p=2$ if $N=1$ and $p=N/(N-1)$ if $N\ge 2$ (this ensures that $\mathrm{BV}(\Omega) \hookrightarrow L^p(\Omega)$ by Sobolev embedding), and

\begin{equation} \tag{H1} \label{H1}
    A:L^p(\Omega) \to L^p(\Omega) \text{ is a continuous and linear operator},
\end{equation}
\begin{equation} \tag{H2} \label{H2}
    A(1) \neq 0.
\end{equation}

Notice that, at least a priori, we \emph{do not} have the $\mathrm{NLBV}_\omega(\Omega) \hookrightarrow L^p(\Omega)$ embedding, hence the solution function space in the problem above needs to be $L^p(\Omega) \cap \mathrm{NLBV}_\omega(\Omega)$.

\begin{remark}
In the original paper \cite{chambolle1997image}, the authors take $A(1)=1$ as the main assumption on $A$, rather than \eqref{H2}. While \eqref{H2} is sufficient to establish the existence and uniqueness of solutions, the condition $A(1)=1$ provides a more direct connection with a related constrained problem. Since the focus of the present work is the variational problem, we present condition \eqref{H2} as the main one instead.
\end{remark}

All assumptions except \eqref{H2} seem to still be reasonable and natural in the nonlocal setting. What is more, under assumption \eqref{H1}, the functional $F:L^p(\Omega) \to \mathbb{R}$ given by $F(u):= \frac{\lambda}{2} \int_\Omega |Au - u_0|^2 \mathrm{d}x$ is proper, bounded from below, convex (in fact, strictly convex with respect to $Au$) and lower semi-continuous (more specifically, weakly lower semi-continuous with respect to weak $L^2(\Omega)$ convergence). To use the Direct Method, the compactness of $\mathrm{NLBV}_\omega(\Omega)$ is needed, which is not true in general (consider, for example, the limit case $\omega \equiv 0$ in $\Omega\times\Omega$). We additionally need a Sobolev-type embedding for $\mathrm{NLBV}_\omega(\Omega)$, in order to ensure that the fidelity term of the minimizing sequence also converges. To get these properties, we will pass through $\mathrm{BV}(\Omega)$.

Moreover, coercivity of the functional is also necessary, which in \cite{chambolle1997image} follows from hypothesis \eqref{H2}. Conceptually, the $\mathrm{TV}$ annihilates constant functions, so, without hypothesis \eqref{H2}, the $L^p$ norm of $u+c$ could grow arbitrarily by choosing increasingly large $c>0$, while the functional in \eqref{OgChamLionVar} would stay bounded. In our setting, $\mathrm{NLTV}_\omega$ instead annihilates affine functions (recall Lemma \ref{ZeroNLTVCharLemma}). So, theoretically, even if $A(1)\neq 0$, without control on $x$ one could construct a sequence $u + a_n\cdot x$ which grows arbitrarily in $L^p$ norm and such that both $\mathrm{NLTV}_\omega(u+a_n\cdot x)$ and $\|A(u+a_n\cdot x)-u_0\|$ are bounded, as the next example shows.

\begin{ex}
    In this example we limit ourselves to $N=1$ (and thus $p=2$) and $\Omega = (-1, 1)$. Consider the continuous and linear operator $A:L^2(-1,1) \to L^2(-1,1)$ given by

    $$Au = \frac{1}{2} \int_{-1}^1 u(x) \mathrm{d}x$$
    for any $u \in L^2(-1,1)$. Note that $A(1) = 1 \neq 0$, but one quickly checks that $A(x) = 0$, and so $A(c x) = 0$ independently of $c \in \R$.
\end{ex}

In view of the example above, hypothesis \eqref{H2} will need to be modified to guarantee coercivity.

\subsection{Nonlocal Sobolev and Compact Embeddings and Coercivity of the Functional}

The  $\mathrm{NLBV}_\omega(\Omega) \hookrightarrow \mathrm{BV}(\Omega)$ embedding which we have proven in Subsection \ref{EmbeddingIntoBVSubsec} allows for the nonlocal Sobolev embedding and compactness to hold in this case, as shown next.

\begin{lemma}[Sobolev embedding and compactness for $\mathrm{NLBV}_\omega(\Omega)$] \label{LemmaSobEmbwboundedbelow}
    Let $N\ge 2$, $\Omega \subset \R^N$ be bounded, connected and Lipschitz, $\omega$ be an admissible weight function satisfying \eqref{H0} and $p=N/(N-1)$. Then we have the continuous embedding

    \begin{equation}
        \mathrm{NLBV}_\omega(\Omega) \hookrightarrow L^q(\Omega)
    \end{equation}
    for $1\le q \le p$. Moreover, the embedding is also compact for $1 \le q < p$.
\end{lemma}

\begin{proof}
    For $1 \le q \le p$ we have $\mathrm{BV}(\Omega) \hookrightarrow L^q(\Omega)$. Coupling this with the $\mathrm{NLBV}_\omega(\Omega) \hookrightarrow \mathrm{BV}(\Omega)$ embedding (which is given by Proposition \ref{NLBVintoBVAltCondBV}) yields the first part of the result.

    Moreover, for $1 \le q < p$, let us consider a sequence which is uniformly bounded in $\mathrm{NLBV}_\omega(\Omega)$. By the aforementioned embedding, it must also be bounded in $\mathrm{BV}(\Omega)$, and since $\mathrm{BV}(\Omega) \hookrightarrow L^q(\Omega)$ is compact in this case, it admits a subsequence converging strongly in $L^q(\Omega)$. Hence the embedding $\mathrm{NLBV}_\omega(\Omega) \hookrightarrow L^q(\Omega)$ is compact, which completes the proof.
\end{proof}

\begin{remark} \label{RemarkSobEmbwboundedbelow}
    The proof above can be readily modified to show that if $N=1$, $\Omega \subset \R$ is bounded, connected and Lipschitz, and $\omega$ is an admissible weight function satisfying \eqref{H0}, then
    \begin{equation}
        \mathrm{NLBV}_\omega(\Omega) \hookrightarrow L^p(\Omega)
    \end{equation}
    for $1\le p \le +\infty$, and the embedding is compact for $1\le p < +\infty$. As in \cite{chambolle1997image}, in this case we choose $p=2$.
\end{remark}

With the Sobolev embedding and compactness for $\mathrm{NLBV}_\omega(\Omega)$ and Lemma \ref{ZeroNLTVCharLemma}, the proof of coercivity can be reproduced in the nonlocal setting provided we modify the \eqref{H2} hypothesis on $A$ as follows: 
\begin{equation} \tag{H2'} \label{H2'}
    \forall a\in \R^N, b \in \R: A(a\cdot x + b) = 0 \; \text{ a.e.~in } \Omega \iff a=0, b=0.
\end{equation}
We remark that this is precisely the condition Kindermann et~al.~impose in their paper (see \cite[Theorem 5.5]{Kindermann}).

\begin{lemma}[$L^1$ coercivity of the functional appearing in \eqref{NLChamLionVar}] \label{NLChamLionCoerciv}
    Let $\Omega \subset \R^N$ be bounded, connected and Lipschitz, $p=2$ if $N=1$ and $p=N/(N-1)$ if $N\ge 2$, and $\omega$ be an admissible weight function satisfying \eqref{H0}. Suppose $A$ satisfies \eqref{H1} and \eqref{H2'}. Then there exists $C>0$ such that

    \begin{equation}
        \|u\|_1 \le C (\mathrm{NLTV}_\omega(u) + \|A u\|_2)
    \end{equation}
    for all $u \in \mathrm{NLBV}_\omega(\Omega)$.
\end{lemma}

\begin{proof}
     We argue by contradiction. Suppose there exists $(u_n)_{n\in \N} \subset \mathrm{NLBV}_\omega(\Omega)$ such that for each $n \in \N$

    \begin{equation}
        \|u_n\|_1 > n (\mathrm{NLTV}_\omega(u_n) + \|A u_n\|_2).
    \end{equation}

    In particular, none of the $u_n$ can be zero. As the inequality is homogeneous of degree 1, by dividing each side by $\|u_n\|_1$ and defining $u_n' := u_n/\|u_n\|_1$, we obtain

    \begin{equation} \label{NLCoercAux0}
        \mathrm{NLTV}_\omega(u_n') + \|A u_n'\|_2 < \frac{\|u_n'\|_1}{n} = \frac{1}{n}
    \end{equation}
    for each $n\in\N$. Note that clearly $(u_n')_{n\in\N}$ is uniformly bounded in $\mathrm{NLBV}_\omega(\Omega)$, and thus, by Lemma \ref{LemmaSobEmbwboundedbelow} (or Remark \ref{RemarkSobEmbwboundedbelow} in the $N=1$ case), we can extract a subsequence (which we still denote by $(u_n')_{n\in\N}$) such that $u_n' \to u'$ in $L^1(\Omega)$ for a certain $u' \in \mathrm{NLBV}_\omega(\Omega)$. Moreover, by lower semi-continuity and \eqref{NLCoercAux0}, we have

    \begin{equation}
        \mathrm{NLTV}_\omega(u') \le \liminf_{n\to+\infty} \mathrm{NLTV}_\omega(u_n') \le \liminf_{n\to+\infty} \frac{1}{n} = 0, 
    \end{equation}
    and by Lemma \ref{ZeroNLTVCharLemma} we must have that $u'$ is an affine function.

    Now, recall that we have the continuous embedding $\mathrm{NLBV}_\omega(\Omega) \hookrightarrow L^p(\Omega)$ (Lemma \ref{LemmaSobEmbwboundedbelow} for $N \ge 2$ or Remark \ref{RemarkSobEmbwboundedbelow} for $N=1$).  Thus, as $(u_n')_{n\in\N}$ is uniformly bounded in $\mathrm{NLBV}_\omega(\Omega)$, it is also uniformly bounded in $L^p(\Omega)$ and we can extract a subsequence $(u_{n_l}')_{l\in \N}$ such that $u_{n_l}' \rightharpoonup u''$ in $L^p(\Omega)$ for a certain $u''\in L^p(\Omega)$. By uniqueness of the distributional limit, we must have $u'' = u'$ a.e.~in $\Omega$.

    By continuity of $A$ we have that $A u_{n_l}' \rightharpoonup Au'$ in $L^p(\Omega)$. Furthermore, $A u_n' \to 0$ strongly in $L^2(\Omega)$, as, by \eqref{NLCoercAux0}

    \begin{equation}
        \| A u_n' \|_2 \le \frac{1}{n} \xrightarrow[n\to \infty]{} 0.
    \end{equation}
    Again, by uniqueness of the distributional limit, we must have that $A u' = 0$ a.e.~in $\Omega$. Therefore, as $u'$ is affine, $u'=0$ using hypothesis \eqref{H2'}, which contradicts the fact that $\|u'\|_1 = 1$ (as $u_n' \to u'$ in $L^1(\Omega)$ and $\|u_n'\|_1 = 1$ for all $n \in \N$).
\end{proof}

\subsection{Proof of Existence and Uniqueness of Solutions}

With compactness and Sobolev embedding of $\mathrm{NLBV}_\omega$ and the coercivity result proven, we may reproduce \cite{chambolle1997image} to prove existence of the variational problem \eqref{NLChamLionVar}. The argument is based on a straightforward application of the Direct Method.

\begin{theorem}[Existence and uniqueness of minimizers of \eqref{NLChamLionVar}] \label{ExistUniqNLChamLionVar}
    Let $N \ge 1$, $\Omega\subset \R^N$ be bounded, connected and Lipschitz, $\lambda > 0$, $\omega$ an admissible weight function satisfying \eqref{H0} and $p=2$ if $N=1$ or $p=N/(N-1)$ if $N\ge 2$. Suppose that $u_0 \in L^2(\Omega)$ and that $A$ satisfies \eqref{H1} and \eqref{H2'}. Then there exists a minimizer $u \in \mathrm{NLBV}_\omega(\Omega)$ of \eqref{NLChamLionVar}. Additionally, if $A$ is injective, the minimizer is unique.
\end{theorem}

\begin{proof}
    Let $(u_n)_{n\in\N} \subset \mathrm{NLBV}_\omega(\Omega)$ be a minimizing sequence of \eqref{NLChamLionVar}, i.e., such that
    
    \begin{equation}
        \lim_{n\to +\infty}  \left[\mathrm{NLTV}_\omega(u_n) + \frac{\lambda}{2} \int_\Omega |Au_n - u_0|^2 \mathrm{d}x \right] = \inf_{u \in \mathrm{NLBV}_\omega(\Omega)} \left[  \mathrm{NLTV}_\omega(u) + \frac{\lambda}{2} \int_\Omega |Au - u_0|^2 \mathrm{d}x\right] =: m.
    \end{equation}

    The constant $m$ is finite, as the functional is lower bounded by $0$ and substituting $u$ by the zero function yields a finite upper bound ($(\lambda/2)\|u_0\|_2^2$).
    
    Note that $(u_n)_{n\in\N}$ is uniformly bounded in $\mathrm{NLBV}_\omega(\Omega)$. Indeed, the $\mathrm{NLTV}_\omega$ bound is derived from the boundedness of the functional in \eqref{NLChamLionVar}. Moreover, by Lemma \ref{NLChamLionCoerciv}, there exists $C>0$ such that for all $n\in\N$

    \begin{equation} \label{NLExistVarwboundedbelowAux0}
        \|u_n\|_1 \le C (\mathrm{NLTV}_\omega(u_n) + \|A u_n\|_2),
    \end{equation}
    and

    \begin{equation} \label{NLExistVarwboundedbelowAux1}
        \|A u_n\|_2 \le \|A u_n - u_0\|_2 + \|u_0\|_2 
    \end{equation}
    whose right-hand side is uniformly bounded by the fact that $(u_n)_{n\in\N}$ is a minimizing sequence and $u_0 \in L^2(\Omega)$. Coupling this with the uniform bound in the $\mathrm{NLTV}_\omega$ semi-norm, we deduce that the right-hand side of \eqref{NLExistVarwboundedbelowAux0} is uniformly bounded.

    Hence, $(u_n)_{n\in\N}$ is uniformly bounded in $\mathrm{NLBV}_\omega(\Omega)$ and by compactness (Lemma \ref{LemmaSobEmbwboundedbelow} for $N \ge 2$ or Remark \ref{RemarkSobEmbwboundedbelow} for $N=1$) we can extract a subsequence (which we still denote by $(u_n)_{n\in\N}$) such that $u_n \to u$ in $L^1(\Omega)$ for a certain $u \in L^1(\Omega)$. In fact, we have $u \in \mathrm{NLBV}_\omega(\Omega)$ by lower semi-continuity (Lemma \ref{lemma:lsc_nlbv}).

    Now, by Sobolev embedding (Lemma \ref{LemmaSobEmbwboundedbelow} for $N \ge 2$ or Remark \ref{RemarkSobEmbwboundedbelow} for $N=1$), $(u_n)_{n\in\N}$ (the subsequence) is also uniformly bounded in $L^p(\Omega)$, and we can extract a subsequence (which, again, we still denote by $(u_n)_{n\in\N}$) such that $u_n \rightharpoonup u'$ in $L^p(\Omega)$ for a certain $u' \in L^p(\Omega)$. By uniqueness of the distributional limit we must have $u = u'$ a.e.~in $\Omega$, and by \eqref{H1}, $A u_n \rightharpoonup A u$ in $L^p(\Omega)$. What is more, by the uniform bound of \eqref{NLExistVarwboundedbelowAux1}, we know that there exists $v \in L^2(\Omega)$ such that $A u_n \rightharpoonup v$ in $L^2(\Omega)$ (up to a subsequence). Again, by uniqueness of the distributional limit, we must have that $v = Au$ a.e.~in $\Omega$.

    Finally, by lower semi-continuity of $\mathrm{NLTV}_\omega$ with respect to $L^1(\Omega)$ convergence (Lemma \ref{lemma:lsc_nlbv}), weak lower semi-continuity of the $L^2(\Omega)$ norm with respect to (weak) $L^2(\Omega)$ convergence, and the superadditivity of the limit inferior, we have

    \begin{equation}
         m \le \mathrm{NLTV}_\omega(u) + \frac{\lambda}{2} \int_\Omega |Au - u_0|^2 \mathrm{d}x \le \liminf_{n\to +\infty} \left[ \mathrm{NLTV}_\omega(u_n) + \frac{\lambda}{2} \int_\Omega |Au_n - u_0|^2 \mathrm{d}x\right] = m,
    \end{equation}
    from which $u$ is a minimizer for \eqref{NLChamLionVar}.

    We move onto proving uniqueness. Let $u, v \in \mathrm{NLBV}_\omega(\Omega)$ be minimizers of \eqref{NLChamLionVar}. Then $(u+v)/2 \in \mathrm{NLBV}_\omega(\Omega)$ and, by convexity of $\mathrm{NLTV}_\omega$ and $\|A \cdot - u_0\|_2^2$,
    \begin{align}
        m &\le \mathrm{NLTV}_\omega\left( \frac{u+v}{2} \right) + \frac{\lambda}{2} \left\|A \left( \frac{u+v}{2} \right) - u_0 \right\|_2^2\\
        &\le \frac{1}{2} \left( \mathrm{NLTV}_\omega(u) + \frac{\lambda}{2} \|Au-u_0\|_2^2 \right) + \frac{1}{2} \left( \mathrm{NLTV}_\omega(v) + \frac{\lambda}{2} \|Av-u_0\|_2^2 \right) = \frac{m}{2} + \frac{m}{2} = m.
    \end{align}
    As $\|\cdot - u_0\|_2^2$ is strictly convex, equality in the second inequality implies that $Au=Av$. Therefore, the image $Au$ is the same for every minimizer $u$ of \eqref{NLChamLionVar}. In particular, if $A$ is injective, then $u=v$, and hence the minimizer is unique.
\end{proof}

\subsection{Numerical Implementation} \label{NumAlgSec}

In this subsection, we discretize problem \eqref{NLChamLionVar} and the $\mathrm{NLTV}_\omega$ functional, and present an algorithm for computing its solutions. We restrict our discussion to grayscale images, as the extension to color images is straightforward when channel coupling is fixed \cite{duran2016collaborative}.

We represent a grayscale image by $u \in \R^M$, where $M$ is the number of pixels, using flattened indexing, and the dual variables by $p \in \mathbb{R}^{M \times M \times 2}$. We denote the pixel indices by $i,j \in \{1,...,M\}$, and the component index by $t \in \{1,2\}$. The gradient $\nabla = (\partial_x, \partial_y) \colon \R^M \to \R^{M \times 2}$ is computed by forward differences with Neumann boundary conditions, and the divergence $\mathrm{div} \colon \R^{M \times 2} \to \R^M$ is defined as the negative adjoint of the gradient, $\mathrm{div} = -\nabla^*$. We endow $\R^M$ with the norm
\begin{equation} \label{RMEuclideanNorm}
    \|w\|_2 := \sqrt{ \sum_{i=1}^M \left(w_{i}\right)^2},
\end{equation}
for $w \in \R^{M}$, and $\R^{M\times M \times 2}$ with the norms
$$
    \|p\|_{1} := \sum_{i=1}^M \sum_{j=1}^M | p_{i,j,\cdot}| \quad \text{and}\quad     \|p\|_2 := \sqrt{\sum_{i=1}^M \sum_{j=1}^M | p_{i,j,\cdot}|^2}
$$
for $p \in \R^{ M\times M \times 2}$, where $|\cdot|$ denotes the Euclidean norm in $\R^2$. Both the $\|\cdot\|_2$ norms in $\R^{ M}$ and $\R^{ M\times M \times 2}$ come from the usual inner product in each of these spaces, which we will denote by $\langle \cdot, \cdot \rangle$ interchangeably when it is clear from context. We also define
\begin{equation} \label{WeightedDiscreteNorm}
    \|p\|_{\omega,1} := \sum_{i=1}^M \sum_{j=1}^M \omega_{i,j}\, |p_{i,j,\cdot}|,
\end{equation}
for $p \in \R^{M \times M \times 2}$, where $\omega \in \R^{M \times M}$ is the nonlocal weight, whose entries are assumed to be nonnegative. Since $\omega$ may vanish at some entries, \eqref{WeightedDiscreteNorm} is in general only a semi-norm, so we use $\|\cdot\|_{\omega,1}$ for notational convenience only.

With this, \eqref{NLChamLionVar} is written in the discrete setting as

\begin{equation} \label{NLChamLionDiscrete}
    \min_u \|\nabla_{\mathrm{NL}} u\|_{\omega, 1} + \frac{\lambda}{2} \| A u - u_0 \|_2^2,
\end{equation}
where $A \colon \R^M \to \R^{M_0}$ is a linear operator, $u_0 \in \R^{M_0}$ is the observed image and $\|\cdot\|_2$ is the Euclidean norm in $\R^{M_0}$ as in \eqref{RMEuclideanNorm} (with $M_0 = M$ unless $A$ involves downsampling, as in super-resolution, where $M_0 = M/s^2$). 
The discrete nonlocal gradient $\nabla_{\mathrm{NL}} u \in \R^{M \times M \times 2}$ is defined by
\begin{equation}
    (\nabla_{\mathrm{NL}} u)_{i,j,t} := (\nabla u)_{j,t} - (\nabla u)_{i,t}
\end{equation}
for $i,j \in \{1,...,M\}$ and $t \in \{1,2\}$. Note that $\|\nabla_{\mathrm{NL}} \cdot\|_{\omega,1}$ corresponds to the $\mathrm{NLTV}_\omega$ functional in accordance with formula \eqref{NLTVFormula}.

In the following, we will also need the adjoint operator of $\nabla_{\mathrm{NL}}$, which, as in the local setting, is minus the nonlocal divergence. For $p \in \R^{M \times M \times 2}$, the discrete nonlocal divergence is
\begin{equation}
    \left(\mathrm{div}_{\mathrm{NL}}(p)\right)_i := \sum_{j=1}^M \left(\mathrm{div}\, \Phi_{\cdot,j,\cdot}\right)_i,
\end{equation}
where, for each fixed $j$, the local divergence acts on the vector field $i \mapsto \Phi_{i,j,\cdot}$, and
\begin{equation}
    \Phi_{i,j,t} := p_{j,i,t} - p_{i,j,t}
\end{equation}
for $i,j \in \{1,\dots,M\}$ and $t \in \{1,2\}$.

With the discrete setting in place, we now derive an algorithm for computing solutions to problem \eqref{NLChamLionVar}. Specifically, we use the primal-dual algorithm of Chambolle and Pock \cite{chambolle2011first}. We distinguish two cases.

\medskip
\noindent\textbf{The case $A=\mathrm{Id}$.} In this case, \eqref{NLChamLionDiscrete} is a nonlocal version of the ROF model \cite{Rudin-Osher-Fatemi}, and the primal-dual formulation is obtained in a straightforward way; for details, we refer the reader to \cite{chambolle2011first, chambolle2016introduction}. Let us recall the main steps. We define
\begin{equation}
    B_\omega^\infty(0) := \left\{ p \in \R^{M\times M\times 2} : |p_{i,j,\cdot}| \le \omega_{i,j} \ \text{ for all } i,j \in \{1,\dots,M\} \right\}
\end{equation}
and denote by $\delta_{B_\omega^\infty(0)}$ its indicator function, that is,
\begin{equation}
    \delta_{B_\omega^{\infty}(0)} (p) = \begin{cases}
        0 & \text{ if } p \in B_\omega^{\infty}(0),\\
        +\infty & \text{ if } p \not \in B_\omega^{\infty}(0).
    \end{cases}
\end{equation}
Since each term in $\|\cdot\|_{\omega,1}$ is the Euclidean norm multiplied by $\omega_{i,j}$, whose convex conjugates are the indicators of the closed discs of radius $\omega_{i,j}$, we have
\begin{equation}
    \|l\|_{\omega,1} = \sup_{p} \ \langle l, p\rangle - \delta_{B_\omega^\infty(0)}(p),
\end{equation}
for $l \in \R^{M\times M\times 2}$, where $p \in \R^{M\times M\times 2}$ is the dual variable. Taking $l = \nabla_{\mathrm{NL}} u$ and using that $\mathrm{div}_{\mathrm{NL}}$ is minus the adjoint of $\nabla_{\mathrm{NL}}$, we obtain
\begin{equation}
    \|\nabla_{\mathrm{NL}} u\|_{\omega,1} = \sup_p \ - \langle u , \mathrm{div}_{\mathrm{NL}}(p) \rangle - \delta_{B_\omega^{\infty}(0)} (p).
\end{equation}
Hence, for $A=\mathrm{Id}$ the primal-dual formulation of \eqref{NLChamLionDiscrete} is
\begin{equation} \label{IdentityPrimalDualFormulation}
    \min_u \max_p \ - \langle u, \mathrm{div}_{\mathrm{NL}}(p)\rangle - \delta_{B_\omega^{\infty}(0)} (p) + \frac{\lambda}{2} \|u - u_0\|_2^2.
\end{equation}
After the corresponding calculations, we obtain the following Chambolle--Pock algorithm. Let $\sigma, \tau > 0$, set $u^{(0)} = u_0 \in \R^{M}$, $p^{(0)} = 0 \in \R^{M \times M \times 2}$ and $\hat u^{(0)} = u^{(0)}$, and iterate for $n \ge 0$
\begin{equation} \label{NLPDIdentityAlgorithm}
\left\lbrace\begin{array}{l}
    p^{(n+1)}_{i,j,t} = \dfrac{p^{(n)}_{i,j,t} + \sigma \nabla_{\mathrm{NL}} \hat u_{i,j,t}^{(n)}}{\max\left\{1, | p^{(n)}_{i,j,\cdot} + \sigma \nabla_{\mathrm{NL}} \hat u_{i,j,\cdot}^{(n)} | / \omega_{i,j}\right\}},\\[3ex]
    u^{(n+1)}_{i} = \dfrac{u^{(n)}_i+\tau \left( \left(\mathrm{div}_{\mathrm{NL}}(p^{(n+1)})\right)_{i}+ \lambda (u_0)_{i}\right)}{1+\tau\lambda},\\[2ex]
    \hat u^{(n+1)}_{i} = 2 u^{(n+1)}_{i} - u^{(n)}_{i},
\end{array}\right.
\end{equation}
with the convention that $p_{i,j\cdot}^{(n+1)}=0$ whenever $\omega_{i,j} = 0.$ As usual, the update in $p$ is the pointwise Euclidean projection of $p^{(n)}_{i,j,\cdot} + \sigma \nabla_{\mathrm{NL}} \hat u_{i,j,\cdot}^{(n)}$ onto the closed disc of radius $\omega_{i,j}$. The update in $u$ is a proximal step for the data term $\frac{\lambda}{2}\|u-u_0\|_2^2$, taken after an explicit gradient step on the coupling term $-\langle u,\mathrm{div}_\mathrm{NL} (p)\rangle$.

\medskip
\noindent\textbf{The general case.} For a general operator $A$, the scheme above is not directly applicable, since evaluating the proximal operator of the fidelity term $u \mapsto \frac{\lambda}{2}\|Au-u_0\|_2^2$ requires solving a linear system with a certain matrix depending on $A$, which may be costly. To avoid this, we also ``dualize'' the fidelity term by introducing a second dual variable $q \in \R^{M_0}$. For $v \in \R^{M_0}$, we define
\begin{equation}
    h(v) := \frac{\lambda}{2}\|v-u_0\|_2^2.
\end{equation}
Its convex conjugate is
\begin{equation*}
    h^*(q) = \langle q, u_0\rangle + \frac{\|q\|_2^2}{2\lambda},
\end{equation*}
and, since $h = h^{**}$, we have
\begin{equation}
    \frac{\lambda}{2}\|Au-u_0\|_2^2 = \sup_{q } \ \langle Au, q\rangle - \langle q, u_0\rangle - \frac{\|q\|_2^2}{2\lambda}.
\end{equation}
With this, and arguing as in the case $A=\mathrm{Id}$, the primal-dual formulation of \eqref{NLChamLionDiscrete} is
\begin{equation} \label{qPrimalDualFormulation}
    \min_u \max_{p,q} \ \langle u, -\mathrm{div}_{\mathrm{NL}}(p) + A^* q\rangle - \delta_{B_\omega^{\infty}(0)}(p) - \langle q, u_0\rangle - \frac{\|q\|_2^2}{2\lambda},
\end{equation}
where $p \in \R^{M\times M\times 2}$, $q \in \R^{M_0}$, and $A^* \colon \R^{M_0} \to \R^M$ is the adjoint of $A$. This formulation leads to the following Chambolle--Pock algorithm. Let $\sigma,\tau>0$, choose $u^{(0)} \in \R^M$ (e.g., $u_0$ when $M_0 = M$), set $p^{(0)} = 0 \in \R^{M\times M\times 2}$, $q^{(0)} = 0 \in \R^{M_0}$ and $\hat u^{(0)} = u^{(0)}$, and iterate for $n \ge 0$
\begin{equation} \label{NLPDqAlgorithm}
\left\{\begin{array}{l}
    p^{(n+1)}_{i,j,t} = \dfrac{p^{(n)}_{i,j,t} + \sigma \nabla_{\mathrm{NL}} \hat u^{(n)}_{i,j,t}}{\max\left\{1,\ |p^{(n)}_{i,j,\cdot} + \sigma \nabla_{\mathrm{NL}} \hat u^{(n)}_{i,j,\cdot}| / \omega_{i,j}\right\}}, \\[3ex]
    q^{(n+1)}_i = \dfrac{\lambda \left(q^{(n)}_i + \sigma \left((A\hat u^{(n)})_i - (u_0)_i\right)\right)}{\lambda + \sigma}, \\[3ex]
    u^{(n+1)}_i = u^{(n)}_i - \tau \left((A^* q^{(n+1)})_i - \left(\mathrm{div}_{\mathrm{NL}}\big(p^{(n+1)}\big)\right)_i\right), \\[1.5ex]
    \hat u^{(n+1)}_i = 2 u^{(n+1)}_i - u^{(n)}_i,
\end{array}\right.
\end{equation}
with the same convention as before that $p^{(n+1)}_{i,j,\cdot} = 0$ whenever $\omega_{i,j} = 0$. The update in $p$ is the same as in \eqref{NLPDIdentityAlgorithm}, and the update in $q$ is the proximal map of $\sigma h^*$. Since the primal variable $u$ appears in \eqref{qPrimalDualFormulation} only through the coupling term, the update in $u$ is a gradient step on this term, with the updated dual variables $(p^{(n+1)}, q^{(n+1)})$, followed by the proximal map of the null function, i.e., the identity.

\section{Experimental Results} \label{ExpResSec}

\begin{figure}[t]
    \centering
    \begin{subfigure}[t]{0.18\linewidth}
        \centering
        \includegraphics[width=\linewidth]{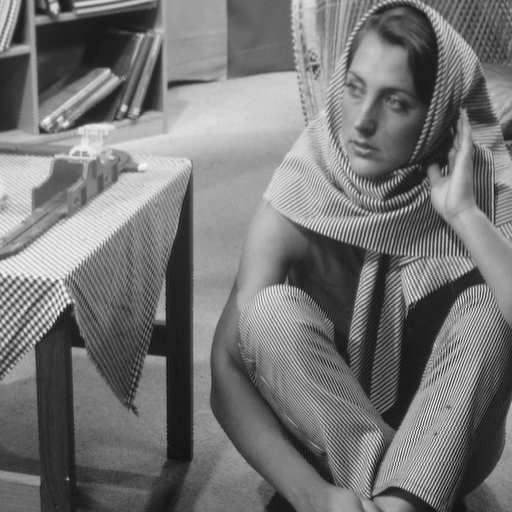}
        \caption*{{Barbara}}
    \end{subfigure}\hfill
    \begin{subfigure}[t]{0.18\linewidth}
        \centering
        \includegraphics[width=\linewidth]{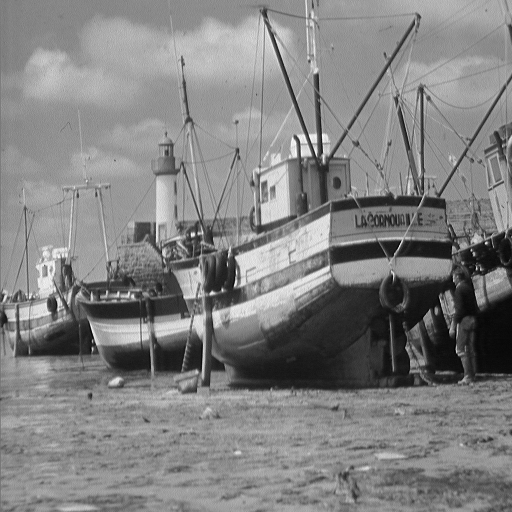}
        \caption*{{Boat}}
    \end{subfigure}\hfill
    \begin{subfigure}[t]{0.175\linewidth}
        \centering
        \includegraphics[width=\linewidth]{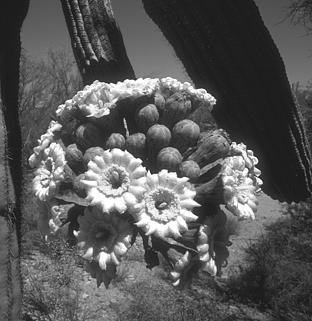}
        \caption*{{Cactus}}
    \end{subfigure}\hfill
    \begin{subfigure}[t]{0.18\linewidth}
        \centering
        \includegraphics[width=\linewidth]{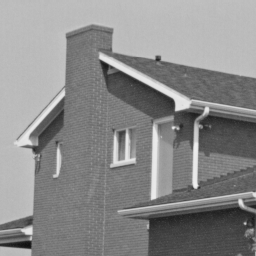}
        \caption*{{House}}
    \end{subfigure}\hfill
    \begin{subfigure}[t]{0.18\linewidth}
        \centering
        \includegraphics[width=\linewidth]{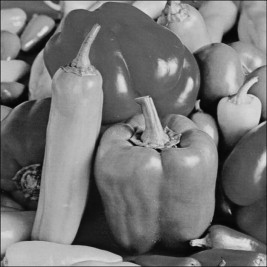}
        \caption*{{Peppers}}
    \end{subfigure}
    \caption{Testing images used in the experiments presented in Section \ref{ExpResSec}.}
    \label{fig:TestImages}
\end{figure}

In this section, we evaluate the performance of the $\mathrm{NLTV}_\omega$ regularizer in four representative image restoration problems: denoising, deblurring, super-resolution and inpainting. The experiments are conducted on five benchmark grayscale test images, shown in Figure \ref{fig:TestImages}.

\begin{table}[t]
\centering
\begin{tabular}{lcc}
\hline
\multirow{2}{*}{\textbf{Model}}  & \multicolumn{2}{c}{\textbf{FLOPs/iter}} \\ 
  & General $\rho$ & $\rho = 15$ \\ \hline
Tikhonov & $13 M$ & $13 M$ \\
$\mathrm{TV}$ & $18M$ & $18M$ \\
Gil.  $\mathrm{NLTV}$ & $(17\rho^2 + 2)M$ & $3827 M$ \\
$\mathrm{TGV}^2$ & $64 M$ & $64 M$ \\
$\mathrm{NLTGV}^2$ & $(63 \rho^2 + 14)M$ & $14189 M$ \\
$\mathrm{NLTV}_\omega$ & $(25 \rho^2 + 4)M$ & $5629 M$ \\
\end{tabular}
\caption{Computational costs of the iteration of the primal-dual algorithm of each model in terms of FLOPs. The model-independent cost $C_{\mathrm{base}}$ in \eqref{CbaseDef} has been subtracted from all reported values. The second column gives the computational cost for a general window size $\rho$, while the third gives the cost for $\rho=15$, which is the value used in our experiments.
}
\label{tab:CompCostComp}
\end{table}

\begin{figure*}[ht]
    \centering
    \begin{subfigure}[t]{0.22\textwidth}
        \centering
        \spybarbara{TestImages/barbara.png}
        \caption*{\centering Ground Truth\\
        PSNR / SSIM}
    \end{subfigure}
    \begin{subfigure}[t]{0.22\textwidth}
        \centering
        \spybarbara{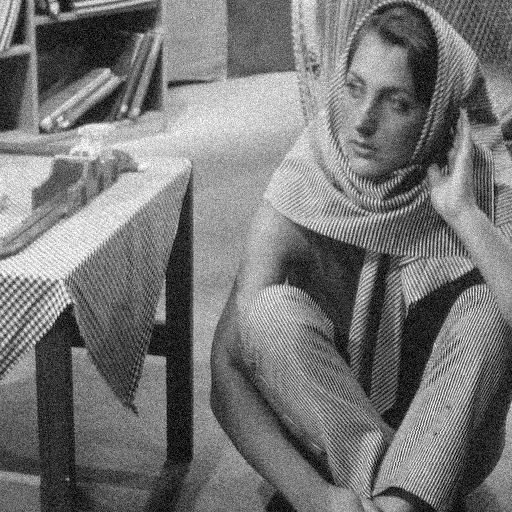}
        \caption*{\centering Noisy Image\\
        24.63 / 0.5772}
    \end{subfigure}
    \begin{subfigure}[t]{0.22\textwidth}
        \centering
        \spybarbara{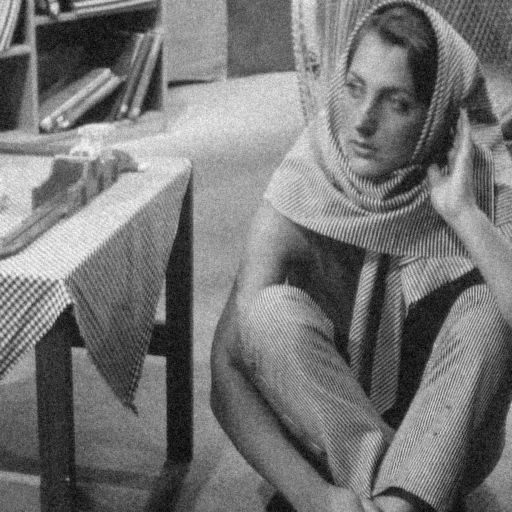}
        \caption*{\centering Tikhonov\\
        27.34 / 0.7203}
    \end{subfigure}
    \begin{subfigure}[t]{0.22\textwidth}
        \centering
        \spybarbara{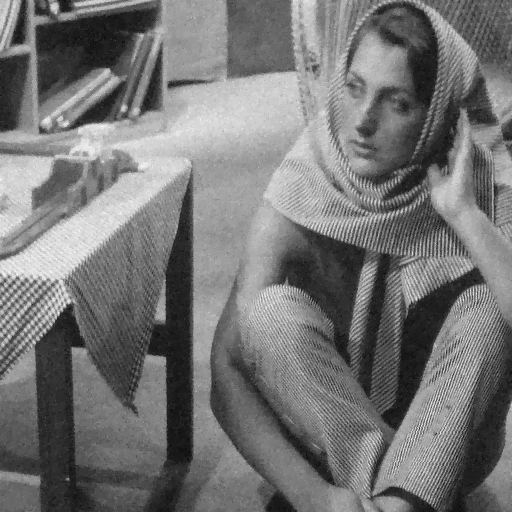}
        \caption*{\centering TV\\
        28.53 / 0.8095}
    \end{subfigure}

    \medskip
    
    \begin{subfigure}[t]{0.22\textwidth}
        \centering
        \spybarbara{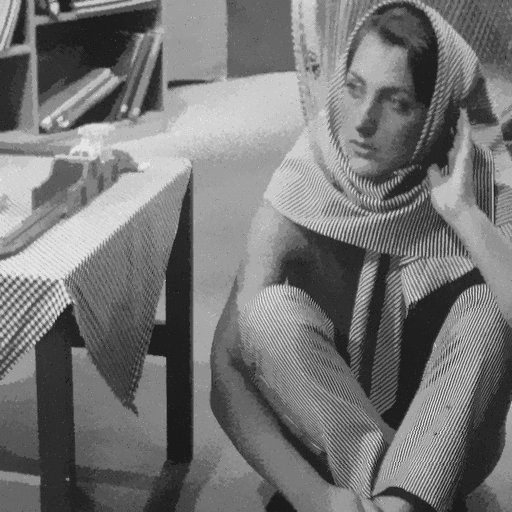}
        \caption*{\centering Gil. $\mathrm{NLTV}$\\
        30.05 / 0.8411}
    \end{subfigure}
    \begin{subfigure}[t]{0.22\textwidth}
        \centering
        \spybarbara{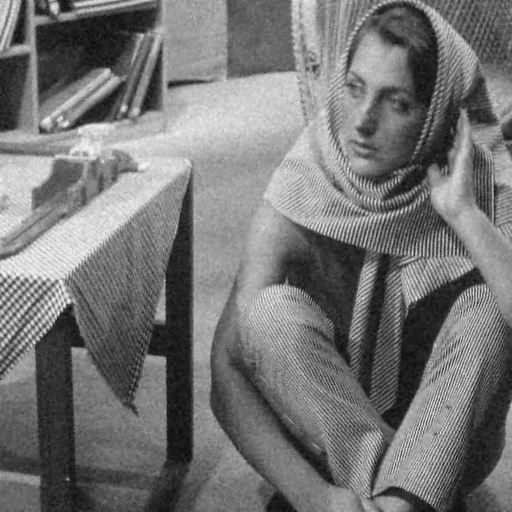}
        \caption*{\centering $\mathrm{TGV}^2$\\
        28.90 / 0.8134}
    \end{subfigure}
    \begin{subfigure}[t]{0.22\textwidth}
        \centering
        \spybarbara{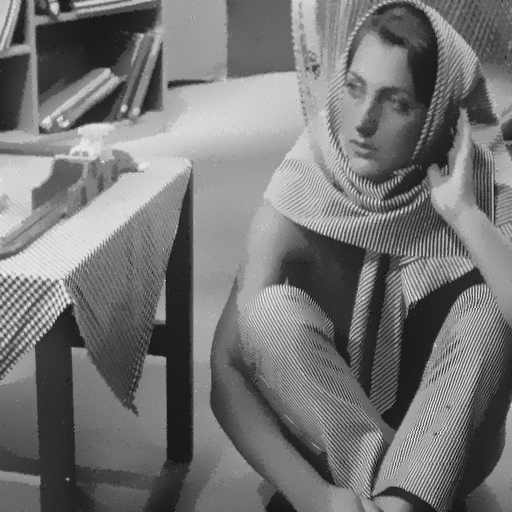}
        \caption*{\centering $\mathrm{NLTGV}^2$\\
        30.13 / 0.8569}
    \end{subfigure}
    \begin{subfigure}[t]{0.22\textwidth}
        \centering
        \spybarbara{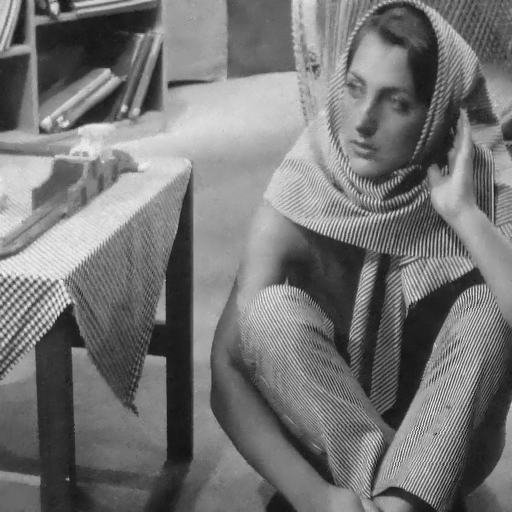}
        \caption*{\centering $\mathrm{NLTV}_\omega$\\
        \textbf{31.21} / \textbf{0.8882}}
    \end{subfigure}

    \caption{Visual comparison of denoising methods on the Barbara test image for $\sigma_{\mathrm{noise}}=15$. As observed in the zoomed-in crop, the local methods struggle to effectively remove noise, while Gilboa and Osher's $\mathrm{NLTV}$ produces undesirable grainy artifacts and $\mathrm{NLTGV}^2$ oversmooths texture and fine detail. Overall, $\mathrm{NLTV}_\omega$ achieves the best reconstruction, both visually and quantitatively.}
    \label{fig:QualCompBarbaraDen15}
\end{figure*}

We compare $\mathrm{NLTV}_\omega$ with the following image regularizers: Tikhonov regularization, i.e., the squared $L^2$ norm of the gradient (Tikhonov) \cite{tikhonov1977solutions}, Total Variation ($\mathrm{TV}$) \cite{Rudin-Osher-Fatemi}, Gilboa and Osher's anisotropic nonlocal total variation (Gil.  $\mathrm{NLTV}$) \eqref{GilboaSOTA} \cite{Gilboa-Osher}, second-order Total Generalized Variation ($\mathrm{TGV}^2$) \eqref{TGVDef} \cite{bredies2010total} and nonlocal Total Generalized Variation ($\mathrm{NLTGV}^2$) \eqref{NLTGVDef} \cite{ranftl2014non}. For each regularizer, we formulate a minimization problem analogous to \eqref{NLChamLionVar}, replacing $\mathrm{NLTV}_\omega$ with the corresponding regularization functional. We then derive and implement the associated primal-dual algorithm, following the approach described in Section \ref{NumAlgSec}. The model parameters are tuned separately for each image to maximize the PSNR.

Throughout this section, we identify flattened indexes $i \in \{1,...,M\}$ with two-dimensional indexes $(i_1,i_2) \in \mathcal{I}$, with $\mathcal{I} = \{1,...,N_1\} \times \{1,...,N_2\}$, where $N_1$ and $N_2$ denote the height and width of the image, respectively, so that $M=N_1N_2$. Hence, a grayscale image now lives in $\R^\mathcal{I}$.

For all nonlocal methods we use the NLMeans weights \cite{Buades-Coll-Morel,Buades-Coll-Morel-bis}, which for an initial image
$u_0 \in \mathbb{R}^{\mathcal{I}}$ are given by
\begin{equation} \label{NLMeansWeightDefDiscrete}
    \omega_{i,j} = \frac{1}{\Gamma_{i}} \exp\left( - \frac{d\left((u_0)_{i}, (u_0)_{j}\right)}{h^2} \right) 
\end{equation}
for $i=(i_1,i_2), j=(j_1,j_2) \in \mathcal{I}$ such that $|i_1 - j_1|\le \left\lfloor \rho / 2 \right \rfloor$ and $|i_2 - j_2| \le \left\lfloor \rho / 2 \right \rfloor$, and $0$ otherwise, with $\rho\in \Z_{>0}$ being the search window size. Here, $\Gamma_{i}$ is a normalization constant ensuring that $\sum_{j} \omega_{i,j} = 1$ for all $i \in \mathcal{I}$, $h>0$ is a filtering parameter that controls the sensitivity of the weight to patch dissimilarity, and $d$ measures the squared Euclidean distance between patches of size $(2\nu + 1) \times (2\nu + 1)$, $\nu \in \Z_{>0}$, that is,
\begin{equation}
    d\left((u_0)_{i}, (u_0)_{j}\right) = d\left((u_0)_{(i_1,i_2)}, (u_0)_{(j_1,j_2)}\right) = \sum_{z_1 = -\nu}^\nu \sum_{z_2 = -\nu}^\nu |(u_0)_{i_1+z_1,i_2+z_2} - (u_0)_{j_1+z_1,j_2+z_2}|^2,
\end{equation}
for all $i, j \in \mathcal{I}$, with out-of-bound indices treated by reflection. To prevent excessive influence of the reference pixel, we set $\omega_{i,i}$, before normalization, as the maximum of the weights within the search window for $j \neq i$. In Appendix
\ref{app:NLMeansWeights} we formulate a continuous version of \eqref{NLMeansWeightDefDiscrete} and show that it is consistent with the theoretical setting, namely that it is an admissible weight function satisfying \eqref{H0}.

The window and patch sizes were fixed to $\rho=15$ and $\nu=3$ for every nonlocal regularizer, while $h$ was optimized together with the remaining trade-off parameters $\lambda$, or $\alpha_0$ and $\alpha_1$. For $\mathrm{NLTGV}^2$, the weight functions $\omega_0$ and $\omega_1$ were chosen as NLMeans weights with the same $\rho$, $\nu$ and $h$, multiplied by two constants $\alpha_0$ and $\alpha_1$ which were also optimized. The primal-dual step sizes $\sigma$ and $\tau$ were chosen so as to satisfy the convergence criterion of \cite{chambolle2011first}.

Quantitative evaluation is reported in terms of the PSNR and the structural similarity index measure (SSIM). The former measures pixelwise error directly, while the latter estimates reconstruction quality in terms of structure preservation.

\subsection{Computational Costs}

In this subsection we discuss the computational costs of all restoration models. While they are induced by image regularizers, directly comparing the computational cost of the regularizers themselves is not meaningful. Instead, we compare the computational cost of one iteration of the primal-dual algorithm associated with each of them, as this is the algorithm we actually use to obtain each result. Recall that $M$ and $M_0$ denote the number of pixels in the restored and corrupted images, respectively, while $\rho$ denotes the window size used to compute the nonlocal weights.
When $A \neq \mathrm{Id}$, the $q$-update appearing in \eqref{NLPDqAlgorithm} is the same for all models and entails a computational cost of $C_A + 4M_0$, where $C_A$ denotes the computational cost associated with the operator $A$. The $u$-update in \eqref{NLPDqAlgorithm} is also similar for all models apart from the computation associated with the regularizer. Its model-independent part has computational cost  $3M + C_{A^*} $, where $C_{A^*}$ denotes the computational cost associated with $A^*$, the adjoint of $A$. 
Instead, when $A=\mathrm{Id}$, the $q$ update is no longer necessary, and the $u$-update takes a different form \eqref{NLPDIdentityAlgorithm}, which is again similar for all models and whose model-independent part entails a computational cost of $5M$.

Thus, the computational cost of each method can be decomposed into a model-independent baseline cost and a model-dependent contribution. We define the former as

\begin{equation} \label{CbaseDef}
    C_{\mathrm{base}} = \begin{cases}
        4M_0 + 3M + C_A + C_{A^*} & \mbox{ if } A \neq \mathrm{Id},\\
        5M & \mbox{ if } A = \mathrm{Id}.
    \end{cases}
\end{equation}

We omit $C_{\mathrm{base}}$ from the subsequent comparison for clarity, as it is common to all techniques considered in the respective setting. The resulting (model-dependent) computational cost of one iteration of the primal-dual algorithm for each model, expressed in floating-point operations (FLOPs), is reported in Table \ref{tab:CompCostComp}. Iteration cost for general window size $\rho$ is shown in the second column, and iteration cost for $\rho = 15$ is shown in the third, as this is the choice we have used in our experiments.

As expected, even for moderate values of $\rho$, the nonlocal methods entail a substantially higher computational cost than their local counterparts, as their cost scales quadratically with the window size. Gilboa and Osher's anisotropic $\mathrm{NLTV}$ is the least expensive of the nonlocal methods, requiring approximately $8\rho^2 M$ fewer FLOPs per iteration than the next most expensive method, namely, our proposed $\mathrm{NLTV}_\omega$ model. Finally, $\mathrm{NLTGV}^2$ is by far the most computationally expensive, with its leading-order $\rho^2 M$ cost being approximately 2.5 times larger than ours. We note that in the original $\mathrm{NLTGV}^2$ article \cite{ranftl2014non}, smaller $5\times 5$ search windows were used, but we have fixed $\rho = 15$ in order to obtain a fair comparison among nonlocal methods.

\begin{figure*}[ht]
    \centering
    \begin{subfigure}[t]{0.22\textwidth}
        \centering
        \spypeppers{TestImages/peppers.png}
        \caption*{\centering Ground Truth\\
        PSNR / SSIM}
    \end{subfigure}
    \begin{subfigure}[t]{0.22\textwidth}
        \centering
        \spypeppers{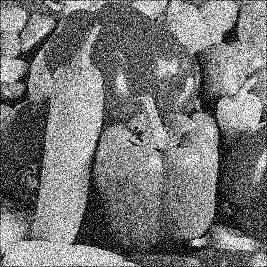}
        \caption*{\centering Noisy Image\\
        14.69 / 0.1804}
    \end{subfigure}
    \begin{subfigure}[t]{0.22\textwidth}
        \centering
        \spypeppers{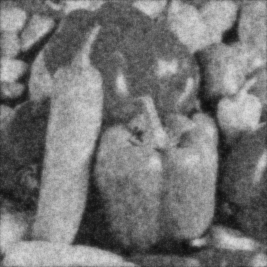}
        \caption*{\centering Tikhonov\\
        24.07 / 0.6114}
    \end{subfigure}
    \begin{subfigure}[t]{0.22\textwidth}
        \centering
        \spypeppers{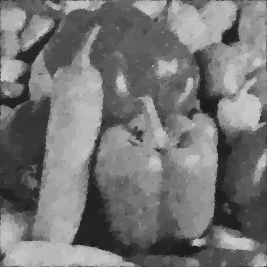}
        \caption*{\centering TV\\
        25.04 / 0.7292}
    \end{subfigure}

    \medskip
    
    \begin{subfigure}[t]{0.22\textwidth}
        \centering
        \spypeppers{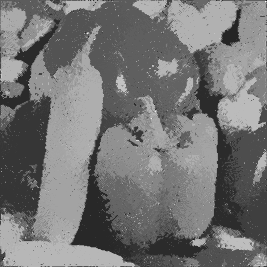}
        \caption*{\centering Gil. $\mathrm{NLTV}$\\
        22.70 / 0.5994}
    \end{subfigure}
    \begin{subfigure}[t]{0.22\textwidth}
        \centering
        \spypeppers{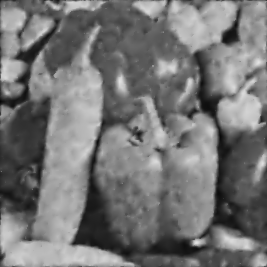}
        \caption*{\centering $\mathrm{TGV}^2$\\
        25.53 / \textbf{0.7507}}
    \end{subfigure}
    \begin{subfigure}[t]{0.22\textwidth}
        \centering
        \spypeppers{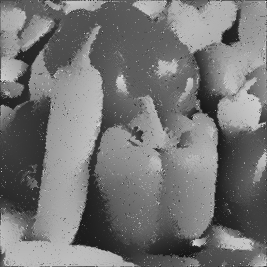}
        \caption*{\centering $\mathrm{NLTGV}^2$\\
        23.72 / 0.6378}
    \end{subfigure}
    \begin{subfigure}[t]{0.22\textwidth}
        \centering
        \spypeppers{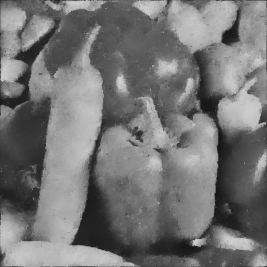}
        \caption*{\centering $\mathrm{NLTV}_\omega$\\
        \textbf{25.59} / 0.7487}
    \end{subfigure}

    \caption{Denoising results for the Peppers test image corrupted by Gaussian noise with standard deviation $\sigma_{\mathrm{noise}}=50$. oves excessive texture and introduces salt-and-pepper-like artifacts. Both $\mathrm{TV}$ and Gilboa and Osher's $\mathrm{NLTV}$ produce undesirable piecewise constant reconstructions, whereas $\mathrm{TGV}^2$ retains a considerable amount of residual noise. By contrast, the proposed $\mathrm{NLTV}_\omega$ regularizer yields the most natural-looking reconstruction, as particularly evident in the enlarged crop.}
    \label{fig:QualCompPeppersDen50}
\end{figure*}

\begin{figure*}[h]
    \centering
    \begin{subfigure}[t]{0.22\textwidth}
        \centering
        \spycactus{TestImages/cactus.png}
        \caption*{\centering Ground Truth\\
        PSNR / SSIM}
    \end{subfigure}
    \begin{subfigure}[t]{0.22\textwidth}
        \centering
        \spycactus{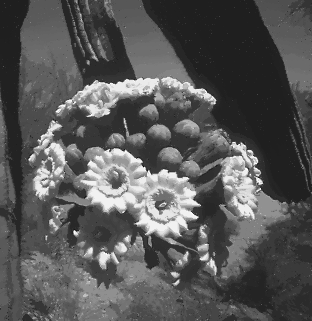}
        \caption*{\centering Gil. $\mathrm{NLTV}$\\
        w. GT weights\\
        28.65 / 0.8490}
    \end{subfigure}
    \begin{subfigure}[t]{0.22\textwidth}
        \centering
        \spycactus{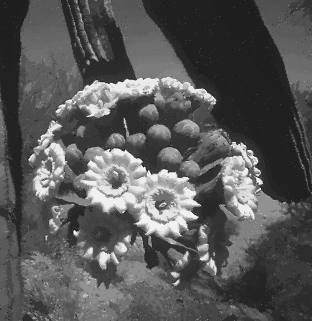}
        \caption*{\centering $\mathrm{NLTGV}^2$\\
        w. GT weights\\
        28.67 / 0.8463}
    \end{subfigure}
    \begin{subfigure}[t]{0.22\textwidth}
        \centering
        \spycactus{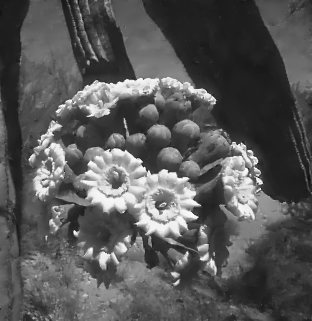}
        \caption*{\centering $\mathrm{NLTV}_\omega$\\
        w. GT weights\\
        \textbf{28.92} / \textbf{0.8604}}
    \end{subfigure}
    \medskip
    \begin{subfigure}[t]{0.22\textwidth}
        \centering
        \spycactus{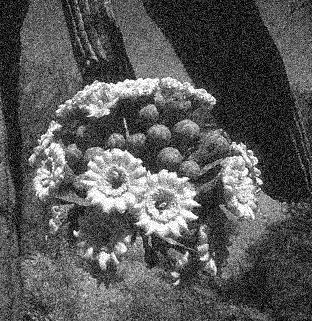}
        \caption*{\centering Noisy Image\\
        20.42 / 0.4287}
    \end{subfigure}
    \begin{subfigure}[t]{0.22\textwidth}
        \centering
        \spycactus{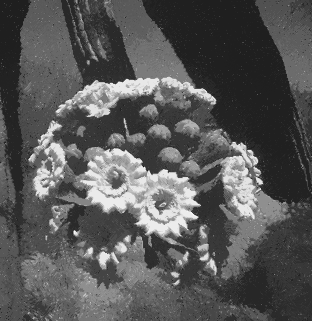}
        \caption*{\centering Gil. $\mathrm{NLTV}$\\
        w. Noisy weights\\
        26.25 / 0.7462}
    \end{subfigure}
    \begin{subfigure}[t]{0.22\textwidth}
        \centering
        \spycactus{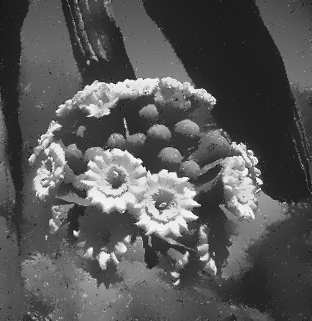}
        \caption*{\centering $\mathrm{NLTGV}^2$\\ 
        w. Noisy weights\\
        26.20 / 0.7511}
    \end{subfigure}
    \begin{subfigure}[t]{0.22\textwidth}
        \centering
        \spycactus{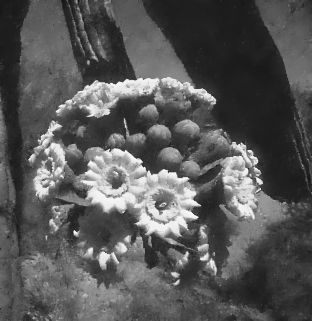}
        \caption*{\centering $\mathrm{NLTV}_\omega$\\
        w. Noisy weights\\
        \textbf{27.28} / \textbf{0.8000}}
    \end{subfigure}

    \caption{Denoising results for the Cactus test image corrupted by Gaussian noise with standard deviation $\sigma_{\mathrm{noise}}=25$, comparing the performance of the nonlocal regularizers under ideal and realistic weight scenarios. The top row shows the ground truth and the denoising results obtained using weights computed from the ground-truth image. The bottom row shows the noisy input image and the corresponding results obtained using weights computed directly from that image. In the latter scenario, Gilboa and Osher's $\mathrm{NLTV}$ leaves substantial residual noise, whereas $\mathrm{NLTGV}^2$ removes excessive detail. By contrast, the proposed $\mathrm{NLTV}_\omega$ regularizer achieves a better balance between noise suppression and texture preservation. When the weights are computed from the ground-truth image, all three nonlocal methods produce substantially better results. Nevertheless, $\mathrm{NLTV}_\omega$ still provides the best texture preservation, as particularly evident in the enlarged crop. Overall, $\mathrm{NLTV}_\omega$ is the model that is least affected by weights derived in non-ideal scenarios.}
    \label{fig:QualCompCactusDen25GTw}
\end{figure*}

\subsection{Image Denoising}

We start the experimental comparison with the task of image denoising, which remains one of the most representative image restoration problems. In our setting, this task corresponds to the $A=\mathrm{Id}$ case, which clearly satisfies \eqref{H1} and \eqref{H2'}. Moreover, we may use the simplified primal-dual algorithm \eqref{NLPDIdentityAlgorithm} that does not dualize the fidelity term. For the other regularizers, we also use similarly simplified algorithms. The noisy observations are generated by adding zero-mean Gaussian noise with standard deviations $\sigma_{\mathrm{noise}}\in\{15, 25,50\}$.

Table \ref{tab:QuantCompDen} shows the average metrics over the testing set for all methods at the three considered noise levels. $\mathrm{NLTV}_\omega$ regularization achieves the best performance in terms of PSNR and SSIM across all noise levels. 
The runner-up is $\mathrm{TGV}^2$, in both metrics for $\sigma_{\mathrm{noise}}\in\{25, 50\}$ and in SSIM for $\sigma_{\mathrm{noise}}= 15$, closely followed by $\mathrm{TV}$, while $\mathrm{NLTGV}^2$ takes second place in PSNR for $\sigma_{\mathrm{noise}}= 15$. As the noise level increases, the performance of Gilboa and Osher's $\mathrm{NLTV}$ and $\mathrm{NLTGV}^2$ drops substantially, even falling behind Tikhonov, which is often considered a lower baseline, in terms of PSNR for $\sigma_{\mathrm{noise}}= 50$.

\begin{table}[t]
\centering
\begin{tabular}{lcccccc}
\hline
& \multicolumn{2}{c}{$\sigma_{\mathrm{noise}}=15$}
& \multicolumn{2}{c}{$\sigma_{\mathrm{noise}}=25$}
& \multicolumn{2}{c}{$\sigma_{\mathrm{noise}}=50$} \\
\cline{2-7}
\textbf{Model}
& \textbf{PSNR $\uparrow$} & \textbf{SSIM $\uparrow$}
& \textbf{PSNR $\uparrow$} & \textbf{SSIM $\uparrow$}
& \textbf{PSNR $\uparrow$} & \textbf{SSIM $\uparrow$} \\
\hline
Noisy
& 24.63 & 0.5367
& 20.29 & 0.3616
& 14.72 & 0.1795 \\

Tikhonov
& 28.83 & 0.7417
& 26.55 & 0.6581
& 23.85 & 0.5524 \\

$\mathrm{TV}$
& 30.49 & 0.8429
& 27.98 & 0.7780
& 24.76 & 0.6764 \\

Gil. $\mathrm{NLTV}$
& 30.54 & 0.8314
& 27.58 & 0.7496
& 23.19 & 0.5968 \\

$\mathrm{TGV}^2$
& 30.74 & \underline{0.8460}
& \underline{28.16} & \underline{0.7817}
& \underline{24.95} & \underline{0.6848} \\

$\mathrm{NLTGV}^2$
& \underline{30.85} & 0.8430
& 27.97 & 0.7664
& 23.75 & 0.6114 \\

$\mathrm{NLTV}_\omega$
& \textbf{31.91} & \textbf{0.8768}
& \textbf{29.21} & \textbf{0.8291}
& \textbf{25.18} & \textbf{0.7006} \\
\hline
\end{tabular}
\caption{Quantitative comparison of denoising methods on the testing set from Figure \ref{fig:TestImages} for different noise levels. Best results are highlighted in \textbf{bold} and second-best results are \underline{underlined}. $\mathrm{NLTV}_\omega$ regularization outperforms the other models at every  noise level and for both metrics.}
\label{tab:QuantCompDen}
\end{table}

Figure~\ref{fig:QualCompBarbaraDen15} shows the Barbara test image corrupted by Gaussian noise with standard deviation  $\sigma_{\mathrm{noise}}=15$, together with the denoising results obtained using each of the considered regularizers. The local methods struggle to effectively remove the noise, as is particularly evident in the enlarged crop. Gilboa and Osher's $\mathrm{NLTV}$ leaves noticeable residual grainy noise, whereas $\mathrm{NLTGV}^2$ tends to oversmooth the image, removing fine textures and shadow details. Overall, the proposed $\mathrm{NLTV}_\omega$ regularizer achieves the best balance between noise removal and texture recovery.

At the most challenging noise level, $\sigma_{\mathrm{noise}}=50$, the limitations of each regularizer become apparent in Figure~\ref{fig:QualCompPeppersDen50}, where the corrupted Peppers test image is displayed together with the corresponding reconstructions. Tikhonov regularization produces an overly smooth result while failing to remove the noise completely, whereas $\mathrm{TV}$ introduces pronounced staircasing artifacts. Gilboa and Osher's $\mathrm{NLTV}$ also yields an undesirable piecewise constant reconstruction, with unnatural, grainy texture. Similarly, $\mathrm{NLTGV}^2$ removes excessive texture and leaves high-frequency residual noise resembling a salt-and-pepper pattern. Among the methods discussed so far, $\mathrm{TGV}^2$ preserves texture more effectively, as reflected in the highest SSIM score, although noise is still noticeable. Finally, the proposed $\mathrm{NLTV}_\omega$ regularizer produces the most natural-looking reconstruction, suppressing noise while preserving fine details such as shadows and reflections, as the enlarged crop shows. 

The nonlocal regularizers deteriorate substantially as the noise level increases, particularly Gilboa and Osher's $\mathrm{NLTV}$ and $\mathrm{NLTGV}^2$. From $\sigma_{\mathrm{noise}}=15$ to $\sigma_{\mathrm{noise}}=50$, these two lose between $7.10$ and $7.35$ dB in PSNR and $0.23$ in SSIM, whereas the local methods show maximum decreases of $5.79$ dB and $0.19$, and $\mathrm{NLTV}_\omega$ loses about $6.7$ dB and $0.18$, the latter comparable to the local methods. We attribute this to the reduced reliability of the nonlocal weights when the input image is heavily corrupted, and we argue that the relative robustness of $\mathrm{NLTV}_\omega$ stems from the fact that it promotes gradient similarity, thereby allowing the regularization to be guided by information from surrounding pixels rather than relying solely on the weight function. Its advantage does not rest on this robustness alone, as Figure~\ref{fig:QualCompCactusDen25GTw} shows for the Cactus image at $\sigma_{\mathrm{noise}}=25$. When the weights are computed from the noisy input, Gilboa and Osher's $\mathrm{NLTV}$ leaves substantial residual noise and unnatural piecewise constant regions, $\mathrm{NLTGV}^2$ removes more noise at the cost of excessive smoothing and a more pronounced loss of detail, and $\mathrm{NLTV}_\omega$ balances noise removal and texture preservation. Computing the weights from the ground truth improves all three methods, especially the first two, and the nonlocal models then clearly outperform the local ones, the best of which is $\mathrm{TGV}^2$ with a PSNR of $26.72$ dB and $0.7630$ SSIM. What separates them in this ideal setting is the design of the regularizer rather than the reliability of the weights. Gilboa and Osher's $\mathrm{NLTV}$ and $\mathrm{NLTGV}^2$ behave almost identically, both losing patterns and structures and staircasing in homogeneous regions such as the sky, whereas $\mathrm{NLTV}_\omega$ preserves fine textures and shadows and avoids staircasing, at the price of faint gradient-like patterns in some homogeneous regions. These minor artifacts are outweighed by its better structure preservation, as the higher metrics reflect.

\subsection{Image Deblurring}

We now consider the task of image deblurring, for both Gaussian and motion blur. In the theoretical setting, the operator corresponding to the first kind of blur is $A = K_{\mathrm{Gauss}} * \cdot$, where $K_\mathrm{Gauss}$ is a Gaussian kernel of the form
\begin{equation} \label{DiscreteGaussKerDef}
    K_{\mathrm{Gauss}}(z) := \frac{1}{C}\exp\left(-\frac{|z|^2}{2 \sigma_{\mathrm{blur}}^2}\right) \mathcal{X}_{B_\kappa(0)}(z),
\end{equation}
while the operator corresponding to the second kind of blur is $A = K_{\mathrm{Mot}} * \cdot$, where $K_{{\mathrm{Mot}}}:\R^2 \to [0,+\infty)$ is a motion kernel of the form
\begin{equation} \label{DiscreteMotKerDef}
    K_{{\mathrm{Mot}} (z)} := \frac{1}{C} \exp\left(-\frac{|\eta^\perp \cdot z |^2}{2 \sigma_{\mathrm{blur}}^2}\right) \mathcal{X}_{B_\kappa(0)}(z),
\end{equation}
with $\eta \in \mathbb{S}^1$ being the direction of the blur, and $\eta^\perp$ being a unit vector orthogonal to $\eta$. For both kernels \eqref{DiscreteGaussKerDef} and \eqref{DiscreteMotKerDef}, $\sigma_{\mathrm{blur}}>0$ denotes the standard deviation, $\kappa > 0$ is a truncation constant and $C>0$ ensures normalization with respect to $L^1$ norm. Moreover, we treat out of boundary points by reflection in both operators. We note that, in the numerical implementation, the kernels are truncated to square domains rather than balls.

\begin{remark}
    The operators associated with the kernels \eqref{DiscreteGaussKerDef} and \eqref{DiscreteMotKerDef} are continuous and linear from $L^p(\Omega)$ to $L^p(\Omega)$ and satisfy \eqref{H2'}, provided that the truncation constant $\kappa>0$ is small enough. Indeed, they are of the form $A = K * \cdot$, with $K$ a {nonnegative} 
    and even kernel with $L^1$ norm equal to $1$. We also assume that the support is small enough so that the subset $\Omega_K := \{ z \in \Omega : z + \mathrm{supp}(K) \subset \Omega \} \subset \Omega$ has positive measure. Given any affine function $g(x) = a \cdot x + b$ with $a\in\R^2$ and $b\in\R$, for $x \in \Omega_K$ we have

    \begin{align}
        A g (x) &= \int_{\mathrm{supp}(K)} \left(a \cdot (x - y) + b\right) K(y) \mathrm{d}y
        \\
        &= (a \cdot x + b) \int_{\mathrm{supp}(K)} K(y) \mathrm{d}y - a \cdot \int_{\mathrm{supp}(K)} y K(y) \mathrm{d}y = a \cdot x + b, \label{ContBlurrComupation}
    \end{align}
    where in the last equality we have used that $K$ has integral equal to $1$ over its support, $y K(y)$ is an odd vector field and $\mathrm{supp}(K)$ is symmetric with respect to the origin. Hence, if $A g(x) = 0$ in $\Omega$, then $A g(x) = 0$ in $\Omega_K$, or $a \cdot x + b = 0$ in $\Omega_K$ by the equality above, from which $a=0$ and $b=0$ because $\Omega_K$ has positive measure by assumption. 
\end{remark}

For Gaussian blur, a standard deviation of $\sigma_{\mathrm{blur}} = 2.0$ and a kernel size $13 \times 13$ was used (corresponding to $\kappa=6$ in \eqref{DiscreteGaussKerDef}), while, for Motion blur, a standard deviation of $\sigma_{\mathrm{blur}} = 0.5$ and a kernel length of $13$ was used (corresponding to $\kappa=6$ in \eqref{DiscreteMotKerDef}). Additionally, zero-mean Gaussian noise with standard deviation $\sigma_{\mathrm{noise}}=5$ was added to the blurred images to create more realistic corrupted data. 

Table \ref{tab:QuantCompBlur} reports the mean metrics of the Gaussian and motion deblurring results over the testing set for all methods. The proposed $\mathrm{NLTV}_\omega$ regularizer achieves the highest scores for both types of blur, followed by $\mathrm{TGV}^2$. Notably, in these experiments, Gilboa and Osher's $\mathrm{NLTV}$ and $\mathrm{NLTGV}^2$ are outperformed by $\mathrm{TV}$ in both metrics. We attribute this result, once again, to the limited reliability of the nonlocal weights under these degradation conditions.

\begin{table}[t]
\centering
\begin{tabular}{lcccc}
\hline
& \multicolumn{2}{c}{Gaussian}
& \multicolumn{2}{c}{Motion} \\
\cline{2-5}
\textbf{Model}
& \textbf{PSNR $\uparrow$} & \textbf{SSIM $\uparrow$}
& \textbf{PSNR $\uparrow$} & \textbf{SSIM $\uparrow$} \\
\hline
Blurry
& 24.46 & 0.6078
& 22.94 & 0.5497 \\

Tikhonov
& 26.54 & 0.7199
& 25.92 & 0.6756 \\

$\mathrm{TV}$
& 27.03 & 0.7627
& 26.99 & 0.7630 \\

Gil. $\mathrm{NLTV}$
& 26.57 & 0.7411
& 25.76 & 0.6920 \\

$\mathrm{TGV}^2$
& \underline{27.27} & \underline{0.7710}
& \underline{27.10} & \underline{0.7667} \\

$\mathrm{NLTGV}^2$
& 26.94 & 0.7524
& 26.03 & 0.7075 \\

$\mathrm{NLTV}_\omega$
& \textbf{27.62} & \textbf{0.7801}
& \textbf{27.36} & \textbf{0.7719} \\
\hline
\end{tabular}
\caption{Quantitative comparison of deblurring methods on the testing
set from Figure \ref{fig:TestImages} for Gaussian and motion blur.
Best results are highlighted in \textbf{bold} and second-best results are \underline{underlined}. $\mathrm{NLTV}_\omega$ regularization outperforms the other methods for both types of blur.}
\label{tab:QuantCompBlur}
\end{table}

\begin{figure*}[t]
    \centering
    \begin{subfigure}[t]{0.22\textwidth}
        \centering
        \spycactusblur{TestImages/cactus.png}
        \caption*{\centering Ground Truth\\
        PSNR / SSIM}
    \end{subfigure}
    \begin{subfigure}[t]{0.22\textwidth}
        \centering
        \spycactusblur{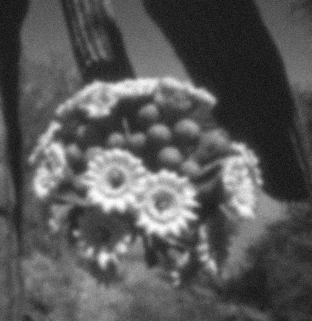}
        \caption*{\centering Blurry\\
        22.73 / 0.5668}
    \end{subfigure}
    \begin{subfigure}[t]{0.22\textwidth}
        \centering
        \spycactusblur{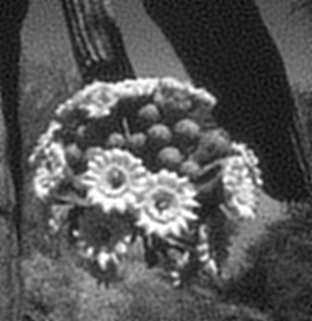}
        \caption*{\centering Tikhonov\\
        24.64 / 0.6717}
    \end{subfigure}
    \begin{subfigure}[t]{0.22\textwidth}
        \centering
        \spycactusblur{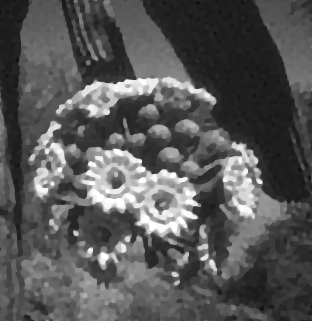}
        \caption*{\centering TV\\
        25.03 / 0.7282}
    \end{subfigure}

    \medskip
    
    \begin{subfigure}[t]{0.22\textwidth}
        \centering
        \spycactusblur{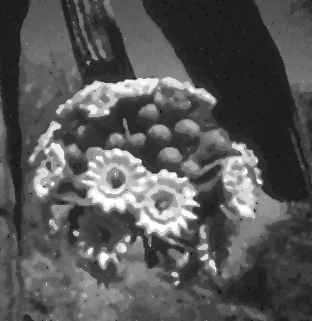}
        \caption*{\centering Gil. $\mathrm{NLTV}$\\
        24.50 / 0.7024}
    \end{subfigure}
    \begin{subfigure}[t]{0.22\textwidth}
        \centering
        \spycactusblur{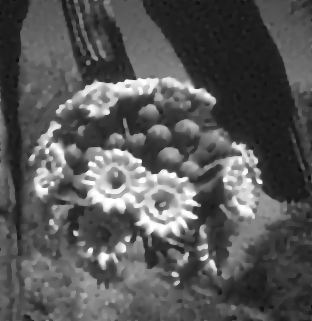}
        \caption*{\centering $\mathrm{TGV}^2$\\
        25.12 / 0.7305}
    \end{subfigure}
    \begin{subfigure}[t]{0.22\textwidth}
        \centering
        \spycactusblur{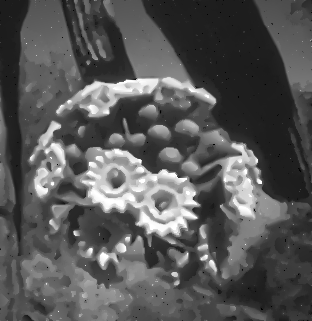}
        \caption*{\centering $\mathrm{NLTGV}^2$\\
        24.52 / 0.7018}
    \end{subfigure}
    \begin{subfigure}[t]{0.22\textwidth}
        \centering
        \spycactusblur{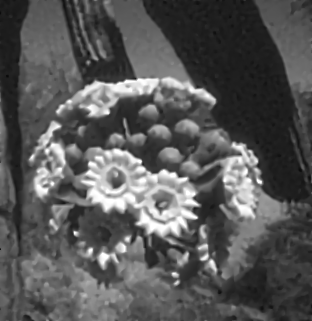}
        \caption*{\centering $\mathrm{NLTV}_\omega$\\
        \textbf{25.20} / \textbf{0.7330}}
    \end{subfigure}

    \caption{Deblurring results for the Cactus test image corrupted by Gaussian blurring induced by the kernel \eqref{DiscreteGaussKerDef} and Gaussian noise with $\sigma_\mathrm{noise}=5$. Tikhonov regularization presents pronounced residual noise, while $\mathrm{TV}$, Gilboa and Osher's $\mathrm{NLTV}$ and $\mathrm{NLTGV}^2$ result in unnatural piecewise constant reconstructions. As particularly evident in the enlarged crop, the proposed $\mathrm{NLTV}_\omega$ regularizer preserves textures more effectively than $\mathrm{TGV}^2$, yielding the most convincing visual result overall.}
    \label{fig:QualCompCactusBlurGauss}
\end{figure*}

Figure~\ref{fig:QualCompCactusBlurGauss} shows the results of Gaussian deblurring on the Cactus test image. The Tikhonov reconstruction retains a substantial amount of residual noise, whereas $\mathrm{TV}$, Gilboa and Osher's $\mathrm{NLTV}$, and $\mathrm{NLTGV}^2$ all produce undesirably piecewise constant reconstructions. By contrast, both $\mathrm{TGV}^2$ and the proposed $\mathrm{NLTV}_\omega$ regularizer effectively suppress noise while better preserving image texture. Their results are visually
close, although the enlarged crop reveals a cartoon-like appearance in the $\mathrm{TGV}^2$ reconstruction, which renders textured regions as flat patches, whereas $\mathrm{NLTV}_\omega$ better retains their fine detail.

\begin{figure*}[ht]
    \centering
    \begin{subfigure}[t]{0.22\textwidth}
        \centering
        \spypeppersblur{TestImages/peppers.png}
        \caption*{\centering Ground Truth\\
        PSNR / SSIM}
    \end{subfigure}
    \begin{subfigure}[t]{0.22\textwidth}
        \centering
        \spypeppersblur{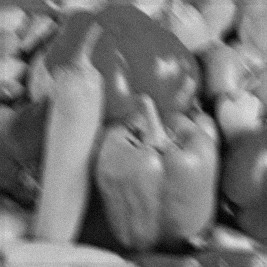}
        \caption*{\centering Blurry\\
        22.55 / 0.6009}
    \end{subfigure}
    \begin{subfigure}[t]{0.22\textwidth}
        \centering
        \spypeppersblur{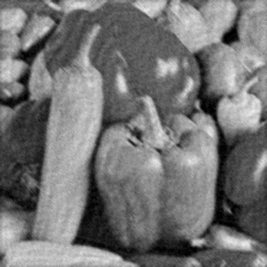}
        \caption*{\centering Tikhonov\\
        26.22 / 0.7190}
    \end{subfigure}
    \begin{subfigure}[t]{0.22\textwidth}
        \centering
        \spypeppersblur{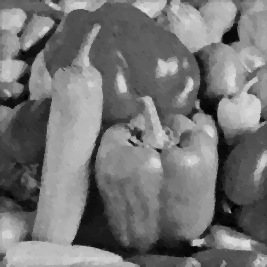}
        \caption*{\centering $\mathrm{TV}$\\
        27.42 / 0.8246}
    \end{subfigure}

    \medskip
    
    \begin{subfigure}[t]{0.22\textwidth}
        \centering
        \spypeppersblur{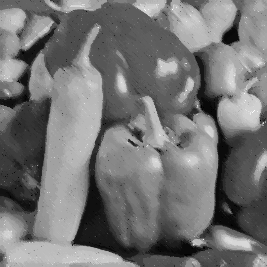}
        \caption*{\centering Gil. $\mathrm{NLTV}$\\
        26.06 / 0.7564}
    \end{subfigure}
    \begin{subfigure}[t]{0.22\textwidth}
        \centering
        \spypeppersblur{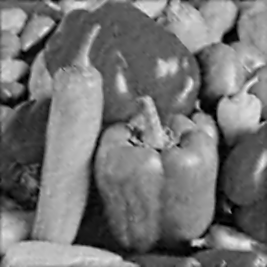}
        \caption*{\centering $\mathrm{TGV}^2$\\
        27.74 / 0.8292}
    \end{subfigure}
    \begin{subfigure}[t]{0.22\textwidth}
        \centering
        \spypeppersblur{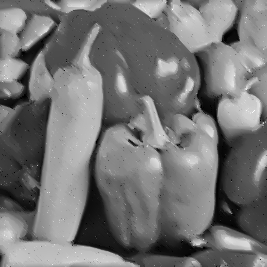}
        \caption*{\centering $\mathrm{NLTGV}^2$\\
        26.67 / 0.7805}
    \end{subfigure}
    \begin{subfigure}[t]{0.22\textwidth}
        \centering
        \spypeppersblur{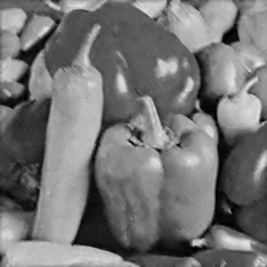}
        \caption*{\centering $\mathrm{NLTV}_\omega$\\
        \textbf{28.03} / \textbf{0.8408}}
    \end{subfigure}

    \caption{Deblurring results for the Peppers test image corrupted by motion blurring induced by the kernel \eqref{DiscreteMotKerDef} and Gaussian noise with standard deviation $\sigma_\mathrm{noise}=5$. The reconstructions obtained using Tikhonov regularization and Gilboa and Osher's $\mathrm{NLTV}$ and $\mathrm{NLTGV}^2$ retain noticeable residual noise, which manifests as salt-and-pepper-like artifacts in the latter two methods. $\mathrm{TV}$ presents pronounced staircasing artifacts. While $\mathrm{TGV}^2$ offers better noise suppression, it is unable to effectively reconstruct edges. On the other hand, $\mathrm{NLTV}_\omega$ offers the best balance of noise removal and edge and texture recovery, as seen in the zoomed-in crop.}
    \label{fig:QualCompPeppersBlurMot}
\end{figure*}

Motion deblurring results on the Peppers test image are collected in Figure~\ref{fig:QualCompPeppersBlurMot}. The Tikhonov reconstruction once again retains a substantial amount of residual noise, while $\mathrm{TV}$ exhibits pronounced staircasing effects. Gilboa and Osher's $\mathrm{NLTV}$ and $\mathrm{NLTGV}^2$, in turn, produce undesirable salt-and-pepper artifacts. Although $\mathrm{TGV}^2$ preserves texture effectively, it fails to recover edges accurately. Overall, the proposed $\mathrm{NLTV}_\omega$ regularizer produces the sharpest edges, suppresses residual noise effectively, and achieves the best recovery of textures, shadows, and reflections. These differences are particularly apparent in the enlarged crop and are also reflected in the quantitative results.

\subsection{Image Super-Resolution}

We now consider the task of single image super-resolution. In particular, we downsample the test images in Figure \ref{fig:TestImages} by a sampling factor $s$ and try to reconstruct the original images at their original, higher resolution. To facilitate discussion, we conduct the analysis of the $A$ operator in the discrete setting. Modeling the downsampling operator in the continuous setting would entail considering a lattice of Dirac deltas and changing the codomain appropriately, as done in \cite{ballester2006variational}. To downsample an image, we first convolve it with a shifted Gaussian kernel of the form \eqref{DiscreteGaussKerDef} and then perform an $s$-decimation, which keeps only one pixel out of every $s \times s$ block of the high-resolution image. The kernel is shifted by $\frac{s-1}{2}$ pixels along each axis. As a result, the value kept from each block is the Gaussian average of the image centered at the center of that block, rather than at its top-left pixel. For odd $s$, the center of each block is a pixel, and this amounts to keeping the central pixel of each block of the blurred image. For even $s$, the center falls between pixels, so the Gaussian is evaluated at half-integer distances. Centering the samples this way makes the low-resolution grid consistent with the one assumed by bicubic interpolation.

More precisely, the downsampling operator is
\begin{equation} \label{DownSamplingOp}
    A = D \circ B,    
\end{equation}
where $B : \R^{\mathcal{I}} \to \R^{\mathcal{I}}$ is the shifted Gaussian blur operator and $D : \R^{\mathcal{I}} \to \R^{\mathcal{I}_s}$ is the $s$-decimation operator, with $\mathcal{I}_s = \{1,\dots,N_1/s\} \times \{1,\dots,N_2/s\}$ being the low-resolution pixel domain. Setting
\begin{equation}
    \mu_- := \left\lfloor \tfrac{s-1}{2} - \kappa \right\rfloor, \qquad
    \mu_+ := \left\lceil \tfrac{s-1}{2} + \kappa \right\rceil,
\end{equation}
the blurred image $B u$ at pixel $i = (i_1,i_2) \in \mathcal{I}$ is given by
\begin{equation}
    (B u)_{i_1,i_2} := \sum_{t_1,t_2 = \mu_-}^{\mu_+} K_{t_1,t_2}\, u_{i_1 + t_1,\, i_2 + t_2},
\end{equation}
where $\kappa \in \Z_{>0}$ is a truncation constant and $K$ is a nonnegative kernel supported on $\{\mu_-,\dots,\mu_+\} \times \{\mu_-,\dots,\mu_+\}$ (of size $(2\kappa+1)\times(2\kappa+1)$ for odd $s$ and $(2\kappa+2)\times(2\kappa+2)$ for even $s$). The kernel is symmetric with respect to $(s-1)/2$ in each component, and satisfies
\begin{equation}
    \sum_{t_1,t_2 = \mu_-}^{\mu_+} K_{t_1,t_2} = 1.
\end{equation}
Out-of-bound indices of $u$ are treated by reflection. The decimated image $D u$ at pixel $j = (j_1,j_2) \in \mathcal{I}_s$ is
\begin{equation}
    (D u)_{j_1,j_2} = u_{s j_1,\, s j_2}.
\end{equation}

\begin{remark}
    Let us see that the $A = D \circ B$ operator satisfies \eqref{H1} and \eqref{H2'} in the discrete setting. The linearity of $A$ follows from the linearity of $B$ and $D$ and the continuity by the finiteness of the considered spaces. To verify \eqref{H2'}, let $a_1, a_2, b \in \R$ and define $g \in \R^{\mathcal{I}}$ by
    \begin{equation*}
        g_{i_1,i_2} = a_1 i_1 + a_2 i_2 + b, \qquad (i_1,i_2) \in \mathcal{I}.
    \end{equation*}
    Let $\mathcal{J}_0 \subset \mathcal{I}_s$ be the set of low-resolution pixels $j$ for which the kernel window around $sj$ lies inside $\mathcal{I}$, so that no reflection is involved:
    \begin{equation*}
        \mathcal{J}_0 = \left\{ \left\lceil \tfrac{1-\mu_-}{s} \right\rceil, \dots, \left\lfloor \tfrac{N_1-\mu_+}{s} \right\rfloor \right\} \times \left\{ \left\lceil \tfrac{1-\mu_-}{s} \right\rceil, \dots, \left\lfloor \tfrac{N_2-\mu_+}{s} \right\rfloor \right\}.
    \end{equation*}
    We additionally assume that $\mathcal{J}_0$ contains at least two indices along each axis (this holds in all our experiments), so that $\mathcal{J}_0$ contains a $2\times 2$ block and, in particular, three non-collinear points. Since $K$ sums to one and is symmetric with respect to $(s-1)/2$, we have
    \begin{equation}
         \sum_{t_1,t_2 = \mu_-}^{\mu_+} K_{t_1,t_2}\, t_1 =  \sum_{t_1,t_2 = \mu_-}^{\mu_+} K_{t_1,t_2}\, t_2 = \frac{s-1}{2}.    
    \end{equation}
    Hence, by adapting the continuous computation \eqref{ContBlurrComupation}, we obtain that for $j \in \mathcal{J}_0$
    \begin{equation}
        (B g)_{sj} =  \sum_{t_1,t_2 = \mu_-}^{\mu_+} K_{t_1,t_2}\, g_{s j_1 + t_1,\, s j_2 + t_2} = g_{s j_1, s j_2} + \tfrac{s-1}{2}(a_1 + a_2).
    \end{equation}
    Note that, because of the shift, $Bg \neq g$ in general; instead $B$ translates $g$ by $\frac{s-1}{2}$ along each axis, which does not affect the argument. Therefore, for every $j = (j_1,j_2) \in \mathcal{J}_0$,
    \begin{equation*}
        (A g)_{j_1,j_2} = (B g)_{s j_1, s j_2} = s\,(a_1 j_1 + a_2 j_2) + b + \tfrac{s-1}{2}(a_1 + a_2).
    \end{equation*}
    If $A g = 0$, this affine function of $j$ vanishes on $\mathcal{J}_0$, and in particular at three non-collinear points. Thus all its coefficients vanish, i.e., $s a_1 = s a_2 = 0$ and $b + \frac{s-1}{2}(a_1+a_2) = 0$, which gives $a_1 = a_2 = b = 0$. Hence $g = 0$, and $A$ satisfies \eqref{H2'} in the discrete setting.
\end{remark}

\begin{table}[t]
\centering
\begin{tabular}{lcccc}
\hline
& \multicolumn{2}{c}{$s=2$}
& \multicolumn{2}{c}{$s=4$} \\
\cline{2-5}
\textbf{Model}
& \textbf{PSNR $\uparrow$} & \textbf{SSIM $\uparrow$}
& \textbf{PSNR $\uparrow$} & \textbf{SSIM $\uparrow$} \\
\hline
Bicubic
& 25.69 & 0.6951
& 24.58 & 0.6688 \\

Tikhonov
& 26.50 & 0.7150
& 24.82 & 0.6660 \\

$\mathrm{TV}$
& 27.07 & 0.7662
& 25.03 & 0.7022 \\

Gil. $\mathrm{NLTV}$
& 26.85 & 0.7482
& 25.03 & 0.6910 \\

$\mathrm{TGV}^2$
& 27.27 & \underline{0.7732}
& 25.21 & 0.7108 \\

$\mathrm{NLTGV}^2$
& \underline{27.30} & 0.7651
& \underline{25.39} & \underline{0.7137} \\

$\mathrm{NLTV}_\omega$
& \textbf{27.65} & \textbf{0.7813}
& \textbf{25.62} & \textbf{0.7254} \\
\hline
\end{tabular}
\caption{Quantitative comparison of super-resolution methods on the testing set from Figure \ref{fig:TestImages} for different sampling rates. Best results are highlighted in \textbf{bold} and second-best results are \underline{underlined}. $\mathrm{NLTV}_\omega$ outperforms the other regularizers at both sampling rates and for both metrics.}
\label{tab:QuantCompSR}
\end{table}

\begin{figure*}[t]
    \centering
    \begin{subfigure}[t]{0.22\textwidth}
        \centering
        \spyboat{TestImages/boat.png}
        \caption*{\centering Ground Truth\\
        PSNR / SSIM}
    \end{subfigure}
    \begin{subfigure}[t]{0.22\textwidth}
        \centering
        \spyboat{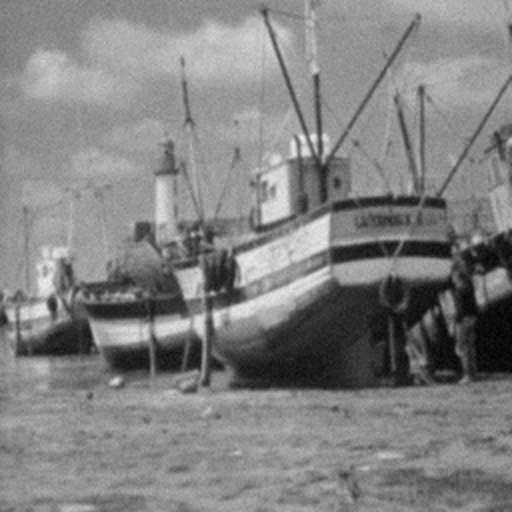}
        \caption*{\centering Bicubic\\
        26.29 / 0.6792}
    \end{subfigure}
    \begin{subfigure}[t]{0.22\textwidth}
        \centering
        \spyboat{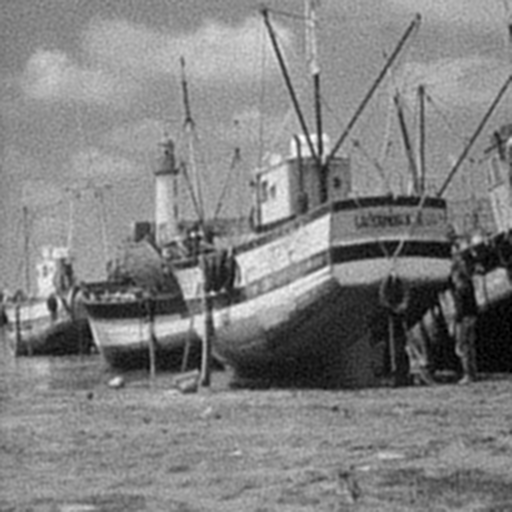}
        \caption*{\centering Tikhonov\\
        27.03 / 0.6974}
    \end{subfigure}
    \begin{subfigure}[t]{0.22\textwidth}
        \centering
        \spyboat{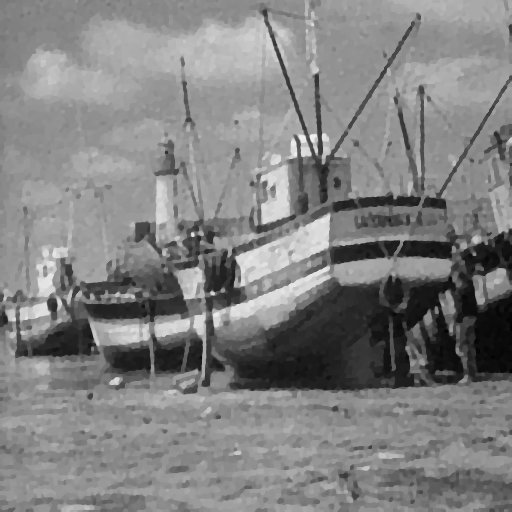}
        \caption*{\centering TV\\
        27.37 / 0.7310}
    \end{subfigure}

    \medskip
    
    \begin{subfigure}[t]{0.22\textwidth}
        \centering
        \spyboat{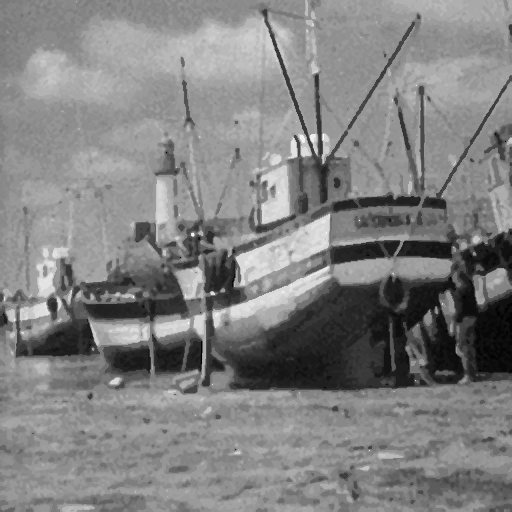}
        \caption*{\centering Gil. $\mathrm{NLTV}$\\
        27.24 / 0.7156}
    \end{subfigure}
    \begin{subfigure}[t]{0.22\textwidth}
        \centering
        \spyboat{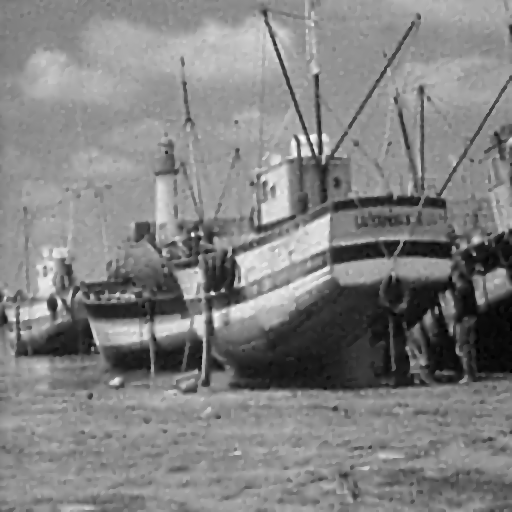}
        \caption*{\centering $\mathrm{TGV}^2$\\
        27.58 / 0.7399}
    \end{subfigure}
    \begin{subfigure}[t]{0.22\textwidth}
        \centering
        \spyboat{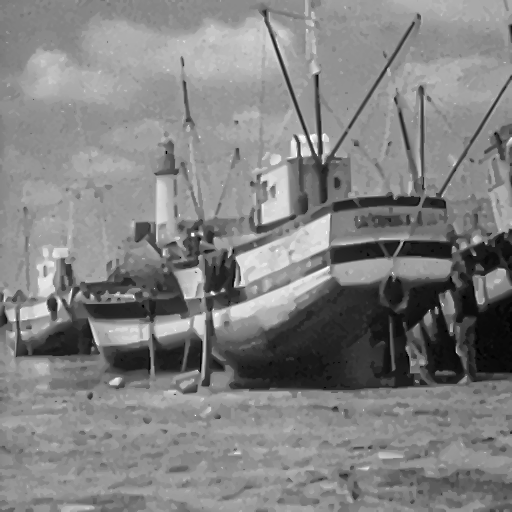}
        \caption*{\centering $\mathrm{NLTGV}^2$\\
        27.51 / 0.7306}
    \end{subfigure}
    \begin{subfigure}[t]{0.22\textwidth}
        \centering
        \spyboat{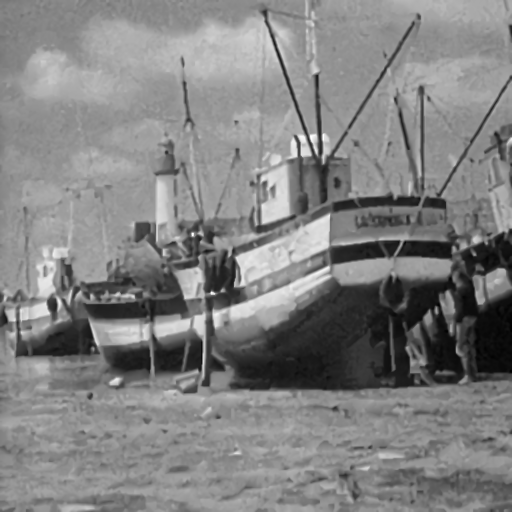}
        \caption*{\centering $\mathrm{NLTV}_\omega$\\
        \textbf{27.96} / \textbf{0.7813}}
    \end{subfigure}

    \caption{Super-resolution results for the Boat test image corrupted by the downsampling operator \eqref{DownSamplingOp} with sampling factor $s=2$ and Gaussian noise with standard deviation $\sigma_\mathrm{noise}=5$. As in the previous experiments, Tikhonov regularization presents a substantial amount of residual noise, and $\mathrm{TV}$ and Gilboa and Osher's $\mathrm{NLTV}$ exhibit pronounced staircasing patterns. $\mathrm{TGV}^2$ fails to effectively reconstruct edges, while $\mathrm{NLTGV}^2$ oversmooths certain regions. By contrast, $\mathrm{NLTV}_\omega$ achieves a better balance between noise suppression and edge and texture preservation. In particular, it yields the most clearly legible text on the boat's hull, as can be seen in the enlarged crop.}
    \label{fig:QualCompBoatSR2}
\end{figure*}

We consider sampling factors $s\in\{2,4\}$ and fix the standard deviation of the Gaussian kernel as $1.5$ and its size as $12$, which corresponds to $\kappa = 5$. As the sizes of the cactus and peppers test images are not divisible by the considered sampling factors, we extract a $256 \times 256$ crop out of each image to obtain suitable references. We finally add zero-mean Gaussian noise with standard deviation $\sigma_{\mathrm{noise}} = 5$ to the downsampled data.

The nonlocal regularizers require NLMeans weights defined on the high-resolution pixel grid, so these weights cannot be computed directly from the low-resolution observations. Instead, we compute them from an upsampled version obtained by bicubic interpolation.

Table \ref{tab:QuantCompSR} reports the average metrics over the testing set at sampling rates $s=2$ and $4$. $\mathrm{NLTV}_\omega$ again attains the best PSRN and SSIM in all cases. The runner-up is $\mathrm{NLTGV}^2$, in PSNR for both factors and in SSIM for $s=4$, while $\mathrm{TGV}^2$ takes the second places in terms of SSIM for $s=2$.

Figure~\ref{fig:QualCompBoatSR2} shows the Boat test image together with the super-resolution results obtained for $s=2$ by bicubic upsampling and each of the considered regularizers. As in the previous experiments, Tikhonov regularization fails to remove the noise effectively, while $\mathrm{TV}$ and Gilboa and Osher's $\mathrm{NLTV}$ produce unnatural piecewise constant reconstructions. These limitations are particularly evident in the enlarged crop. Both $\mathrm{TGV}^2$ and $\mathrm{NLTGV}^2$ perform better in terms of noise removal and structure preservation. However, $\mathrm{TGV}^2$ struggles to recover edges accurately, whereas $\mathrm{NLTGV}^2$ tends to oversmooth certain regions, giving rise to block-like artifacts. By contrast, the proposed $\mathrm{NLTV}_\omega$ regularizer suppresses noise as well as the two latter methods, while providing the best recovery of image structures and pixel intensities. In particular, its reconstruction preserves the text on the boat's hull more clearly than any of the other models.

\subsection{Image Inpainting}

As a last experiment, we consider image inpainting, which consists in reconstructing the content of a region where the image data has been lost, using only the information available around it. In the continuous setting, this task corresponds to the operator $A:L^p(\Omega) \to L^p(\Omega)$ given by
\begin{equation} \label{AInpaitingDef}
    Au = \mathcal{X}_D u,
\end{equation}
where $D \subset \Omega$ is the region of the image in which the pixel values are known.

\begin{remark}
    Provided that $D \subset \Omega$ has positive measure, it is immediate to check that operator $A$ given by \eqref{AInpaitingDef} satisfies \eqref{H1} and \eqref{H2'}.
\end{remark}

In this experiment, we corrupt the house test image with random scribbles, removing approximately $11\%$ of its pixels. For the nonlocal regularizers, we compute the nonlocal weights from an initial inpainting obtained via harmonic filling, which serves as a baseline. 

\begin{figure*}[t]
    \centering
    \begin{subfigure}[t]{0.25\textwidth}
        \centering
        \includegraphics[width=\linewidth]{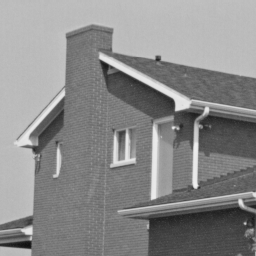}
        \caption*{\centering Ground Truth\\
        PSNR / SSIM}
    \end{subfigure}
    \begin{subfigure}[t]{0.25\textwidth}
        \centering
        \includegraphics[width=\linewidth]{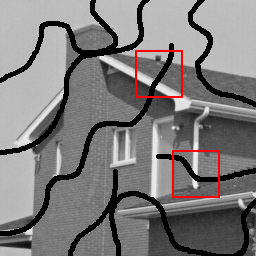}
        \caption*{\centering Masked Image\\
        - / -}
    \end{subfigure}
    \begin{subfigure}[t]{0.25\textwidth}
        \centering
        \includegraphics[width=\linewidth]{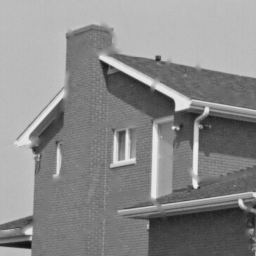}
        \caption*{\centering Harmonic Fill\\
        35.14 / 0.9705}
    \end{subfigure}

    \medskip

    \begin{subfigure}[t]{0.25\textwidth}
        \centering
        \includegraphics[width=\linewidth]{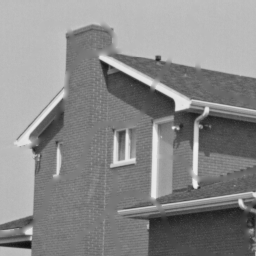}
        \caption*{\centering Tikhonov\\
        35.40 / 0.9720}
    \end{subfigure}
    \begin{subfigure}[t]{0.25\textwidth}
        \centering
        \includegraphics[width=\linewidth]{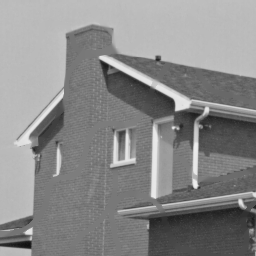}
        \caption*{\centering $\mathrm{TV}$\\
        34.57 / 0.9712}
    \end{subfigure}
    \begin{subfigure}[t]{0.25\textwidth}
        \centering
        \includegraphics[width=\linewidth]{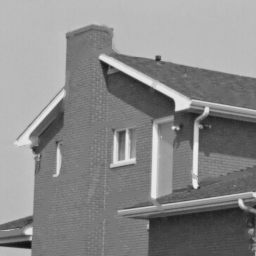}
        \caption*{\centering Gil. $\mathrm{NLTV}$\\
        35.91 / 0.9732}
    \end{subfigure}

    \medskip

    \begin{subfigure}[t]{0.25\textwidth}
        \centering
        \includegraphics[width=\linewidth]{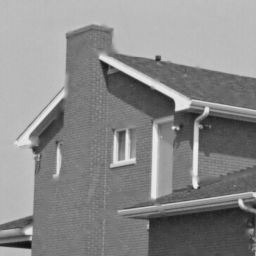}
        \caption*{\centering $\mathrm{TGV}^2$\\
        \textbf{36.44} / \underline{0.9748}}
    \end{subfigure}
    \begin{subfigure}[t]{0.25\textwidth}
        \centering
        \includegraphics[width=\linewidth]{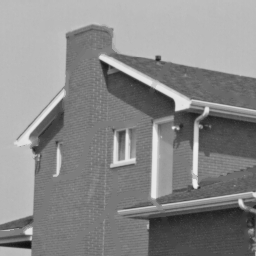}
        \caption*{\centering $\mathrm{NLTGV}^2$\\
        35.89 / 0.9732}
    \end{subfigure}
    \begin{subfigure}[t]{0.25\textwidth}
        \centering
        \includegraphics[width=\linewidth]{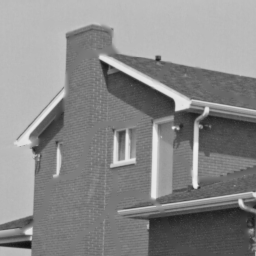}
        \caption*{\centering $\mathrm{NLTV}_\omega$\\
        \underline{36.38} / \textbf{0.9760}}
    \end{subfigure}

    \caption{Inpainting results for the house test image corrupted by random scribbles, with approximately $11\%$ of pixels missing. The red boxes indicate the regions shown in detail in Figure~\ref{fig:QualCompHouseInpaintingCrops}. Quantitative performance metrics (PSNR and SSIM) are reported below each reconstruction.}
    \label{fig:QualCompHouseInpainting}
\end{figure*}

\begin{figure*}[t]
    \centering
    \begin{subfigure}[t]{0.15\textwidth}
        \centering
        \includegraphics[width=\linewidth]{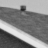}
        \caption*{\centering Ground Truth}
    \end{subfigure}
    \begin{subfigure}[t]{0.15\textwidth}
        \centering
        \includegraphics[width=\linewidth]{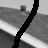}
        \caption*{\centering Masked Image}
    \end{subfigure}
    \begin{subfigure}[t]{0.15\textwidth}
        \centering
        \includegraphics[width=\linewidth]{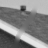}
        \caption*{\centering Harmonic Fill}
    \end{subfigure}\hfill
    \begin{subfigure}[t]{0.15\textwidth}
        \centering
        \includegraphics[width=\linewidth]{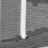}
        \caption*{\centering Ground Truth}
    \end{subfigure}
    \begin{subfigure}[t]{0.15\textwidth}
        \centering
        \includegraphics[width=\linewidth]{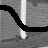}
        \caption*{\centering Masked Image}
    \end{subfigure}
    \begin{subfigure}[t]{0.15\textwidth}
        \centering
        \includegraphics[width=\linewidth]{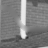}
        \caption*{\centering Harmonic Fill}
    \end{subfigure}

    \medskip

    \begin{subfigure}[t]{0.15\textwidth}
        \centering
        \includegraphics[width=\linewidth]{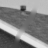}
        \caption*{\centering Tikhonov}
    \end{subfigure}
    \begin{subfigure}[t]{0.15\textwidth}
        \centering
        \includegraphics[width=\linewidth]{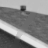}
        \caption*{\centering $\mathrm{TV}$}
    \end{subfigure}
    \begin{subfigure}[t]{0.15\textwidth}
        \centering
        \includegraphics[width=\linewidth]{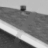}
        \caption*{\centering Gil. $\mathrm{NLTV}$}
    \end{subfigure}\hfill
    \begin{subfigure}[t]{0.15\textwidth}
        \centering
        \includegraphics[width=\linewidth]{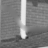}
        \caption*{\centering Tikhonov}
    \end{subfigure}
    \begin{subfigure}[t]{0.15\textwidth}
        \centering
        \includegraphics[width=\linewidth]{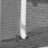}
        \caption*{\centering $\mathrm{TV}$}
    \end{subfigure}
    \begin{subfigure}[t]{0.15\textwidth}
        \centering
        \includegraphics[width=\linewidth]{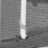}
        \caption*{\centering Gil. $\mathrm{NLTV}$}
    \end{subfigure}
    \medskip

    \begin{subfigure}[t]{0.15\textwidth}
        \centering
        \includegraphics[width=\linewidth]{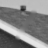}
        \caption*{\centering $\mathrm{TGV}^2$}
    \end{subfigure}
    \begin{subfigure}[t]{0.15\textwidth}
        \centering
        \includegraphics[width=\linewidth]{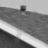}
        \caption*{\centering $\mathrm{NLTGV}^2$}
    \end{subfigure}
    \begin{subfigure}[t]{0.15\textwidth}
        \centering
        \includegraphics[width=\linewidth]{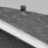}
        \caption*{\centering $\mathrm{NLTV}_\omega$}
    \end{subfigure}\hfill
    \begin{subfigure}[t]{0.15\textwidth}
        \centering
        \includegraphics[width=\linewidth]{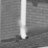}
        \caption*{\centering $\mathrm{TGV}^2$}
    \end{subfigure}
    \begin{subfigure}[t]{0.15\textwidth}
        \centering
        \includegraphics[width=\linewidth]{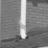}
        \caption*{\centering $\mathrm{NLTGV}^2$}
    \end{subfigure}
    \begin{subfigure}[t]{0.15\textwidth}
        \centering
        \includegraphics[width=\linewidth]{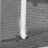}
        \caption*{\centering $\mathrm{NLTV}_\omega$}
    \end{subfigure}

    \caption{Close-up comparison of the inpainting results in two selected regions of the house test image. The selected regions are highlighted in red in Figure~\ref{fig:QualCompHouseInpainting}. The proposed $\mathrm{NLTV}_\omega$ regularizer better preserves edge sharpness and distinguishes adjacent regions, yielding more accurate structural and intensity reconstruction.}
    \label{fig:QualCompHouseInpaintingCrops}
\end{figure*}

Figure~\ref{fig:QualCompHouseInpainting} shows the ground truth image and the mask used to corrupt it, together with the inpainting results of the harmonic filling and all considered regularizers, while Figure~\ref{fig:QualCompHouseInpaintingCrops} zooms in on the two regions highlighted in red there. Quantitative metrics are reported below each reconstruction.
 
As Tikhonov regularization is essentially an unconstrained counterpart of harmonic filling, the two methods give very comparable results, both visually and quantitatively. $\mathrm{TV}$ regularization tends to fill the missing regions with nearly constant patterns, resulting in poor overall inpainting performance, as reflected in its low PSNR. Among the local methods, $\mathrm{TGV}^2$ regularization performs best, achieving the highest PSNR and the second-highest SSIM among all the methods considered.

As for the nonlocal regularizers, Gilboa and Osher's $\mathrm{NLTV}$ and $\mathrm{NLTGV}^2$ yield very similar results, both visually and in terms of metrics. These regularizers recover edges well, although they tend to fill the missing regions with nearly constant patterns, to a lesser extent than $\mathrm{TV}$. Furthermore, they struggle to correct inaccuracies inherited from the harmonic filling used to compute the nonlocal weights.

Finally, our proposed $\mathrm{NLTV}_\omega$ regularizer achieves the second-highest PSNR and the highest SSIM among all the methods considered. Visually, it yields the sharpest edges within the missing regions and the most accurate intensities, as Figure~\ref{fig:QualCompHouseInpaintingCrops}
highlights. Moreover, unlike Gilboa and Osher's $\mathrm{NLTV}$ and $\mathrm{NLTGV}^2$, our regularizer is not entirely constrained by the nonlocal weights, since it can deviate from the harmonic filling initialization and thereby recover more natural-looking structures.

\section{Conclusions and Future Work} \label{ConclusionsSec}

In this work we have introduced the nonlocal total variation $\mathrm{NLTV}_\omega$ regularizer for image restoration, a weighted second-order nonlocal functional that builds upon the work of Kindermann et~al. \cite{Kindermann}. Moreover, we have expanded our previous theoretical study \cite{alcover2026nonlocal} by relaxing the smoothness condition on the weight function, removing the symmetry requirement, and adopting a weaker definition of $\mathrm{NLTV}_\omega$. We have shown that $\mathrm{NLTV}_\omega$ and its corresponding space $\mathrm{NLBV}_\omega$ retain some desirable properties in this new relaxed setting. Moreover, we have established condition \eqref{H0} on the weight functions, which is satisfied by a broad class of weights relevant to image processing and allows for the embedding of $\mathrm{NLBV}_\omega$ into the classical $\mathrm{BV}$ space to hold. Building on these properties and the resulting embedding, we have considered a nonlocal variant of the variational image restoration model of Chambolle and Lions \cite{chambolle1997image} and proved that, under suitable assumptions, the associated minimization problem admits a unique solution.

On the computational side, we have solved the discretized problem with a first-order primal-dual algorithm \cite{chambolle2011first}, and have implemented analogous schemes for several local and nonlocal regularizers in order to conduct an extensive experimental evaluation across a range of image restoration tasks. The results have shown that $\mathrm{NLTV}_\omega$ achieves the best performance in nearly all the considered experiments, both visually and quantitatively. In particular, we have demonstrated its robustness to nonlocal weights derived in non-ideal scenarios, and its ability to preserve image textures and correctly recover pixel intensity.

Overall, these theoretical and experimental results establish $\mathrm{NLTV}_\omega$ as a promising regularization approach for image restoration, combining a rigorous variational framework with strong practical performance.

We conclude by outlining several directions for future work. First, unlike $\mathrm{TGV}^2$ \cite{bredies2010total} and $\mathrm{NLTGV}^2$ \cite{ranftl2014non}, the proposed $\mathrm{NLTV}_\omega$ regularizer does not explicitly account for first-order variations. We are currently working on models that combine $\mathrm{NLTV}_\omega$ with a first-order regularization term, potentially alleviating some of its limitations, particularly the appearance of gradient-like artifacts in smooth regions. A second direction worth exploring is the integration of  $\mathrm{NLTV}_\omega$ into a deep unfolding architecture. This would consist on considering an optimization algorithm involving $\mathrm{NLTV}_\omega$, such as \eqref{NLPDIdentityAlgorithm} or \eqref{NLPDqAlgorithm}, and replacing some of its constituent operators with learnable modules. Such an approach could combine the strong empirical performance of data-driven methods with the mathematical guarantees established in the present work. Finally, more expressive strategies for computing the nonlocal weights could be investigated. These could include hand-crafted weights specifically designed to measure gradient differences, as well as machine-learning-based approaches to refine the weight computation and identify the image features most relevant to a given restoration task.

\bibliographystyle{abbrv}
\bibliography{Biblio}

\begin{appendix}
\section{Technical Proofs} \label{app:TechProofs}

In this appendix we give the detailed proofs of Lemma \ref{OldNLTVDefEquivalence} and Proposition \ref{C1NLTV}.

\subsection{Proof of Lemma \ref{OldNLTVDefEquivalence}}

We restate the lemma for convenience.

\begin{lemma} \label{OldNLTVDefEquivalenceapp}
    Let $\omega$ be an admissible weight function which, in addition, is symmetric and belongs to $\mathcal{C}^1(\Omega\times\Omega)$. Then

    $$\mathrm{NLTV}_\omega(u)= \sup \left\{  \int_\Omega u(x) \; \mathrm{div}_\omega (\phi)(x) \mathrm{d}x : \phi \in \mathcal{C}^1_c (\Omega \times \Omega; \mathbb{R}^N), \| \phi \|_{L^\infty(\Omega\times\Omega;\R^N)} \le 1 \right\}$$
    for any $u \in L^1(\Omega)$, where $\| \cdot \|_{L^\infty(\Omega\times\Omega;\R^N)}$ denotes the uniform norm in $\mathcal{C}_c(\Omega\times\Omega;\R^N)$ and

    \begin{equation}
        \mathrm{div}_{\omega}(\phi)(x) := \int_\Omega \mathrm{div}_x \left[ \omega(x,y) ( \phi(y,x) - \phi (x,y) ) \right] \mathrm{d}y
    \end{equation}
    for $x \in \Omega$.
\end{lemma}

\begin{proof}
    Fix $u \in L^1(\Omega)$. Given $\phi \in \mathcal{C}_c^1(\Omega\times\Omega;\mathbb{R}^N)$ such that $\| \phi \|_{L^\infty(\Omega\times\Omega;\R^N)} \le 1$, define $\phi^\omega : \Omega\times\Omega \to \mathbb{R}^N$ as $\phi^\omega := \omega \phi$. By definition, $\phi^\omega \in \mathcal{C}_c^1(\Omega\times\Omega;\mathbb{R}^N)$, $|\phi^\omega|\le \omega$ in $\Omega\times\Omega$ (note that since $\omega$ is now continuous, the inequality can be imposed in the whole $\Omega\times\Omega$ instead of in the a.e. sense) and, using the symmetry of $\omega$, $\mathrm{div}_\omega(\phi) = \mathrm{div}_{\mathrm{NL}} (\phi^\omega)$ in $\Omega$. Therefore,

    $$\int_\Omega u(x) \;\mathrm{div}_\omega (\phi)(x) \mathrm{d}x = \int_\Omega u(x)\; \mathrm{div}_{\mathrm{NL}} (\phi^\omega)(x) \mathrm{d}x \le \mathrm{NLTV}_\omega(u),$$
    from which
    \begin{equation}
        \sup \left\{  \int_\Omega u(x) \; \mathrm{div}_\omega (\phi)(x) \mathrm{d}x : \phi \in \mathcal{C}^1_c (\Omega \times \Omega; \mathbb{R}^N), \| \phi \|_{L^\infty(\Omega\times\Omega;\R^N)} \le 1 \right\} \le \mathrm{NLTV}_\omega(u).
    \end{equation}

    On the other hand, given $\phi \in \mathcal{C}_c^1(\Omega\times\Omega;\mathbb{R}^N)$ such that $|\phi|\le \omega$ in $\Omega\times\Omega$ for each $\epsilon>0$ define $\phi_{\omega,\epsilon}:\Omega\times\Omega\to \mathbb{R}^N$ as
    \begin{equation}
        \phi_{\omega,\epsilon} := \frac{\phi}{\omega + \epsilon}.
    \end{equation}
    Since $\phi_{\omega,\epsilon}$ is the quotient of two $\mathcal{C}^1$ functions, $\omega+\epsilon>0$ and $\mathrm{supp}(\phi_{\omega,\epsilon}) = \mathrm{supp}(\phi)$, we have $\phi_{\omega,\epsilon}\in C_c^1(\Omega\times\Omega;\R^N)$ for each $\epsilon>0$. Moreover,

    \begin{equation}
        |\phi_{\omega,\epsilon}| = \frac{|\phi|}{\omega + \epsilon} \le \frac{\omega}{\omega + \epsilon} \le 1
    \end{equation}
    in $\Omega\times\Omega$, from which $\|\phi_{\omega,\epsilon}\|_{L^\infty(\Omega\times\Omega;\R^N)} \le 1$. Hence, we have that
    \begin{align} \label{BaldiAux0}
        \int_\Omega u(x) \;\mathrm{div}_{\omega} (\phi_{\omega,\epsilon})(x) \mathrm{d}x\le \sup \left\{  \int_\Omega u(x) \; \mathrm{div}_\omega (\phi)(x) \mathrm{d}x : \phi \in \mathcal{C}^1_c (\Omega \times \Omega; \mathbb{R}^N), \| \phi \|_{L^\infty(\Omega\times\Omega;\R^N)} \le 1 \right\}
    \end{align}
    for each $\epsilon>0$. Let us see that the left-hand side converges precisely to $\int_\Omega u \;\mathrm{div}_{\mathrm{NL}}(\phi) \mathrm{d}x$ as $\epsilon \to 0$. To do so, we note that
    \begin{equation}
        \omega \phi_{\omega,\epsilon} = \phi - \frac{\epsilon \phi}{\omega + \epsilon} =: \phi - g_\epsilon,
    \end{equation}
    and hence, using again the symmetry of $\omega$,
    \begin{equation} \label{divergencesequality}
        \mathrm{div}_\omega(\phi_{\omega,\epsilon}) = \mathrm{div}_{\mathrm{NL}}(\omega \phi_{\omega,\epsilon}) = \mathrm{div}_{\mathrm{NL}}(\phi) - \mathrm{div}_{\mathrm{NL}}(g_\epsilon).
    \end{equation}

    A direct computation gives
    \begin{align} \label{divg_epsilon}
        \mathrm{div}_x (g_\epsilon(x,y)) &= \frac{-\epsilon \nabla_x \omega(x,y)}{(\omega(x,y) + \epsilon)^2} \cdot \phi(x,y) + \frac{\epsilon}{\omega(x,y)+\epsilon} \mathrm{div}_x(\phi(x,y))
    \end{align}
    for $(x,y) \in \Omega\times\Omega$. In particular, $\mathrm{div}_x (g_\epsilon) = 0$ outside $\mathrm{supp}(\phi)$. For $(x,y) \in \mathrm{supp}(\phi)$, since $\epsilon/(t+\epsilon) \le 1$ and $t\epsilon/(t+\epsilon)^2 \le 1/4$ for $t\ge 0$, we get
    \begin{align}
        |\mathrm{div}_x (g_\epsilon(x,y))| &\le |\nabla_x \omega(x,y)| \frac{\epsilon|\phi(x,y)|}{(\omega(x,y) + \epsilon)^2} + |\mathrm{div}_x(\phi(x,y))|\\
        &\le |\nabla_x \omega(x,y)| \frac{\epsilon\omega(x,y)}{(\omega(x,y) + \epsilon)^2} + |\mathrm{div}_x(\phi(x,y))|\\
        &\le \frac{1}{4} \|\nabla \omega\|_{L^\infty(\mathrm{supp}(\phi);\R^N\times\R^N)} + C \|D \phi\|_{L^\infty(\Omega\times\Omega;\R^N\times\R^{2N})},
    \end{align}
    where $\nabla \omega$ denotes the ``full" gradient of $\omega$, $D \phi$ denotes the full derivative of $\phi$, and $C>0$ is a constant depending on $N$ and the choice of norm. This is a finite uniform bound because $\omega$ and $\phi$ are $\mathcal{C}^1$, and $\phi$ has compact support. From this, we deduce that $\mathrm{div}_x (g_\epsilon)$ is uniformly bounded in $\Omega\times\Omega$. Analogously, defining $\tilde \phi(x,y) := \phi(y,x)$ for $(x,y)\in\Omega\times\Omega$, we see that
    \begin{equation}
         |\mathrm{div}_x (g_\epsilon(y,x))| \le \frac{1}{4} \|\nabla \omega\|_{L^\infty(\mathrm{supp}( \phi);\R^N\times\R^N)} + C \|D \phi\|_{ L^\infty(\Omega\times\Omega;\R^N\times\R^{2N})}
    \end{equation}
     for $(x,y) \in \mathrm{supp}(\tilde \phi)$, and the left-hand side is zero outside this set. Therefore, for any $(x,y) \in \Omega\times\Omega$,
    \begin{equation} \label{BaldiAux1}
         |\mathrm{div}_x (g_\epsilon(y,x) - g_\epsilon(x,y))| \le \Phi(x,y),
    \end{equation}
    where
    \begin{equation}
        \Phi(x,y) = \begin{cases}
            \frac{1}{2} \|\nabla \omega\|_{L^\infty(\mathrm{supp}(\phi);\R^N\times\R^N)} + 2 C\| D\phi\|_{L^\infty(\Omega\times\Omega;\R^N\times\R^{2N})} & \mbox{ for } (x,y) \in \mathrm{supp}(\phi) \cup \mathrm{supp}(\tilde \phi),\\
            0 & \mbox{ otherwise}.
        \end{cases}
    \end{equation}
    
    Now, for $(x,y) \in \Omega\times\Omega$ such that $\omega(x,y)>0$, from \eqref{divg_epsilon} it is deduced that
    \begin{equation}
        \mathrm{div}_x (g_\epsilon(y,x) - g_\epsilon(x,y)) \xrightarrow[\epsilon\to 0]{} 0.
    \end{equation}
    For $(x,y) \in \Omega\times\Omega$ such that $\omega(x,y) = 0$, since $\omega$ is nonnegative and differentiable everywhere in $\Omega\times\Omega$, we must also have $\nabla \omega(x,y) = 0$ (as $\omega$ has a minimum there). A first-order Taylor expansion of $\omega$ at $(x,y)$ then gives $\omega((x,y)+z) = o(|z|)$ as $z\to 0$. Moreover, $|\phi(x,y)| \le \omega(x,y)=0$, so that
    \begin{equation}
        |\phi((x,y) + z) - \phi(x,y)| = |\phi((x,y)+z)|\le \omega((x,y)+z) = o(|z|)
    \end{equation}
    as $z\to 0$. Dividing by $|z|$ and taking the limit $z \to 0$ in the inequality above shows that $D \phi(x,y) = 0$. By symmetry of $\omega$, we can repeat the same argument to obtain that $D\phi(y,x)=0$. This yields
    \begin{equation}
        \mathrm{div}_x (g_\epsilon(y,x) - g_\epsilon(x,y)) = 0 \xrightarrow[\epsilon\to 0]{} 0,
    \end{equation}
    for such points. From this, using the dominated convergence theorem, we obtain that
    \begin{equation}
        \mathrm{div}_\mathrm{NL}(g_\epsilon) = \int_\Omega  \mathrm{div}_x (g_\epsilon(y,\cdot) - g_\epsilon(\cdot,y)) \mathrm{d}y \xrightarrow[\epsilon \to 0]{} \int_\Omega 0 \; \mathrm{d}y = 0.
    \end{equation}
    point-wise in $\Omega$. Note that, by \eqref{BaldiAux1}, $u\; \mathrm{div}_\mathrm{NL}(g_\epsilon)$ is dominated a.e.~in $\Omega$ by $|u| \int_\Omega \Phi(\cdot, y) \mathrm{d}y$, which is integrable, and dominated convergence theorem once more yields that
    \begin{equation}
        \int_\Omega u \; \mathrm{div}_\mathrm{NL}(g_\epsilon) \mathrm{d}x \xrightarrow[\epsilon\to 0]{}\int_\Omega 0 \; \mathrm{d}x = 0.
    \end{equation}
    Together with \eqref{divergencesequality} this shows that
    \begin{equation}
        \int_\Omega u\;\mathrm{div}_\omega(\phi_{\omega,\epsilon}) \mathrm{d}x \xrightarrow[\epsilon\to 0]{} \int_\Omega u\;\mathrm{div}_\mathrm{NL}(\phi) \mathrm{d}x.
    \end{equation}
    
    Letting $\epsilon\to 0$ in \eqref{BaldiAux0} we obtain
    \begin{equation}
        \int_\Omega u\;\mathrm{div}_\mathrm{NL}(\phi) \mathrm{d}x \le \sup \left\{  \int_\Omega u(x) \; \mathrm{div}_\omega (\phi)(x) \mathrm{d}x : \phi \in \mathcal{C}^1_c (\Omega \times \Omega; \mathbb{R}^N), \| \phi \|_{L^\infty(\Omega\times\Omega;\R^N)} \le 1 \right\},
    \end{equation}
    from which the remaining inequality follows by taking supremum over $\phi$.
\end{proof}

\subsection{Proof of Proposition \ref{C1NLTV}}

We restate the proposition for convinience.

\begin{prop}
\label{C1NLTVapp}
    Let $\omega$ be an admissible weight function. For any $u \in L^1(\Omega)$ with weak first-order derivatives we have
    \begin{equation} \label{NLTVFormulaapp}
        \mathrm{NLTV}_\omega(u) = \int_{\Omega \times \Omega} \omega(x,y) |\nabla u (y) - \nabla u(x)| \mathrm{d}x \mathrm{d}y,
    \end{equation}
    where $\nabla u$ is the weak gradient of $u$.
\end{prop}

As the proof is very technical, we give some insight into the overall argument. Inequality
\begin{equation}
    \mathrm{NLTV}_\omega(u) \le \int_{\Omega \times \Omega} \omega(x,y) |\nabla u (y) - \nabla u(x)| \mathrm{d}x \mathrm{d}y
\end{equation}
is straightforward. We discuss the other inequality.

In our previous article \cite{alcover2026nonlocal}, $\omega$ was embedded into the nonlocal divergence, and $\mathrm{NLTV}_\omega$ searched its supremum on $\phi \in \mathcal{C}_c^1(\Omega\times\Omega;\R^N)$ functions satisfying $|\phi| \le 1$ in $\Omega\times\Omega$ (see Lemma \ref{OldNLTVDefEquivalence}). Similarly to the local case, the argument was based on approximating
\begin{equation}
    \frac{\nabla_\omega u}{|\nabla_\omega u|} \mathcal{X}_{\{\nabla_\omega u \neq 0\}}
\end{equation}
by test functions bounded by $1$ in an appropriate sense, where, simplifying for the sake of this discussion, the nonlocal gradient read as $\nabla_\omega u (x,y) = \omega(x,y) (\nabla u(y) - \nabla u(x))$ for almost every $(x,y) \in \Omega\times \Omega$.

In this new setting, $\omega$ is embedded instead on the conditions of the test functions, which are $\phi \in \mathcal{C}_c^1(\Omega\times\Omega;\R^N)$ and $|\phi| \le \omega$ a.e.~in $\Omega\times\Omega$. Moreover, the nonlocal gradient is simply $\nabla_{\mathrm{NL}} u(x,y) = \nabla u(y) - \nabla u(x)$ for almost every $(x,y) \in \Omega\times\Omega$. Hence, what we need now is to approximate
\begin{equation}
    \frac{\nabla_{\mathrm{NL}} u}{|\nabla_{\mathrm{NL}} u|} \mathcal{X}_{\{\nabla_{\mathrm{NL}} u \neq 0\}} \omega
\end{equation}
by $\mathcal{C}_c^1(\Omega\times\Omega;\R^N)$ functions which are bounded a.e.~by $\omega$. If $\omega$ is only continuous a.e., the function above may not even belong to $L^1_{\mathrm{loc}}(\Omega\times\Omega;\R^N)$, and approximating it directly while keeping the bound is very difficult.

What we actually do is approximate
\begin{equation}
    \frac{\nabla_{\mathrm{NL}} u}{|\nabla_{\mathrm{NL}} u|} \mathcal{X}_{\{\nabla_{\mathrm{NL}} u \neq 0\}}
\end{equation}
by $\mathcal{C}_c^1$ functions bounded by $1$ in $L^1$, and argue that $\omega$ can be approximated from below by nonnegative Lipschitz functions, and these in turn can be approximated uniformly by nonnegative $\mathcal{C}^1$ functions bounded by $\omega$. All these approximations are done locally, i.e.~on a relatively compact subset of $\Omega\times\Omega$, and we use monotone convergence theorem to pass from the Lipschitz functions to $\omega$, and to obtain the inequality in the whole $\Omega\times\Omega$.  

\begin{proof}[Proof of Proposition \ref{C1NLTVapp}]
    For almost every $(x,y) \in \Omega\times\Omega$ define
    \begin{equation}
        \nabla_{\mathrm{NL}} u(x,y) := \nabla u(y) - \nabla u(x),
    \end{equation}
    which belongs to $L^1_{\mathrm{loc}}(\Omega\times\Omega;\R^N)$ by the fact that $\nabla u \in L^1_{\mathrm{loc}}(\Omega;\R^N)$.

    We first prove the inequality $\mathrm{NLTV}_\omega(u) \le \int_\Omega \int_\Omega \omega |\nabla_{\mathrm{NL}} u| \mathrm{d}x \mathrm{d}y$. A straightforward integration by parts argument shows that
    \begin{equation} \label{C1NLTVAux0}
        \int_\Omega u(x) \;\mathrm{div}_{\mathrm{NL}}(\phi)(x) \mathrm{d}x = - \int_\Omega \int_\Omega \phi(x,y) \cdot \nabla_{\mathrm{NL}} u(x,y) \mathrm{d}x \mathrm{d}y
    \end{equation}
    for all $\phi \in \mathcal{C}_c^1(\Omega\times\Omega;\R^N)$. Then, if $|\phi| \le \omega$ a.e.~in $\Omega\times\Omega$, we have
    \begin{equation}
        \left|- \int_\Omega \int_\Omega \phi(x,y) \cdot \nabla_{\mathrm{NL}} u(x,y) \mathrm{d}x \mathrm{d}y\right| \le  \int_\Omega \int_\Omega |\phi(x,y)| |\nabla_{\mathrm{NL}} u(x,y)| \mathrm{d}x \mathrm{d}y \le \int_\Omega \int_\Omega \omega(x,y) |\nabla_{\mathrm{NL}} u(x,y)| \mathrm{d}x \mathrm{d}y,
    \end{equation}
    from which the desired inequality is proven by taking supremum over all such $\phi$. 

    Let us prove the other inequality. Let $U = \{ (x,y) \in \Omega\times\Omega : \omega \mbox{ continuous at } (x,y) \}$. Note that $|(\Omega\times\Omega) \backslash U| = 0$ because $\omega$ is continuous a.e. Consider $(\omega_m)_{m\in \N}$ given by
    \begin{equation}
        \omega_m(x,y) := \inf_{(x',y') \in U} \left( \omega(x',y') + m |(x,y) - (x',y')| \right)
    \end{equation}
    for $(x,y) \in \Omega \times\Omega$. For each $m \in \N$, it is standard to prove that $\omega_m$ is $m$-Lipschitz and $0 \le \omega_m$ and $\omega_m \le \omega_{m+1}$ in $\Omega\times\Omega$. Moreover, we also have that $\omega_m \le \omega$ in $U$ for each $m \in \N$. Let us prove that $\omega_m \to \omega$ a.e.~in $\Omega\times\Omega$. Let $(x,y) \in U$, and for $\epsilon>0$ choose $\delta>0$ so that, for $(x',y') \in B_\delta(x,y) \cap U$,
    \begin{equation}
        |\omega(x',y') - \omega(x,y)| < \epsilon \implies \omega(x',y') > \omega(x,y) - \epsilon.
    \end{equation}
    Choosing $m\in \N$ such that $m > \omega(x,y) / \delta$, then $\omega_m(x,y) \ge \min \{\omega(x,y) - \epsilon, m\delta \} = \omega(x,y) - \epsilon$. As $\omega_m \le \omega$ in $U$, this yields $|\omega(x,y) - \omega_m(x,y)|  \le \epsilon$, from which convergence at $(x,y)$ is proven. As $U$ is of full measure, this gives convergence a.e.~in $\Omega\times\Omega$.

    Now, fix $V \subset \subset \Omega$ open, $m \in \N$ and $\epsilon'>0$. We claim that we can find $\phi \in \mathcal{C}_c^1(\Omega\times\Omega;\R^N)$  such that $|\phi| \le \omega$ a.e.~in $\Omega\times\Omega$ and

    \begin{equation} \label{C1NLTVClaim}
        \int_\Omega\int_\Omega \phi \cdot \nabla_{\mathrm{NL}} u \; \mathrm{d}x \mathrm{d}y \ge \int_V \int_V \omega_m |\nabla_{\mathrm{NL}} u|\mathrm{d}x \mathrm{d}y - \epsilon'.
    \end{equation}
    We assume that $\nabla_{\mathrm{NL}} u$ is not equivalent to zero a.e.~in $V\times V$, as otherwise the claim is trivial. We start by considering
    \begin{equation}
        \psi = \frac{\nabla_{\mathrm{NL}} u}{|\nabla_{\mathrm{NL}} u|} \mathcal{X}_{\{\nabla_{\mathrm{NL}} u \neq 0\}}
    \end{equation}
    a.e.~in $\Omega\times\Omega$, so that $|\psi| \le 1$ and $\psi \cdot \nabla_{\mathrm{NL}} u = |\nabla_{\mathrm{NL}} u|$ a.e.~in $\Omega\times\Omega$. Then, clearly $\psi \in L^1(V\times V;\R^N)$, as $\psi$ is bounded on $V\times V$, a bounded set, and by density of $\mathcal{C}_c^1(V\times V;\R^N)$ in $L^1(V\times V;\R^N)$, there exists $(\psi_n)_{n\in\N} \subset \mathcal{C}_c^1(V\times V;\R^N)$ such that $\psi_n \to \psi$ in $L^1(V\times V;\R^N)$.
    Additionally, $(\psi_{n})_{n\in \N}$ can be chosen so that $|\psi_{n}| \le 1$ in $\Omega\times\Omega$ for every $n\in \N$ also (the standard mollification and multiplication by a smooth cut-off function already does the job).
    
    By $L^1$ convergence we know that, up to a subsequence, $\psi_n \to \psi$ a.e.~in $V\times V$. Hence, $|\psi - \psi_n| \omega_m |\nabla_{\mathrm{NL}} u| \to 0$ a.e.~in $V\times V$, and this sequence is dominated a.e.~by $2 \omega_m |\nabla_{\mathrm{NL}} u|$, which is integrable because $\omega_m$ is bounded in $V\times V$ and $\nabla_{\mathrm{NL}} u \in L^1_{\mathrm{loc}}(\Omega\times\Omega)$. Therefore, dominated convergence theorem yields that $|\psi - \psi_n| \omega_m |\nabla_{\mathrm{NL}} u| \to 0$ in $L^1(V\times V)$, and we can find $\psi_{\epsilon'} \in \mathcal{C}_c^1(V\times V;\R^N)$ such that

    \begin{equation} \label{C1NLTVAux1}
        \left|\int_V \int_V \omega_m (\psi - \psi_{\epsilon'}) \cdot \nabla_{\mathrm{NL}} u \; \mathrm{d}x \mathrm{d}y\right| \le \int_V\int_V |\omega_m| |\psi - \psi_{\epsilon'}| |\nabla_{\mathrm{NL}} u| \; \mathrm{d}x \mathrm{d}y < \frac{\epsilon'}{2}.
    \end{equation}
    Moreover, for $0<\delta_1 < \mathrm{dist}(V\times V, \partial(\Omega\times\Omega))$ and $\delta_2 > 0$ consider
    \begin{equation}
        \omega_{m,\delta_1,\delta_2}(x,y) := f_{\delta_2} (\omega_{m,\delta_1}(x,y)),
    \end{equation}
    for $(x,y) \in V\times V$ where
    \begin{equation}
        f_{\delta_2} (t) := \begin{cases}
            \frac{t^2}{t + \delta_2} & \mbox{ for } t\ge 0,\\
            0 & \mbox{ otherwise,}
        \end{cases}
    \end{equation}
    for $t \in \R$, and
    \begin{equation}
        \omega_{m,\delta_1}(x,y) := (J_{\delta_1} * \omega_m)(x,y) - m \delta_1,
    \end{equation}
    with $J_{\delta_1}$ being a standard mollifier in $\R^N \times \R^N$. Note that $\omega_{m,\delta_1} \in \mathcal{C}^\infty(V\times V)$, and since $\omega_{m}$ is $m$-Lipschitz we have that $\omega_{m} - 2m\delta_1 \le \omega_{m,\delta_1} \le \omega_m$. For $t\ge 0$, $f_{\delta_2}$ satisfies that
    \begin{equation}
        0 \le f_{\delta_2}(t) \le t, \qquad t - f_{\delta_2}(t) \le \delta_2.
    \end{equation}
    Hence, for every $(x,y) \in V\times V$, if $\omega_{m,\delta_1}(x,y) \ge 0$, then
    \begin{align}
        0 \le \omega_{m}(x,y) - \omega_{m,\delta_1}(x,y) &\le \omega_{m}(x,y) - f_{\delta_2} (\omega_{m,\delta_1}(x,y))\\
        &= (\omega_{m}(x,y) - \omega_{m,\delta_1}(x,y)) + (\omega_{m,\delta_1}(x,y) - f_{\delta_2} (\omega_{m,\delta_1}(x,y))) \le 2m\delta_1 + \delta_2,
    \end{align}
    and if $\omega_{m,\delta_1}(x,y) < 0$, $\omega_{m}(x,y) - f_{\delta_2} (\omega_{m,\delta_1}(x,y)) = \omega_{m}(x,y)$ and
    \begin{equation}
        0\le \omega_{m}(x,y) \le \omega_{m,\delta_1}(x,y) +2m\delta_1 \le 2m\delta_1  \le 2 m \delta_1 + \delta_2. 
    \end{equation}
    With these facts, we have that $\omega_{m,\delta_1,\delta_2} \in \mathcal{C}^1(V\times V)$ and
    \begin{equation}
        0 \le \omega_m - \omega_{m,\delta_1,\delta_2} \le 2m \delta_1 + \delta_2 \mbox{ in } V\times V \implies \| \omega_{m,\delta_1,\delta_2} - \omega_m\|_{L^\infty(V\times V)} \le 2m \delta_1 + \delta_2.
    \end{equation}
    Now, let $\omega_{m,\epsilon'} \in \mathcal{C}^1(V\times V)$ be associated with $\delta_1, \delta_2$ small enough so that
    \begin{equation}
        2m \delta_1 + \delta_2 < \frac{\epsilon'}{2\|\nabla_{\mathrm{NL}} u\|_{L^{1}(V\times V;\R^N)}}.
    \end{equation}
    Then,
    \begin{align} 
        \left|\int_V \int_V (\omega_{m,\epsilon'} - \omega_m) \psi_{\epsilon'} \cdot \nabla_{\mathrm{NL}} u \; \mathrm{d}x \mathrm{d}y\right| &\le \int_V \int_V |\omega_{m,\epsilon'} - \omega_m| |\psi_{\epsilon'}| |\nabla_{\mathrm{NL}} u| \mathrm{d}x \mathrm{d}y\\
        &\le \int_V \int_V |\omega_{m,\epsilon'} - \omega_m| |\nabla_{\mathrm{NL}} u| \mathrm{d}x \mathrm{d}y\\
        &\le \| \omega_{m,\epsilon'} - \omega_m\|_{L^\infty(V\times V)} \|\nabla_{\mathrm{NL}} u\|_{L^1(V\times V;\R^N)}< \frac{\epsilon'}{2}. \label{C1NLTVAux2}
    \end{align}
    With this, choose $\phi = \omega_{m,\epsilon'} \psi_{\epsilon'}$. As $\phi \in \mathcal{C}_c^1(V\times V;\R^N)$, we can extend by zero to a function $\phi \in \mathcal{C}_c^1(\Omega\times \Omega;\R^N)$ which additionally satisfies $|\phi| \le \omega_m \le \omega$ a.e.~in $\Omega\times\Omega$. Moreover, by the estimates \eqref{C1NLTVAux1} and \eqref{C1NLTVAux2},

    \begin{align}
        \int_\Omega\int_\Omega \phi \cdot \nabla_{\mathrm{NL}} u \; \mathrm{d}x \mathrm{d}y  &= \int_V \int_V \omega_{m,\epsilon'} \psi_{\epsilon'} \cdot \nabla_{\mathrm{NL}} u \; \mathrm{d}x \mathrm{d}y\\
        &= \int_V \int_V \omega_m \psi \cdot \nabla_{\mathrm{NL}} u \; \mathrm{d}x \mathrm{d}y\\
        &\qquad \qquad - \int_V \int_V (\omega_m - \omega_{m,\epsilon'}) \psi_{\epsilon'} \cdot \nabla_{\mathrm{NL}} u \; \mathrm{d}x \mathrm{d}y - \int_V \int_V \omega_m (\psi - \psi_{\epsilon'}) \cdot \nabla_{\mathrm{NL}} u \; \mathrm{d}x \mathrm{d}y\\
        &\ge \int_V \int_V \omega_m \psi \cdot \nabla_{\mathrm{NL}} u \; \mathrm{d}x \mathrm{d}y - \frac{\epsilon'}{2} - \frac{\epsilon'}{2}= \int_V \int_V \omega_m  |\nabla_{\mathrm{NL}} u| \mathrm{d}x \mathrm{d}y - \epsilon',
    \end{align}
    from which the claim is proven. By \eqref{C1NLTVAux0} and the claim, by choosing $\tilde \phi = - \phi$, where $\phi$ is the $\mathcal{C}_c^1(\Omega\times\Omega;\R^N)$ function that satisfies \eqref{C1NLTVClaim}, we deduce that

    \begin{equation}
        \mathrm{NLTV}_{\omega}(u) \ge \int_V \int_V \omega_m  |\nabla_{\mathrm{NL}} u| \mathrm{d}x \mathrm{d}y - \epsilon'
    \end{equation}
    for every $\epsilon'>0$, $m\in \N$ and $V \subset \subset \Omega$. Hence
    \begin{equation}
        \mathrm{NLTV}_{\omega}(u) \ge \int_V \int_V \omega_m  |\nabla_{\mathrm{NL}} u| \mathrm{d}x \mathrm{d}y
    \end{equation}
    for every $m\in \N$ and $V \subset \subset \Omega$. By monotone convergence theorem we may pass the $m \to \infty$ limit to obtain that

    \begin{equation}
        \mathrm{NLTV}_{\omega}(u) \ge \int_V \int_V \omega  |\nabla_{\mathrm{NL}} u| \mathrm{d}x \mathrm{d}y
    \end{equation}
    for every $V \subset \subset \Omega$. Cover $\Omega$ by $(V_n)_{n\in \N}$ such that $V_n \subset \subset \Omega$ and $V_n \subset V_{n+1}$ for each $n\in \N$, and monotone convergence theorem again yields that 
    \begin{equation}
        \mathrm{NLTV}_{\omega}(u) \ge \int_\Omega \int_\Omega \omega  |\nabla_{\mathrm{NL}} u| \mathrm{d}x \mathrm{d}y
    \end{equation}
    (independently of the finiteness of the right-hand side).
\end{proof}

\section{NLMeans weights in the continuous setting} \label{app:NLMeansWeights}

In this appendix, we formally formulate the NLMeans weights defined by \eqref{NLMeansWeightDefDiscrete} in the continuous setting, and show that it is consistent with the theoretical setting described in Section \ref{NLTVPriorSec}.

Let $\Omega \subset \R^N$ be nonempty, open, bounded and convex, and $u : \Omega \to \R$ denote the initial image from which the weight is computed. We define $f: \Omega\times \Omega \to [0,+\infty)$ as

\begin{equation} \label{NLMeansDefCont}
    f(x,y) = \frac{1}{\Gamma(x)} \exp \left( -\frac{D_{u}(x,y)}{h^2} \right),
\end{equation}
where $h>0$,
\begin{equation} \label{GammaDef}
    \Gamma(x) = \int_\Omega \exp \left( -\frac{D_{u}(x,y)}{h^2} \right) \mathrm{d}y
\end{equation}
and

\begin{equation} \label{D_uDef}
    D_u(x,y) = \int_{B_\nu(0)} |\tilde u(x+z) - \tilde u (y+z)|^2 \mathrm{d}z,
\end{equation}
with $\nu>0$ and $\tilde u$ denoting the extension of $u$ to $\Omega + B_\nu(0)$ to treat out of boundary values, which we assume to be in $L^2(\Omega + B_\nu(0))$. 

\begin{prop} \label{NLMeansContProp}
    The function $f$ given by \eqref{NLMeansDefCont} is continuous and strictly positive in $\Omega\times\Omega$ and continuously extendable to $\overline{\Omega\times\Omega}$ and strictly positive there.
\end{prop}

\begin{proof}
    Let $\hat u \in L^2(\R^N)$ be the extension of $\tilde u$ by zero to the whole $\R^N$. Hence, both $\Gamma$ \eqref{GammaDef} and $D_u$ \eqref{D_uDef} admit an extension to $\R^N$ and $\R^N\times \R^N$, respectively.
    
    Let us see that $D_u$ (its extension) is continuous in $\R^N \times \R^N$. Note that we may write

    \begin{equation}
        D_u(x,y) = F(x) + F(y) - 2 G(x,y),
    \end{equation}
    where $F: \R^N \to \R$ is defined by
    \begin{equation}
        F(x) := \|\hat u(x + \cdot)\|_{L^2(B_\nu(0))}^2 = \int_{B_\nu(0)} |\hat  u(x+z)|^2\mathrm{d}z
    \end{equation}
    and $G: \R^N \times \R^N\to \R$ is defined by
    \begin{equation}
        G(x,y) := \langle \hat u(x+\cdot), \hat u(y + \cdot) \rangle_{L^2(B_\nu(0))} = \int_{B_\nu(0)} \hat  u(x+z) \hat u (y+z) \mathrm{d}z.
    \end{equation}

    To see that both $F$ and $G$ are continuous, we use the continuity of translations with respect to $L^p$ norm in $\R^N$ \cite[Proposition 8.5]{folland1999real}, which states that the mapping $\Phi_{\hat u}: \R^N \to L^2(\R^N)$ given by
    \begin{equation}
        \Phi_{\hat u} (a) := \hat u(a + \cdot) 
    \end{equation}
    is continuous. In particular, for a fixed $x_0 \in \R^N$ we have that $\hat u(x+\cdot) \to \hat u(x_0+\cdot)$ in $L^2(B_{\nu}(0))$ as $x \to x_0$, from which
    \begin{equation}
        \lim_{x\to x_0} \|\hat u(x+\cdot)\|_{L^2(B_\nu(0))} = \|\hat u(x_0+\cdot)\|_{L^2(B_\nu(0))} \implies \lim_{x\to x_0} \|\hat u(x+\cdot)\|_{L^2(B_\nu(0))}^2 = \|\hat u(x_0+\cdot)\|_{L^2(B_\nu(0))}^2,
    \end{equation}
    so $F$ is continuous in $\R^N$. Moreover, for $x_0, y_0 \in \R^N$, we analogously deduce that
    \begin{equation}
        \lim_{x \to x_0} \|\hat u(x_0+\cdot) - \hat u(x+\cdot)\|_{L^2(B_\nu(0))} = 0, \qquad \lim_{y \to y_0} \| \hat u(y_0+\cdot) - \hat u(y+\cdot)\|_{L^2(B_\nu(0))} = 0, 
    \end{equation}
    and
    \begin{align}
        |&G(x, y) - G(x_0,y_0)|\\
        &\le |\langle \hat u(x+\cdot) - \hat u(x_0+\cdot), \hat u (y + \cdot) \rangle_{L^2(B_\nu(0))}| + |\langle \hat u (x_0+ \cdot), \hat u(y+\cdot) - \hat u(y_0 +\cdot)\rangle_{L^2(B_\nu(0))}|\\
        &\le \|\hat u(x_0+\cdot) - \hat u(x+\cdot)\|_{L^2(B_\nu(0))} \|\hat u (y+\cdot)\|_{L^2(B_\nu(0))} + \|\hat u (x_0+\cdot)\|_{L^2(B_\nu(0))} \| \hat u(y_0+\cdot) - \hat u(y+\cdot)\|_{L^2(B_\nu(0))},
    \end{align}
    which tends to zero as $x \to x_0$ and $y \to y_0$ by continuity of $F$ and the fact that $\hat u \in L^2(\R^N)$. From this we deduce the continuity of $D_u$ in $\R^N\times\R^N$.
    
    Now, note that $\Gamma > 0$ in $\R^N$. Indeed, the integrand in \eqref{GammaDef} is strictly positive because $D_u$ is finite due to the fact that $\hat u \in L^2(\R^N)$. Since the expression appearing in the exponential is always negative, the integrand is bounded by $1$ (which is integrable in $\Omega$ by the boundedness of the domain). Hence, the continuity of $D_u$ and the dominated convergence theorem yield that $\Gamma$ is also continuous in $\R^N$.
    
    We conclude that $f$ has a continuous extension to $\R^N\times \R^N$ in which it is strictly positive, from which it is continuous and strictly positive in $\Omega\times\Omega$ and continuously extendable to $\overline{\Omega\times\Omega}$ and strictly positive there.
\end{proof}

With Proposition \ref{NLMeansContProp}, we consider $\rho >0$ and
\begin{equation}
    \omega := \mathcal{X}_{\{|x-y| < \rho\}} f,
\end{equation}
and reproduce the arguments described in Remark \ref{SearchWindowSatisfiesH0} to obtain that $\omega$ is an admissible weight function that satisfies \eqref{H0}.

\end{appendix}

\end{document}